\documentclass[paper,revised]{geophysics}		

\usepackage{amssymb,amsmath,amsthm,amsthm}
\allowdisplaybreaks 
\usepackage{amsbsy}
\usepackage{mathtools}

\usepackage[english]{babel}
\usepackage{csquotes}
\usepackage{accents}
\usepackage{dsfont}
\usepackage{bm}
\usepackage{booktabs}
\usepackage{etoolbox}
\apptocmd{\sloppy}{\hbadness 10000\relax}{}{}

\usepackage[shortlabels]{enumitem}

\usepackage{graphicx}
\usepackage{array,multirow,hhline}
\usepackage{subfig}     
\usepackage[dvipsnames,x11names]{xcolor}
\usepackage{graphicx}
\usepackage{amsmath}
\usepackage{tikz}
\usepackage{amsfonts}
\usetikzlibrary{arrows.meta, shapes.geometric, positioning, decorations.pathreplacing, shadows, fit, calc}
\usetikzlibrary{positioning}

\definecolor{azure}{rgb}{0.0, 0.5, 1.0}
\definecolor{offwhite}{RGB}{248, 248, 248}
\definecolor{DarkBlueNavy}{RGB}{0, 0, 128}
\definecolor{nicegreen}{RGB}{76, 175, 80}
\definecolor{LightGrey}{RGB}{211, 211, 211} 
\definecolor{niceblue}{RGB}{0, 102, 204} 
\definecolor{nicegreen}{RGB}{76, 175, 80} 
\definecolor{LightGrey}{RGB}{211, 211, 211} 

\definecolor{darkred}{rgb}{0.0, 0.0, 0.0} 

\definecolor{myblue}{RGB}{0,85,164}
\definecolor{Tomato}{rgb}{0.8, 0.0, 0.0}
\definecolor{Goldenrod}{RGB}{218,165,32}
\definecolor{mywhite}{RGB}{255,255,255}

\newcommand{\red}{\textcolor{darkred}}

\usepackage{pgfplots}
\pgfplotsset{compat=1.18}

\newcommand{\colorbar}[4]{
\begin{tikzpicture}
    \begin{axis}[
        axis on top,
        height=#4,
        width=0.5cm,
        scale only axis,
        enlargelimits=false,
        xmin=0,
        xmax=1,
        ymin=#2,
        ymax=#3,
        xtick=\empty,
        yticklabel pos=right,
        y tick label style = {font=\footnotesize},
        try min ticks=5
        ]
      \addplot[] graphics[xmin=0,ymin=#2,xmax=1,ymax=#3] {Figures/colormap_#1.png};
    \end{axis}
  \end{tikzpicture}
}

\newcommand{\horizontalColorbar}[4]{
	\begin{tikzpicture}
		\begin{axis}[
			axis on top,
			height=0.5cm,
			width=#4,
			scale only axis,
			enlargelimits=false,
			xmin=#2,
			xmax=#3,
			ymin=0,
			ymax=1,
			ytick=\empty,
			x tick label style = {font=\footnotesize},
			try min ticks=5
			]
			\addplot[] graphics[xmin=#2,ymin=0,xmax=#3,ymax=1] {Figures/colormap_#1.png}; 
		\end{axis}
	\end{tikzpicture}
}

\newcommand{\marmousiCB}[1]{\colorbar{jet}{-1.8654}{3.6714}{#1}}
\newcommand{\marmousiCBH}[1]{\horizontalColorbar{jet}{-1.8654}{3.6714}{#1}}
\newcommand{\seamCB}[1]{\colorbar{jet}{-2.2221}{1.0064}{#1}}
\newcommand{\seamCBH}[1]{\horizontalColorbar{jet}{-2.2221}{1.0064}{#1}}

\usepackage[hypertexnames=false, hidelinks, colorlinks=true, linkcolor=niceblue, citecolor=azure, urlcolor=nicegreen]{hyperref}
\newtheorem{theorem}{Theorem}

\newtheorem{definition}[theorem]{Definition}

\newcommand{\R}{\mathbb{R}}

\newcommand{\argmin}{\operatorname{argmin}}

\newcommand{\bx}{{\boldsymbol{x}}}
\newcommand{\by}{{\boldsymbol{y}}}
\newcommand{\bz}{{\boldsymbol{z}}}

\newcommand{\Loss}{L}

\usepackage{algorithm}
\usepackage{algpseudocode}

\newcommand{\Hm}[1]{\leavevmode{\marginpar{\tiny%
			$\hbox to 0mm{\hspace*{-0.5mm}$\leftarrow$\hss}%
			\vcenter{\vrule depth 0.1mm height 0.1mm width \the\marginparwidth}%
			\hbox to 0mm{\hss$\rightarrow$\hspace*{-0.5mm}}$\\\relax\raggedright
			#1}}}

\title{Improved impedance inversion \\by the iterated graph Laplacian}

\author{Davide Bianchi\footnote{School of Mathematics (Zhuhai),
Sun Yat-sen University, Zhuhai, 518055, China.\newline Email: bianchid@mail.sysu.edu.cn}, Florian Bo\ss mann\footnote{School of Mathematics, Harbin Institute of Technology, Harbin, 150001, China.\newline Email: f.bossmann@hit.edu.cn \newline Email: wenlong.wang@hit.edu.cn}, Wenlong Wang\footnotemark[2] and Mingming Liu\footnotemark[2]}

\date{}
\begin{document}
\maketitle

\begin{abstract}
\red{ 
We introduce a data-adaptive inversion method that integrates classical or deep learning-based approaches with iterative graph Laplacian regularization, specifically targeting acoustic impedance inversion — a critical task in seismic exploration. Our method initiates from an impedance estimate derived using either traditional inversion techniques or neural network-based methods. This initial estimate guides the construction of a graph Laplacian operator, effectively capturing structural characteristics of the impedance profile. Utilizing a Tikhonov-inspired variational framework with this graph-informed prior, our approach iteratively updates and refines the impedance estimate while continuously recalibrating the graph Laplacian. This iterative refinement shows rapid convergence, increased accuracy, and enhanced robustness to noise compared to initial reconstructions alone. Extensive validation performed on synthetic and real seismic datasets across varying noise levels confirms the effectiveness of our method. Performance evaluations include four initial inversion methods: two classical techniques and two neural networks — previously established in the literature.} 
\end{abstract}


\section{Introduction}

\red{Impedance inversion has become a standard technique used to characterize subsurface reservoirs in many applications \cite{sams2013practical,ray2016building}. Here,}
subsurface properties and structures have to be reconstructed from seismic profiles. Mathematically, a subsurface impedance $\bx\in X$ has to be recovered from a possibly noisy seismic profile $\by\in Y$, i.e., we have to find $\bx$ such that
\begin{equation}\label{eq:seismic}
    \by^\delta = F(\bx)+\boldsymbol{\eta}, \qquad \|\boldsymbol{\eta}\|\leq \delta,
\end{equation}
where $F:X\rightarrow Y$ is the forward operator. \red{The inverse problem of equation \ref{eq:seismic} is ill-posed in many} cases \cite{wu2021deep} which means that even a small amount of noise $0<\delta\ll1$ can lead to extreme approximation errors. This ill-posedness originates from a non-linear operator $F$, heterogeneity in the profile $\bx$, or an insufficient amount of data available.

A classical approach to overcome the ill-posedness is to regularize the inverse problem \cite{wang2010seismic,liu2015impedance}. Many of these methods fall under the broader category of Tikhonov-like variational methods which assumes the form
\begin{equation}\label{genTik}
	\bx_\alpha \coloneqq \underset{\bx \in X}{\argmin}\left\{ \mathcal{D}(F(\bx); \by^\delta) + \alpha \mathcal{R}(\bx) \right\},
\end{equation}
where $\mathcal{D}(\cdot; \cdot)$ is a pseudo-distance that quantifies the fidelity of the reconstruction,  $\mathcal{R}(\cdot)$  serves as regularization term, and $\alpha >0$ balances the trade-off between data fidelity and the regularization effect. A standard choice is
\begin{equation*}
    \mathcal{D}(F(\bx); \by^\delta)=\frac{1}{2}\|F(\bx)-\by^\delta\|_2^2, \qquad \mathcal{R}(\bx)=\|\bx\|_1.
\end{equation*}
As references see \cite{engl1996regularization,scherzer2009variational}.

The regularization operator $\mathcal{R}$ is of crucial importance. If we have access to a priori information on the ground-truth solution, we can tailor a specific $\mathcal{R}$ to incorporate such information and guide the overall regularization towards a narrower subset of approximate solutions which present the features we aim to recover.

Another typical example in image processing is given by $\mathcal{R}(\bx) = \|L\bx\|_q^q$, where $L$ is a linear differential operator, such as the first or second derivative along each axis, and $\|\cdot\|_q$ is the Euclidean $q$-norm, $q\geq 1$, see  \cite{hansen2006deblurring,jain1989fundamentals,ng1999fast}. \red{Differential operators are particularly useful because they detect intensity jumps (edges), widely acknowledged as essential visual cues for perception and image interpretation. See the discussions in \cite[Section 3.2.4]{chan2005image}, \cite[Chapter 2]{marr2010vision}, and \cite{canny1986computational,mumford1989optimal}.}

Digital images, which consist of pixels arranged on a 2D grid, naturally possess a graph structure. For this reason, in recent years, graph-based differential operators~$\Delta$ have been introduced to replace standard Euclidean differential operators~$L$ in various image processing tasks such as denoising  \cite{gilboa2009nonlocal}, image deblurring \cite{bianchi2021graph,bianchi2021graph_approximation,buccini2021graph,aleotti2023fractional,zhang2010bregmanized}, Computed Tomography \cite{lou2010image}, and other applications \cite{Peyre2008nonlocal,arias2009variational,gilboa2007nonlocal}. Graph operators showed a general good performance, due to the fact that they can  more effectively model the complex structures and textures present in images. Unlike traditional Euclidean differential operators that primarily consider the spatial proximity of pixels, graph-based operators can take into account the intensity similarity between pixels as well, allowing for a deeper information extrapolation from image data.

Besides these classical methods, machine learning and deep learning approaches have gained interest over recent years. They have successfully been applied in seismic interpolation \cite{wang2019deep}, full waveform inversion \cite{zhang2022regularized}, impedance inversion \cite{das2019convolutional}, and seismic interpretation \cite{wrona20213d}. See also \cite{yu2021deep,khosro2024machine} for an overview. In many of these cases neural networks generate the best outcome by far. A great advantage of such techniques is that the non-linearity in \red{equation} \ref{eq:seismic} is directly learned from the data. This can also be used to construct a forward operator without requiring a detailed and complicated physical model \cite{Alfarraj2019network1}. However, training a neural network requires a suitable amount and quality of data. A poorly trained network can result in strong artifacts and errors in the reconstruction, especially when dealing with ill-posed inverse problems \cite{antun2020instabilities,colbrook2022difficulty}. This can happen whenever there is insufficient data, biased data, noise, inaccurate labeling, or a shift in the dataset between training and deployment \cite{huot2018jump,zhang2021robust,karimpouli2020physics}. This makes data acquisition a complicated and expensive process.

\begin{figure}[!ht]
    \centering
    \includegraphics [width=0.9\textwidth]{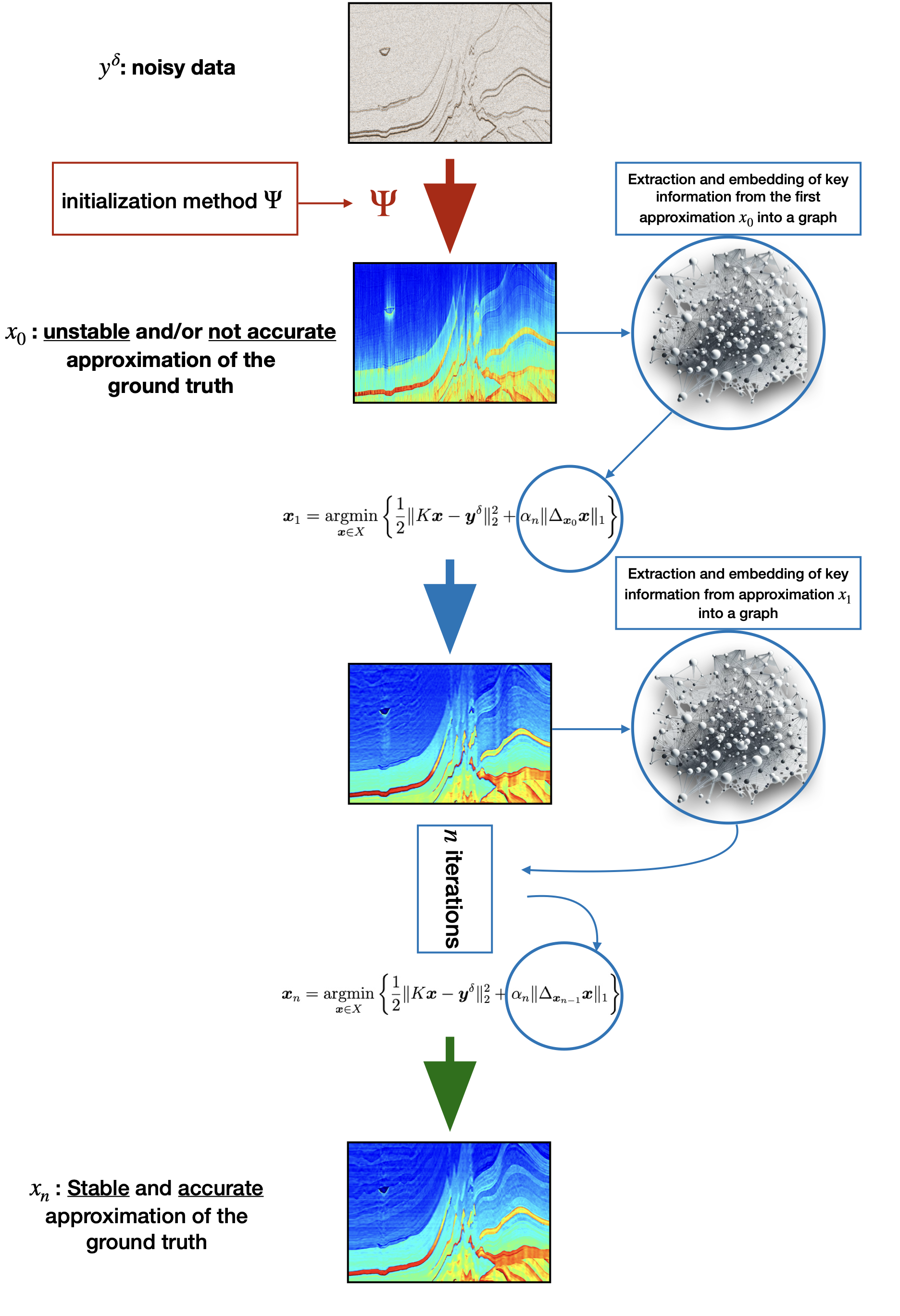}
    \caption{Pipeline of the iterative method \ref{eq:iterative_graphLaNet}.}\label{fig:visual_abstract}
\end{figure}

To overcome said problem hybrid methods that combine neural networks with classical approaches have been studied lately \cite{arridge2019solving,li2020nett}. In \cite{bianchi2023data,bianchi2023graph}, a new `two-step' method called \texttt{graphLa}$\Psi$ was proposed. The idea is to get a first approximate reconstruction of the ground-truth solution via a generic reconstructor operator $\Psi \colon Y \to X$, applied to the observed data $\by^\delta$, and then solve equation~\ref{genTik} replacing $\mathcal{R}(\bx)$ with $\|\Delta_\Psi \bx\|_1$, where $\Delta_\Psi$ is a graph Laplacian built from $\Psi(\by^\delta)$. The first approximate reconstruction $\Psi(\by^\delta)$ helps identifying the pixel connections in the ground-truth signal by both their spatial proximity and intensity similarity. This approach \red{has been empirically seen to} greatly improve the quality of the approximate solution $\Psi(\by^\delta)$ for any given initial reconstructor $\Psi$.

\red{The initial reconstructor can be quite general and may vary depending on the problem and the user's preferences. From a formal mathematical perspective, it only needs to satisfy very mild assumptions. For instance, any regularization method (as per \cite[Definition 3.1]{engl1996regularization}) or even any locally Lipschitz operator can serve as the initial reconstructor. See \cite[Section 3]{bianchi2023data} and the discussion presented therein. In particular, we highlight that any modern neural network architecture satisfies the assumptions that make the \texttt{graphLa}$\Psi$ method regularizing.}

In this work we apply the \texttt{graphLa}$\Psi$ method iteratively to the seismic impedance inversion equation \ref{eq:seismic}, using as initial reconstructor  \red{two classical methods (Sparse Spike Inversion \cite{velis2008stochastic} and Split Bregman \cite{gholami2016fast}) and} two Deep Neural Networks (DNNs) specifically trained on seismic impedance datasets \cite{Alfarraj2019network1,Mingming2024network2}. 

\red{A significant issue with equation~\ref{eq:seismic} is the non-linearity of the involved operator $F$. While  both DNNs can learn the non-linear forward operator \ref{eq:seismic} directly from the given data, the classical methods as well as \texttt{graphLa}$\Psi$ require a linear forward model. Hence, we use a linear operator $K$ as approximation,i.e., $K\bx\approx F(\bx)$. This is a commonly used approach and yields good results under mild assumptions \cite{assis2019colored}.}

\newpage
Our iterative method then reads,

\begin{equation}\tag{\texttt{it-graphLa}$\Psi$}\label{eq:iterative_graphLaNet}
\begin{cases}
\bx_n = \underset{\bx \in X}{\argmin}\left\{ \frac{1}{2}\|K\bx - \by^\delta \|_2^2 + \alpha_n\|\Delta_{\bx_{n-1}} \bx\|_1 \right\} & \mbox{for } n\geq 1,\\
\bx^\delta_0 \coloneqq \Psi(\by^\delta),
\end{cases}
\end{equation}
where $\Psi:Y\rightarrow X$ is \red{the chosen initialization method}. The idea is that at each step~$n$, having at our disposal a better approximation $\bx_{n-1}$ of the ground-truth, we can update the weights of the graph Laplacian $\Delta$ taking into account this new information. The updated weights bring a sharper insight on the pixel connections in the ground-truth signal, driving then the overall reconstruction into a closer neighborhood of the ground truth. A visual abstract of our method in presented in Figure~\ref{fig:visual_abstract}.

The manuscript is organized as follows. In the first next section  we give the basic notation and theory required for this work. We also introduce the standard \texttt{graphLa}$\Psi$ method and shortly summarize its theory. Afterwards, we present our iterated \texttt{graphLa}$\Psi$ algorithm (\texttt{it-graphLa}$\Psi$). \red{The following numerical section discusses different initialization methods and presents results for model as well as field data.} A discussion of the results and outlook on future problems can be found in the conclusion section at the end of the paper.

\section{The graph Laplacian and the \texttt{graphLa}$\Psi$ method}\label{sec:graphLaplacian}
We first introduce some preliminaries on graph theory and theoretical results about the \texttt{graphLa}$\Psi$ regularization method.

\subsection{Images and graphs}
For a modern introduction to graph theory we invite the interested reader to consult \cite{keller2021graphs}. For simplicity, in this work we consider as a graph any pair $G=(P,w)$ where $P$ is a finite set, called \emph{node set}, and $w \colon P\times P \to [0,\infty)$ is a symmetric function, called \emph{edge-weight function}.  The set of the edges of a graph is given by $E\coloneqq \{ (p,q) \in P\times P \mid w(p,q)\neq 0 \}$. Two nodes $p$ and $q$ are connected if $(p,q)\in E$, and we write $p\sim q$. The intensity of the connection is given by $w(p,q)$.

For any function $\bx \colon P \to \R$ we can compute the graph Laplacian $\Delta \bx \colon P \to \R$, defined by the action
\begin{equation}\label{graph_laplacian}
	\Delta \bx (p) \coloneqq \sum_{q\sim p}w(p,q)(\bx(p) - \bx(q)).
\end{equation}

\subsubsection{Induced graph}
Observe that to define a graph we only need a node set $P$ and an edge-weight function~$w$. Let us see now how to generate a graph from an image. 

Since  any image is made by the union of several pixels $p \in P$ disposed on a grid, it is then natural to identify the pixels as ordered pairs $p=(i_p,j_p)$ with $i_p=1,\ldots, n$ and $j_p=1,\ldots, m$, and where $n$ and $m$ indicates the total number of pixels  along the horizontal and vertical axis, respectively. So, we can set the node set as 
$$
P = \{ p \mid p=(i_p, j_p), \; i_p=1,\ldots,n, \, j_p=1,\ldots,m \}.
$$

For the sake of simplicity, consider now a gray-scale image, which is given by the light intensities  of its pixels. That is, a gray-scale image can be represented by a function 
$$
\bx \colon P \to [0,1],
$$
where $0$ means black and $1$ means white. 
A very popular choice to make the connection $w$ of two pixels depending on both their spatial proximity and light intensity is given by
\begin{equation}\label{eq:edge-weight-function:applications}
	w_{\bx}(p,q) = \mathds{1}_{(0,R]}(\operatorname{dist}(p,q))g_\bx(p,q),
\end{equation}
where $\operatorname{dist}(\cdot, \cdot)$ is a (pseudo) distance on $P$ and 
\begin{equation}\label{exp1}
	g_\bx(p,q)\coloneqq \mathrm{e}^{-\frac{|\bx(p)-\bx(q)|^2}{\sigma}}.
\end{equation}
The function $\mathds{1}_{(0,R]}$ is the indicator function of the interval $(0,R]$ and  $R>0$ is a parameter of control which tells the maximum proximity distance allowed for two pixels to be neighbors. If the distance between pixels $p$ and $q$ satisfies $0 < \operatorname{dist}(p,q) \leq R$, then $p$ and $q$ are connected by an edge with magnitude $g_\bx(p,q)$, i.e., $w_{\bx}(p,q) = g_\bx(p,q)$. The second control parameter, $\sigma>0$, determines how sharply the edge weights vary with the difference in pixel light intensities. Common choices for the distance function $\operatorname{dist}(\cdot,\cdot)$ are 
\begin{equation*}
	\begin{aligned}
		&\operatorname{dist}(p,q)	= \operatorname{dist}_1(p,q) \coloneqq |i_p - i_q| + |j_p - j_q|,\\
		&\operatorname{dist}(p,q)	= \operatorname{dist}_\infty(p,q) \coloneqq \max\{  |i_p - i_q|; |j_p - j_q|\}.
	\end{aligned}
\end{equation*}
See Figure \ref{fig:example_image2graph_1} for a simple example of building a graph from an image.

\begin{figure}[tb!]
	\begin{minipage}{0.25\textwidth}
		\centering
		\begin{tikzpicture}[scale=0.4]

    \foreach \x in {0,...,6} {
        \foreach \y in {0,...,6} {
            \fill[myblue] (\x,\y) rectangle ++(1,1);
        }
    }
    
    \foreach \x/\y in {2/2, 3/2, 4/2, 2/3, 3/3, 4/3, 2/4, 3/4, 4/4} {
        \fill[Tomato] (\x,\y) rectangle ++(1,1);
    }

    \foreach \x/\y in {2/1, 3/1, 4/1} {
        \fill[Goldenrod] (\x,\y) rectangle ++(1,1);
    }
    \foreach \x/\y in {2/5, 3/5, 4/5} {
        \fill[Goldenrod] (\x,\y) rectangle ++(1,1);
    }
    \foreach \x/\y in {1/2, 1/3, 1/4, 5/2, 5/3, 5/4} {
        \fill[Goldenrod] (\x,\y) rectangle ++(1,1);
    }

    \foreach \x in {0,...,6} {
        \foreach \y in {0,...,6} {
            \draw (\x,\y) rectangle ++(1,1);
        }
    }
\end{tikzpicture}
	\end{minipage}%
	\begin{minipage}{0.25\textwidth}
		\centering
		\begin{tikzpicture}[scale=0.4]

    \foreach \x in {0,...,6} {
        \foreach \y in {0,...,6} {
            \fill[myblue] (\x,\y) rectangle ++(1,1);
        }
    }
    
    \foreach \x/\y in {2/2, 3/2, 4/2, 2/3, 3/3, 4/3, 2/4, 3/4, 4/4} {
        \fill[Tomato] (\x,\y) rectangle ++(1,1);
    }

    \foreach \x/\y in {2/1, 3/1, 4/1} {
        \fill[Goldenrod] (\x,\y) rectangle ++(1,1);
    }
    \foreach \x/\y in {2/5, 3/5, 4/5} {
        \fill[Goldenrod] (\x,\y) rectangle ++(1,1);
    }
    \foreach \x/\y in {1/2, 1/3, 1/4, 5/2, 5/3, 5/4} {
        \fill[Goldenrod] (\x,\y) rectangle ++(1,1);
    }

    \foreach \x in {0,...,6} {
        \foreach \y in {0,...,6} {
            \draw (\x,\y) rectangle ++(1,1);
        }
    }

    \foreach \x in {0,...,6} {
        \foreach \y in {0,...,6} {
            \draw (\x,\y) rectangle ++(1,1);
            
            \fill[black] (\x+0.5,\y+0.5) circle (0.12);
        }
    }

\end{tikzpicture} 
	\end{minipage}%
	\begin{minipage}{0.25\textwidth}
		\centering
		 \begin{tikzpicture}[scale=0.4]
    \foreach \x in {0,...,6} {
        \foreach \y in {0,...,6} {
            \fill[myblue] (\x,\y) rectangle ++(1,1);
        }
    }
    
    \foreach \x/\y in {2/2, 3/2, 4/2, 2/3, 3/3, 4/3, 2/4, 3/4, 4/4} {
        \fill[Tomato] (\x,\y) rectangle ++(1,1);
    }

    \foreach \x/\y in {2/1, 3/1, 4/1} {
        \fill[Goldenrod] (\x,\y) rectangle ++(1,1);
    }
    \foreach \x/\y in {2/5, 3/5, 4/5} {
        \fill[Goldenrod] (\x,\y) rectangle ++(1,1);
    }
    \foreach \x/\y in {1/2, 1/3, 1/4, 5/2, 5/3, 5/4} {
        \fill[Goldenrod] (\x,\y) rectangle ++(1,1);
    }

    \foreach \x in {0,...,6} {
        \foreach \y in {0,...,5} {
            \draw[line width=1mm, mywhite] (\x+0.5,\y+0.5) -- (\x+0.5,\y+1.5);
        }
    }
    \foreach \y in {0,...,6} {
        \foreach \x in {0,...,5} {
            \draw[line width=1mm, mywhite] (\x+0.5,\y+0.5) -- (\x+1.5,\y+0.5);
        }
    }

    \foreach \x in {0,...,6} {
        \foreach \y in {0,...,6} {
            \draw (\x,\y) rectangle ++(1,1);
        }
    }

    \foreach \x in {0,...,6} {
        \foreach \y in {0,...,6} {
            \fill[black] (\x+0.5,\y+0.5) circle (0.12);
        }
    }
\end{tikzpicture} 
	\end{minipage}%
	\begin{minipage}{0.25\textwidth}
		\centering
		\begin{tikzpicture}[scale=0.4]

  \newcommand{\getcolor}[2]{%
    \ifnum#1<0 myblue\else
    \ifnum#1>6 myblue\else
    \ifnum#2<0 myblue\else
    \ifnum#2>6 myblue\else
      \ifnum#1>1
        \ifnum#1<5
          \ifnum#2>1
            \ifnum#2<5
              Tomato
            \else
              Goldenrod
            \fi
          \else
            Goldenrod
          \fi
        \else
          \ifnum#1=5
            \ifnum#2>1
              \ifnum#2<5
                Goldenrod
              \else
                myblue
              \fi
            \else
              myblue
            \fi
          \else
            myblue
          \fi
        \fi
      \else
        \ifnum#1=1
          \ifnum#2>1
            \ifnum#2<5
              Goldenrod
            \else
              myblue
            \fi
          \else
            myblue
          \fi
        \else
          myblue
        \fi
      \fi
    \fi\fi\fi\fi
  }

  \foreach \x in {0,...,6} {
    \foreach \y in {0,...,6} {
      \fill[myblue] (\x,\y) rectangle ++(1,1);
    }
  }

  \foreach \x/\y in {2/2, 3/2, 4/2, 2/3, 3/3, 4/3, 2/4, 3/4, 4/4} {
    \fill[Tomato] (\x,\y) rectangle ++(1,1);
  }

  \foreach \x/\y in {2/1, 3/1, 4/1} {
    \fill[Goldenrod] (\x,\y) rectangle ++(1,1);
  }
  \foreach \x/\y in {2/5, 3/5, 4/5} {
    \fill[Goldenrod] (\x,\y) rectangle ++(1,1);
  }
  \foreach \x/\y in {1/2, 1/3, 1/4, 5/2, 5/3, 5/4} {
    \fill[Goldenrod] (\x,\y) rectangle ++(1,1);
  }


  \foreach \x in {0,...,6} {
    \foreach \y in {0,...,5} {
      \draw[line width=0.05mm, mywhite] (\x+0.5,\y+0.5) -- (\x+0.5,\y+1.5);
    }
  }

  \foreach \y in {0,...,6} {
    \foreach \x in {0,...,5} {
      \draw[line width=0.1mm, mywhite] (\x+0.5,\y+0.5) -- (\x+1.5,\y+0.5);
    }
  }


  \foreach \y in {0,...,6} {
    \foreach \x in {0,...,5} {
      \edef\currentcolor{\getcolor{\x}{\y}}%
      \edef\nextcolor{\getcolor{\the\numexpr\x+1\relax}{\y}}%
      \ifx\currentcolor\nextcolor
        \draw[line width=1mm, mywhite] (\x+0.5,\y+0.5) -- (\x+1.5,\y+0.5);
      \fi
    }
  }

\draw[line width=1mm, mywhite] (1.5, 0.5) -- (2.5, 0.5);
\draw[line width=1mm, mywhite] (4.5, 0.5) -- (5.5, 0.5);
\draw[line width=1mm, mywhite] (1.5, 6.5) -- (2.5, 6.5);
\draw[line width=1mm, mywhite] (4.5, 6.5) -- (5.5, 6.5);

\draw[line width=1mm, mywhite] (0.5, 0.5) -- (0.5, 6.5);
\draw[line width=1mm, mywhite] (6.5, 0.5) -- (6.5, 6.5);

\draw[line width=1mm, mywhite] (1.5, 2.5) -- (1.5, 4.5);
\draw[line width=1mm, mywhite] (5.5, 2.5) -- (5.5, 4.5);

\draw[line width=1mm, mywhite] (2.5, 2.5) -- (2.5, 4.5);
\draw[line width=1mm, mywhite] (4.5, 2.5) -- (4.5, 4.5);

\draw[line width=1mm, mywhite] (1.5, 0.5) -- (1.5, 1.5);
\draw[line width=1mm, mywhite] (1.5, 5.5) -- (1.5, 6.5);

\draw[line width=1mm, mywhite] (5.5, 0.5) -- (5.5, 1.5);
\draw[line width=1mm, mywhite] (5.5, 5.5) -- (5.5, 6.5);

\draw[line width=1mm, mywhite] (3.5, 2.5) -- (3.5, 4.5);

  \foreach \x in {0,...,6} {
    \foreach \y in {0,...,6} {
      \draw (\x,\y) rectangle ++(1,1);

      \fill[black] (\x+0.5,\y+0.5) circle (0.12);
    }
  }

\end{tikzpicture} 
	\end{minipage}\caption{Simple illustration that demonstrates a step-by-step process for converting an image $\bx$ into a graph. Starting from the left, a $7\times7$ pixel image—composed of red, orange, and blue squares—is shown, where each pixel’s color intensity is determined by evaluating the function $\bx$. Moving right, each pixel is mapped to a graph node (depicted as a black circle) with its position corresponding to an ordered pair in $\mathbb{Z}^2$ based on its location in the grid. Next, nodes are connected by an edge whenever their $\ell^1$-distance is one, with each connection initially assigned a weight of one. Finally, the edge weight is adjusted by the function $g_\bx(p,q)$, which is visually represented by the edge thickness: thicker lines indicate that adjacent pixels have very similar intensities, while thinner lines indicate a larger difference.}\label{fig:example_image2graph_1}
\end{figure}

\begin{figure}[h!tb]
	\centering
	\includegraphics [width=0.8\textwidth]{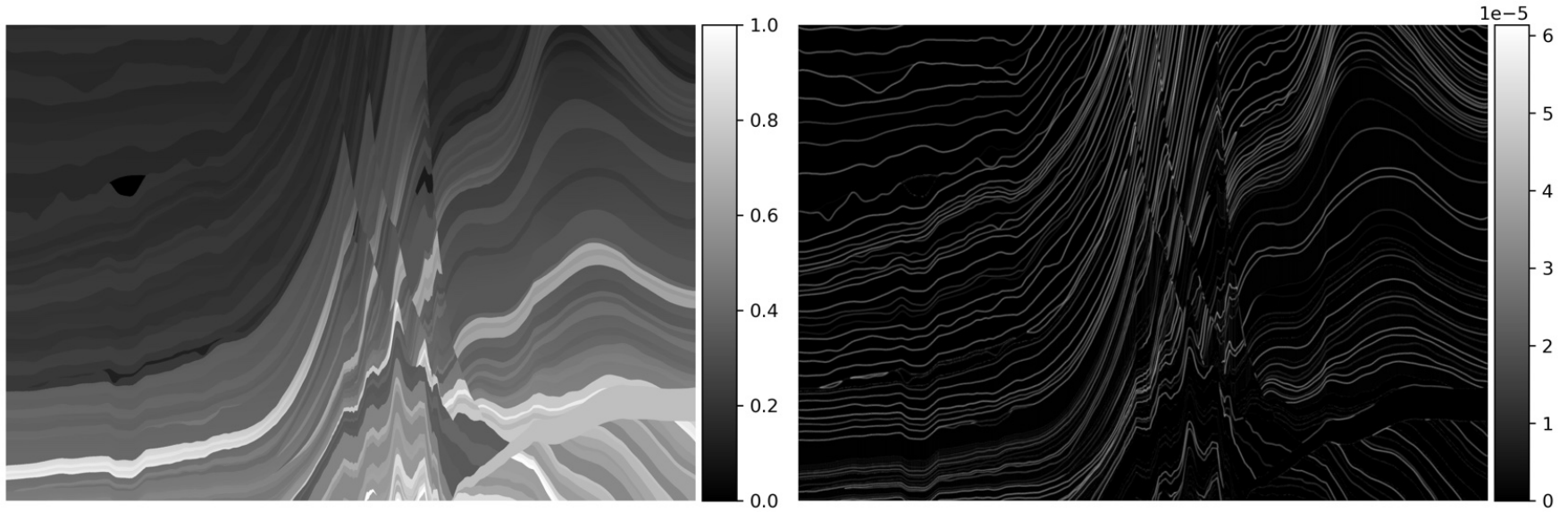}\caption{Left: The original image $\bx$ in grayscale. Right: $|\Delta_\bx\bx|$, where $\Delta_\bx$ is computed using equation~\ref{eq:edge-weight-function:applications}. As it can be seen,  $|\Delta_\bx\bx|\approx \boldsymbol{0}$ and all the edges are well identified.}\label{fig:Delta_xgt}
\end{figure}

Eventually, by all the above considerations, given an image $\bx \colon P \to [0,1]$ we can then define its associated graph $G=(P, w_\bx)$, and consequently the graph Laplacian on $G$, $\Delta = \Delta_\bx$, as per equation \ref{graph_laplacian}.  Let us observe that $\bz\mapsto \|\Delta_\bx\bz\|_1$ is a pseudo-metric. In some sense, it is telling us how close we are to $\bx$. See Figure \ref{fig:Delta_xgt}, where we plot $|\Delta_\bx\bx|$.

\subsection{The standard \texttt{graphLa}$\Psi$ method: convergence and stability results}
For the reader convenience, we report here a convergent result about the \texttt{graphLa}$\Psi$  method which makes it a regularizing method. It was proven in \cite{bianchi2023data}.

By a reconstruction method $\Psi$, we more broadly mean a family of operators of the form $\{\Psi_\Theta : Y \rightarrow X \mid \Theta \in \R^k\}$, which we call reconstructors. The parameters $\Theta$ may depend on $\delta$ and $\by^\delta$.

For instance, any variational method as in equation~\ref{genTik} qualifies as a reconstructor, with  parameter $\Theta=\alpha \in (0,+\infty)$. In this case, some specific parameter choice rules $\alpha=\alpha(\delta,\by^\delta)$ can make the pair $(\Psi_\alpha, \alpha)$ a regularization method, as per \cite[Definition 3.1]{engl1996regularization}.

Another example of interest is when $\Psi_\Theta$ is a trained DNN, where $\Theta$ denotes the parameters of the DNN. Even if DNNs can be theoretically trained to take into account different kind of noise levels $\delta$ and noise distributions, in practice this is impossible due the paramount computational time it would require. Therefore, DNNs are typically trained on just a few fixed level of noise intensities or noise distributions. So, in this case, $\Theta$ does not depend on  $\delta$ nor on $\by^\delta$, that is, $\Theta(\delta, \by^\delta) \equiv \hat{\Theta}$ for a fixed~$\hat{\Theta}$.

Chosen then a reconstructor $\Psi_\Theta$ and defining $\Psi_\Theta^\delta\coloneqq \Psi_\Theta(\by^\delta)$, the approximated solution $\bx_\alpha$ given by the \texttt{graphLa}$\Psi$ method reads


\begin{equation}\label{graphModel}
	\bx_\alpha \in  \underset{\bx \in X}{\argmin} \left\{\frac{1}{2}\|K\bx - \by^\delta\|_2^2 + \alpha \|\Delta_{\Psi^\delta_\Theta}\bx\|_1\right\}.
\end{equation}

Let $\bx_0 \in X$ and $\Theta=\Theta(\delta,\by^\delta)$ be such that 
\begin{equation}\label{eq:limit_solution}
	\bx_0= \lim_{\delta \to 0} \Psi_{\Theta(\delta,\by^\delta)} (\by^\delta).
\end{equation}
We have the following definition of solution.
\begin{definition}
	We call $\bx$ a \emph{graph-minimizing solution} with respect to $\bx_0$,	if  
	\begin{equation}
		\begin{aligned}
			&\textnormal{(i)} \; K\bx = \by, \\
			&\textnormal{(ii)}\; 	 \|\Delta_{\bx_{0}} \bx\|_1 = \min\{  \|\Delta_{\bx_{0}} \bx\|_1 \mid \bx \in X, \; K\bx = \by \}.
		\end{aligned}
	\end{equation}
\end{definition}
The next theorem provides a convergence result. For more details, see \cite[Section 3]{bianchi2023data}.

\begin{theorem}\label{thm:convergence}
	Fix a sequence $\{\delta_k\}$ and $\alpha \colon \R_+ \to \R_+$ be such that 
\begin{equation*}
	\|\by^{\delta_k}-\by\|_2 \leq \delta_k, \quad \lim_{k\to\infty} \delta_k = 0,    \quad 	\lim_{\delta_k \to 0} \alpha(\delta_k) =0, \quad \lim_{\delta_k \to 0} \frac{\delta_k^2}{\alpha(\delta_k)}=0.
\end{equation*}
 Then every sequence $\{\bx_k\}$ of elements that minimize the functional  \ref{graphModel}, with $\delta_k$ and $\Theta(\delta_k, \by^{\delta_k})$, has a convergent subsequence. The limit $\bx$ of the convergent subsequence $\{\bx_{k'}\}$is a graph-minimizing solution with respect to $\bx_0$.	If $\bx$ is unique, then $\bx_k \to \bx$. 
\end{theorem}

\section{The iterated \texttt{graphLa}$\Psi$ method: \texttt{it-graphLa}$\Psi$}\label{sec:iterated_graphLaPsi}
Iterative variants of equation~\ref{genTik} are very popular and have a long story.  Iterated Tikhonov methods can often  converge to   approximated solutions of higher quality in many applications. For example,
we refer the readers to  \cite{hanke1998nonstationary,bachmayr2009iterative,jin2014nonstationary,bianchi2015iterated,bianchi2023uniformly,bianchi2025convergence}.

In this work we propose an iterative version of \ref{graphModel} which reads

\begin{equation}\tag{\texttt{it-graphLa}$\Psi$}\label{iterated_graphLaNet}
	\begin{cases}
		\bx_n = \underset{\bx \in X}{\argmin}\left\{ \frac{1}{2}\|K\bx - \by^\delta \|_2^2 + \alpha_n\|\Delta_{\bx_{n-1}} \bx\|_1 \right\} & \mbox{for } n\geq 1,\\
		\bx^\delta_0 \coloneqq \Psi(\by^\delta),
	\end{cases}
\end{equation}
where \red{$\Psi$ is any initial reconstructor}. That is, at each step $n$ we build a new graph Laplacian $\Delta_{\bx_{n-1}}$ on the previous approximated solution $\bx_{n-1}$. For $n=1$ we recover exactly equation~\ref{graphModel}, with $\bx_0^\delta = \Psi_\Theta^\delta$.

To address the optimization problem of \ref{iterated_graphLaNet}, we employ a Majorization–Minimization technique along with a Generalized Krylov Subspace approach (MMGKS)  at each step, reducing its complexity and computational cost.  \red{The gist of MMGKS is to project into an appropriate low dimensional Krylov subspace the original problem and approximating the $\ell^2-\ell^1$ functional with a smooth and uniformly convex one. We pre-set all the parameters beforehand so the minimizer is computed automatically at each iteration. The dimension $d$ of the projected subspace is fixed at $d=50$ and the regularization parameter $\alpha_n$ is determined by the discrepancy principle.} For an in-depth explanation of the algorithm, we refer to  \cite{lanza2015generalized}.

\red{Clearly, this is not the only feasible approach. For example, another computationally attractive method for $\ell^1$ regularized problems is the variable projection augmented Lagrangian algorithm, \cite{chung2023variable}}.

\red{Notice that, in \ref{iterated_graphLaNet}, the regularization term $\mathcal{R}_n=\|\Delta_{\bx_{n-1}} \bx\|_1$ may change at every iteration. If it were instead fixed $\mathcal{R}_n = \mathcal{R} =\|\Delta_{\bx_{0}} \bx\|_2^2$ for every $n$, then convergence and stability results would follow from  \cite{jin2014nonstationary} and \cite{buccini2017iterated}.}

We report hereafter the pseudo code of the proposed algorithm.

\begin{algorithm}
\caption*{\texttt{it-graphLa}$\Psi$ algorithm}
\begin{algorithmic}[1] 
\State \textbf{Input:} $K$, $\by^\delta$, $\Psi(\by^\delta)$, $\sigma$, $R$, $N$ 
\State \textbf{Output:} $\bx_N$
\State Initialize $\bx^\delta_0 = \Psi(\by^\delta)$
\For{$n\leq N$}
    \State given $\bx_{n-1}$, $\sigma$, $R$, compute $\Delta_{\bx_{n-1}}$ 
    \State find $\bx_n$ minimizer of $ \frac{1}{2}\|K\bx - \by^\delta \|_2^2 + \alpha_n\|\Delta_{\bx_{n-1}} \bx\|_1$ by MMGKS
\EndFor
\State Return $\bx_N$
\end{algorithmic}
\end{algorithm}

\section{Numerical Experiments}\label{sec:numerical-experiments}

\color{darkred}

Following, we present different experiments to evaluate the proposed method. To demonstrate the versatility of our approach, we use four different initialization methods $\Psi$. We perform two numerical experiments, one with artificial data from the SEAM and Marmousi2 model data under different levels of noise, and a second experiment on field data from the Volve oil field. However, we first discuss the error measures used for evaluation.

\subsection{Error measures}

The proposed method is designed to recover and improve structural details in the impedance profile, i.e., it reduces the effects of noise and sharpens edges within the data. On the flipside, a stronger regularization generally leads to overall smaller impedance values, the data gets damped by some factor. This leads to a tradeoff where we can choose a strong regularization with sharp edges but overall smaller values, or a weak regularization with less sharp details but more accurate impedance values. We argue that the first case is generaly preferable since the data can simply be rescaled while the edge details are lost in the second case. However, this leads to a problem with commonly used error measures such as the Mean Squared Error (MSE) or the (Peak-)Signal-to-Noise Ratio (PSNR/SNR). These measures are far less sensitive to blurry or misplaced edges which only influence a small percentage of the data than to impedance values of large layers being slightly off. Hence, these methods will favor a weak regularization instead of sharp contours.

For this reason, we are using two different error measures in the presented experiments. As a first measure, we use the MSE but applied to the data differential along the time axis. Therefore, let $\nabla$ indicate a first order differential operator, i.e., in the discrete case a difference operator of the form

\begin{align*}
    \nabla=\begin{pmatrix}
        1 & -1 & 0 & \hdots & 0 \\
        0 & 1 & -1 & \ddots & \vdots \\
        \vdots & \ddots & \ddots & \ddots & 0 \\
        0 & \hdots & 0 &  1 & -1
    \end{pmatrix}.
\end{align*}

Then we calculate the MSE on the differential (D-MSE) of the reconstructed solution $\bx_\text{rec
}$ compared to the ground trouth $\bx_\text{true}$ by

\begin{align*}
    \text{D-MSE}(\bx_\text{rec})&=\frac{1}{\|\nabla\bx_\text{true}\|_0}\sum\limits_{p,q}(\nabla\bx_\text{rec}(p,q)-\nabla\bx_\text{true}(p,q))^2,\\
    &\mbox{where } \|\nabla\bx_\text{true}\|_0 \coloneqq \left|\{(p,q) \mid \nabla\bx_\text{true}(p,q) \neq 0\} \right|.
\end{align*}

Here we choose to scale the value by the number of nonzero entries in $\nabla\bx_\text{true}$, rather than with the data size, which more accurately scales with the number of layers in the data and reduces the effect of large homogeneous areas. Besides this measure, we also use the Structural Similarity Index Measure (SSIM) \cite{wang2004image}. This index compares the similarity of two images by calculating
\begin{equation}\label{ssim}
    \text{SSIM}(x,y)=\frac{(2\mu_x\mu_y+c_1)(2\sigma_{xy}+c_2)}{(\mu_x^2+\mu_y^2+c_1)(\sigma_x^2+\sigma_y^2+c_2)}
\end{equation}
where $\mu_x$, $\mu_y$ are the (sample) mean value, $\sigma_x$, $\sigma_y$ are the (sample) variance, and $\sigma_{xy}$ is the (sample) covariance. Small constants $c_1,c_2\approx O(10^{-4})$ are used to ensure a non-zero denominator. Equation \ref{ssim} is calculated over all $11\times11$ pixel sized windows and the mean value over all windows is returned. The SSIM index returns a value in $[-1,1]$ where $1$ is a perfect match between both images. To reduce the influence of the damping effect we furthermore normalize both the true and the reconstructed impedance before calculating the SSIM, i.e., both data sets will have mean zero and variance one. (Note that the mean value and variance used in equation \ref{ssim} are calculated on each $11\times11$ window and are thus not necessarily normalized.)

\subsection{Initialization methods}

The \ref{iterated_graphLaNet} requires a starting guess $\bx^\delta_0 \coloneqq \Psi(\by^\delta)$. Here $\Psi$ can in principle be any  inversion method. To show that the proposed algorithm can be combined with a wide variety of inversion methods, we will use four different strategies in our experiments: two neural network based approaches and two classical reconstructors, where one method of each category performs a simple trace-wise inversion while the other integrates inter-trace correlations in the process. 

\subsubsection{Neural network approaches}

Neural networks are a natural choice as an initialization method as they are Lipschitz-continuous maps which aligns well with the theory of  \texttt{graphLa+}$\Psi$ \cite[Example 3.3]{bianchi2023data}. We will use two different network based approaches in this work. A trace-wise reconstruction by  \cite{Alfarraj2019network1} (AA) and an improved version by  \cite{Mingming2024network2} (Liu)  that also considers inter-trace correlations.

In \cite{Alfarraj2019network1} a new impedance inversion approach was suggested, that simultaneously learns both the forward and inverse operator. Two neural networks $\Psi_{\Theta_1}$ and $\Psi_{\Theta_2}$ are trained simultaneously where $\Psi_{\Theta_1}$ models the impedance inversion and $\Psi_{\Theta_2}$ is the forward model (here, $\Theta_1$ and $\Theta_2$ are the network parameters). The loss function used during the training reads as
\begin{equation}\label{eq:lossfct}
    \Loss(\Theta_1,\Theta_2)=\frac{\alpha}{N_p}\left\|\bx_{\text{true}}-\Psi_{\Theta_1}(\by^\delta)\right\|_{F,\Omega}^2+\frac{\beta}{N_s}\left\|\by^\delta-\Psi_{\Theta_2}(\Psi_{\Theta_1}(\by^\delta))\right\|_F^2,
\end{equation}
where $\|\cdot\|_F$ is the Frobenius norm and $N_s$ is the number of seismic traces in the dataset. Furthermore, the training requires a ground truth $\bx_{\text{true}}$ on a subset $\Omega$ of the given traces which is usually obtained from a certain number $N_p$ of well log samples. The loss function \ref{eq:lossfct} consists of two terms, one measuring the mean squared error on the known ground truth and the second one enforcing $\Psi_{\Theta_1}$ to be the inverse operator of $\Psi_{\Theta_2}$. A simple 4-layer convolutional Neural Network usually suffices as forward model $\Psi_{\Theta_2}$. For the inverse operator $\Psi_{\Theta_1}$ a more complicated architecture consisting of a combination of linear, Gated Recurrent Unit, convolutional, and deconvolutional layers is used. For a detailed description on the used architectures and training process, we refer the reader to the original work \cite{Alfarraj2019network1} as well as to the Python code \cite{Alfarraj2019PythonCode}. The results used in this work were achieved with the exact same setup as in the original works.

The above approach was recently improved in \cite{Mingming2024network2} where a different architecture for the inverse network $\Psi_{\Theta_1}$ was suggested. Since the original approach did not take inter-trace correlations into account, the reconstruction can become unstable under noise. To improve this, Liu et al. introduced a network based on attention layers. The forward network architecture as well as the training process remained unchanged. Again we refer the reader to the original work \cite{Mingming2024network2} for details.

In our experiments we use the results obtained with the above methods as starting guess for \ref{iterated_graphLaNet}. The training process is unchanged to the original works, i.e., the weights are chosen as $\alpha=0.2$, $\beta=1$ and $N_p=20$ equally spaced traces are taken from the ground truth impedance to simulate given well log samples.

\subsubsection{Classic approaches}

We use two classical impedance inversion techniques as initialization methods which are both based on a regularized optimization problem that reads as 

\begin{align}\label{eq:classical_model}
    \min\limits_{\bx}\|\by^\delta-K\bx\|_F^2+\alpha\|\nabla\bx\|_{1,1}+\beta\|\nabla\bx^T\|_{1,1}.
\end{align}

The data fidelity term fits the predicted seismic data $K\bx$ to the measurements $\by^\delta$. Since $\|\cdot\|_{1,1}$ is a sparsity promoting norm, the regularizers favor impedance profiles $\bx$ which sparse derivative in time and spatial domain. Both terms have a physical interpretation. A sparse derivative in time $\nabla\bx$ corresponds to only having a few layer boundaries at which the impedance value changes, a sparse derivative in space $\nabla\bx^T$ forces the layers to align along different traces.

Using only the first regularizer (i.e., $\beta=0$) problem~\ref{eq:classical_model} can be solved trace wise. This ansatz is known as sparse spike inversion \cite{velis2008stochastic} (SSI). By replacing $\nabla\bx=\bx_t$ we can rewrite the problem as
\begin{align*}
    \min\limits_{\bx_t}\|\by^\delta-K\nabla^{-1}\bx_t\|_F^2+\alpha\|\bx_t\|_{1,1}
\end{align*}
which is a standard compressed sensing problem. We used an interior point method to solve the above problem (\cite{kim2007interior,l1lscode}).

If both regularizers are used in \ref{eq:classical_model} the problem becomes a 2D linear inversion with total variation regularizer, which can be solved using the split bregman technique \cite{gholami2016fast} (SB). The implementation used in our experiments is based on \cite{goldstein2009split,splitbregmancode}.

Both the SSI and SB require prior knowledge about the linear forward operator $K$. In case $K$ is not known, a common approach for impedance inversion assumes that the seismic impulse is reflected from the subsurface layer boundaries. This leads to a linear forward operator model $K\bx=W\nabla\bx=W\bx_t$ where $W$ is a convolution matrix based on the seismic wave. The impedance profile and seismic wave can be reconstructed simultaneously by solving a blind deconvolution problem. We state in more detail on how the operator $K$ was obtained for each experiment individually.

\subsection{Model data experiments}

In this experiment we use two artificially created datasets, the Marmousi2 \cite{martin2002updated} and SEAM \cite{fehler2011seam}, where the seismic data was obtained using the convolutional model with the Ricker wavelet and an undersampling factor of $4$. The SEAM model is a smaller dataset with $1502$ traces while Marmousi2 is nearly twice as large with $2721$ traces. To simplify the training process of the involved neural networks in our experiments, we added additional time and depth samplings to the SEAM model such that both models have $1880$ rows for the impedance profile and $470$ rows for the seismic data. Because of the applied undersampling, the inverse problem is highly underdetermined. For impedance inversion we use the noiseless seismic data as well as noisy data with four different noise levels with a PSNR of about $39$, $33$, $30$, and $27$. We note that the Marmousi2 model data is overall more complex than the SEAM data, e.g., it has more low impedance regions, a lot of thin layers with complex geometry, and the impedance difference between two neighboring layers is smaller on average. This combination results in a higher impact of the noise on the data and the reconstructions.

\begin{figure}[h!tb]
	\centering
    \begin{tabular}{llll}
    \subfloat[]{
    \includegraphics[height=0.25\textwidth]{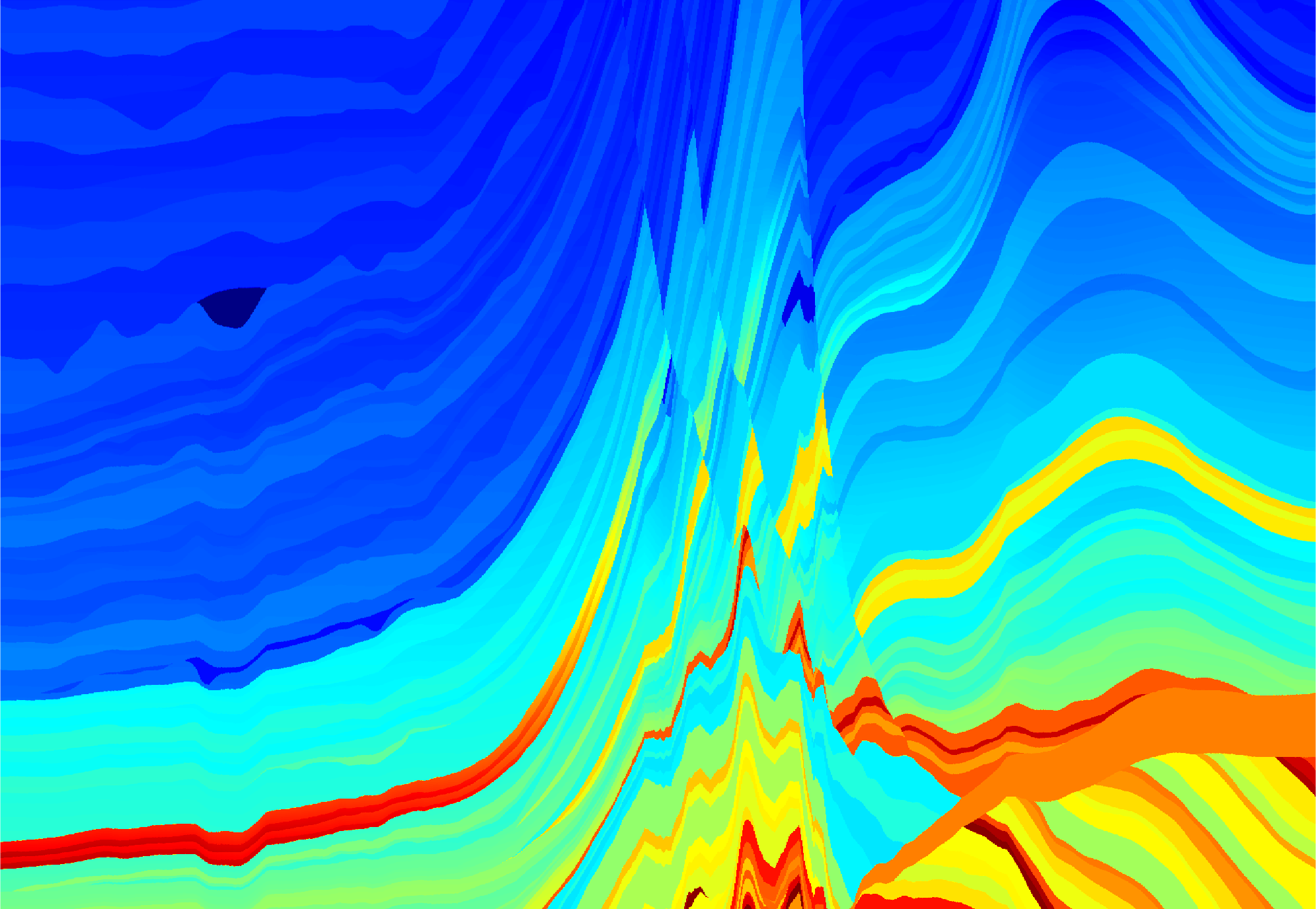}
    } & \marmousiCB{0.25\textwidth} &
    \subfloat[]{
    \includegraphics[height=0.25\textwidth]{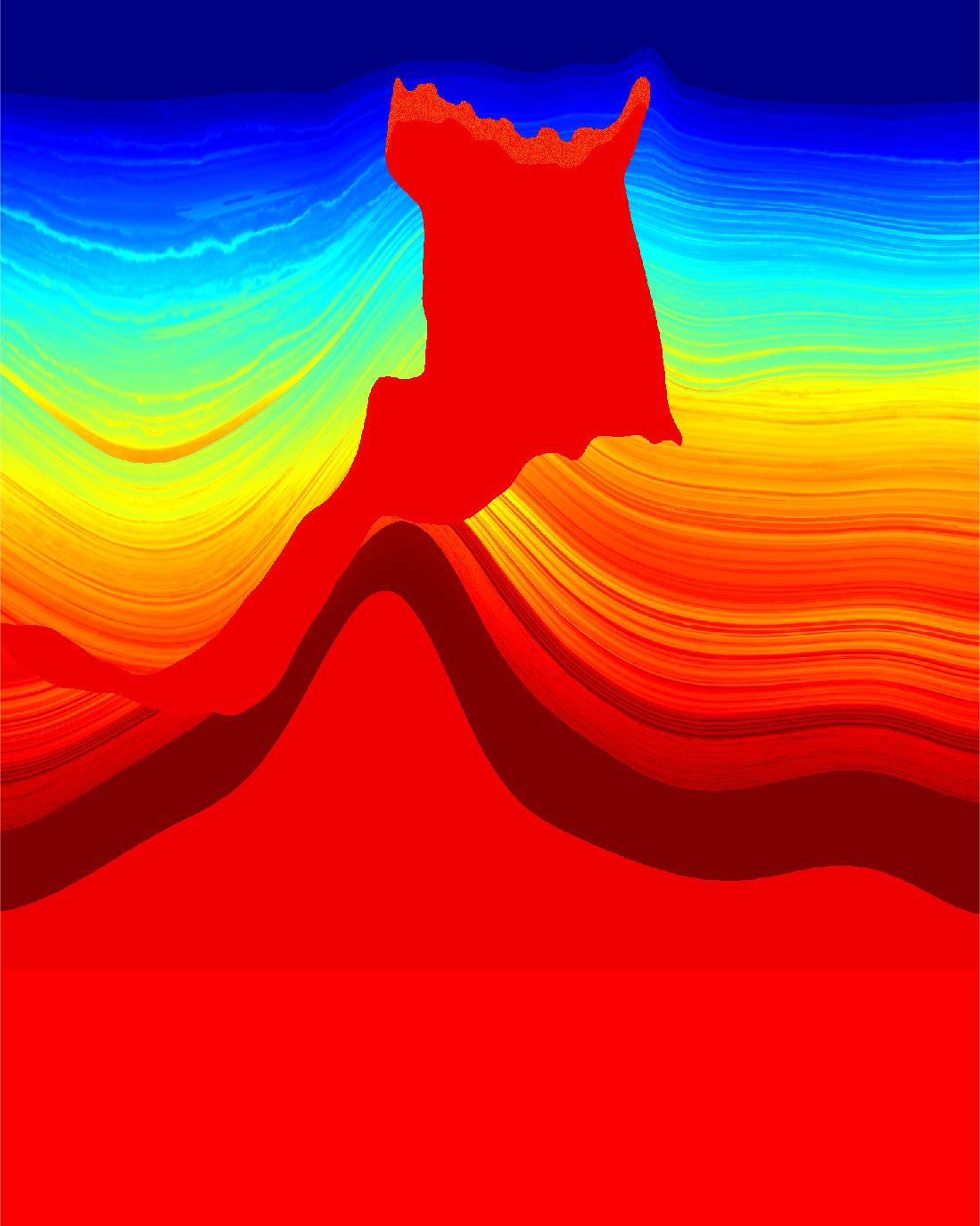}
    } & \seamCB{0.25\textwidth} \\
    \subfloat[]{
    \includegraphics[height=0.25\textwidth]{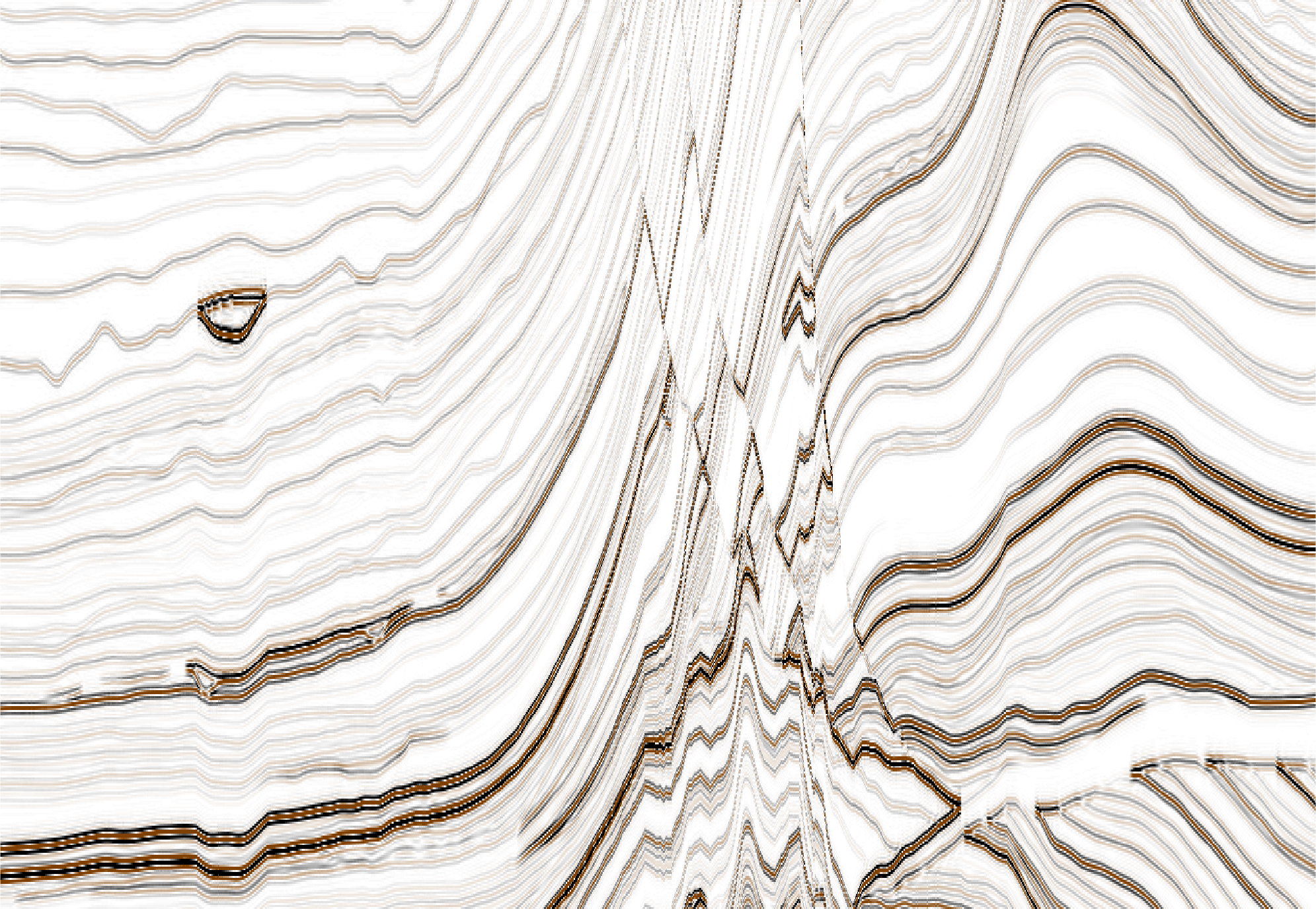}
    } & &
    \subfloat[]{
    \includegraphics[height=0.25\textwidth]{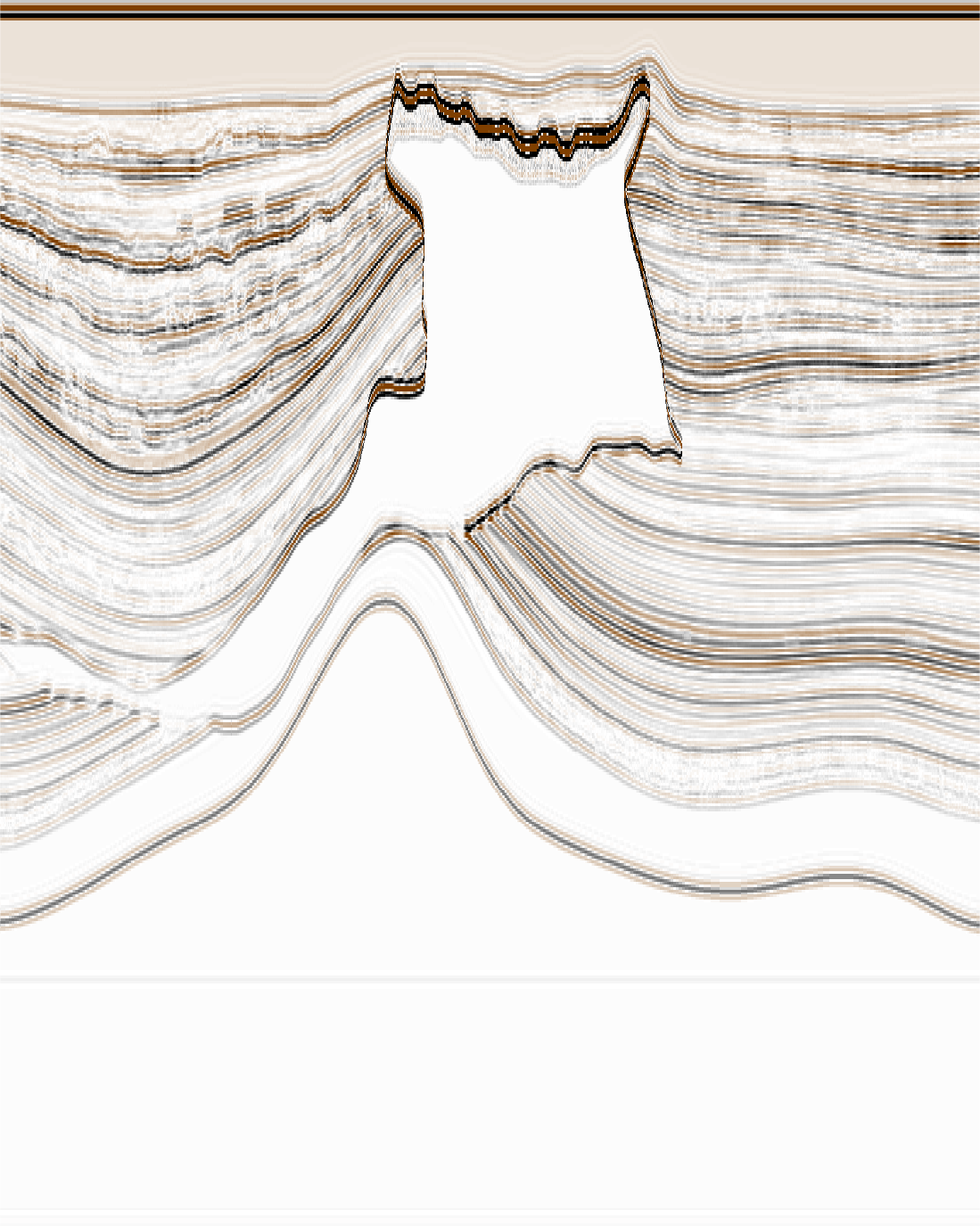}
    } & \colorbar{seis}{-5}{5}{0.25\textwidth} \\
    \subfloat[]{
    \includegraphics[height=0.25\textwidth]{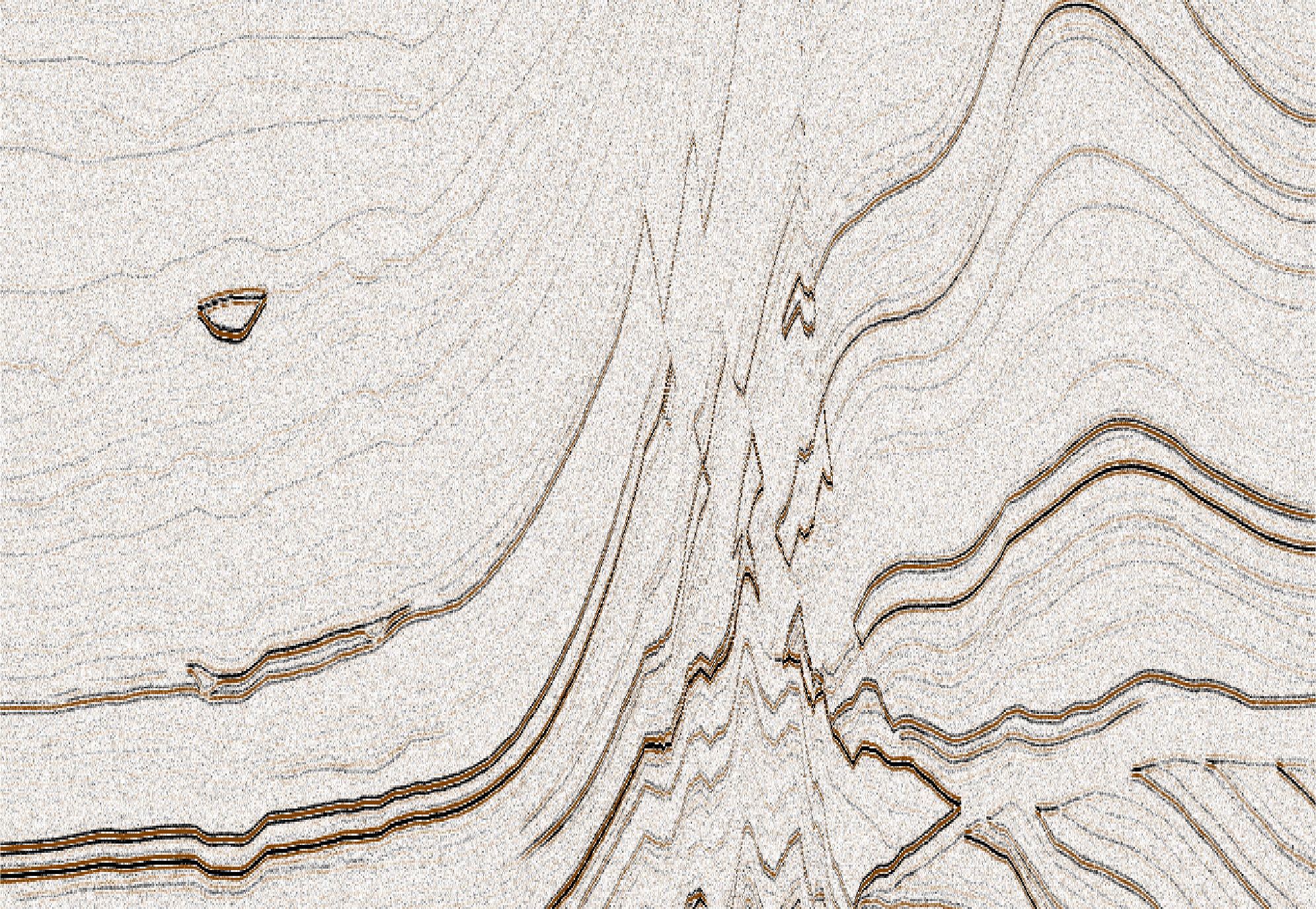}
    } & &
    \subfloat[]{
    \includegraphics[height=0.25\textwidth]{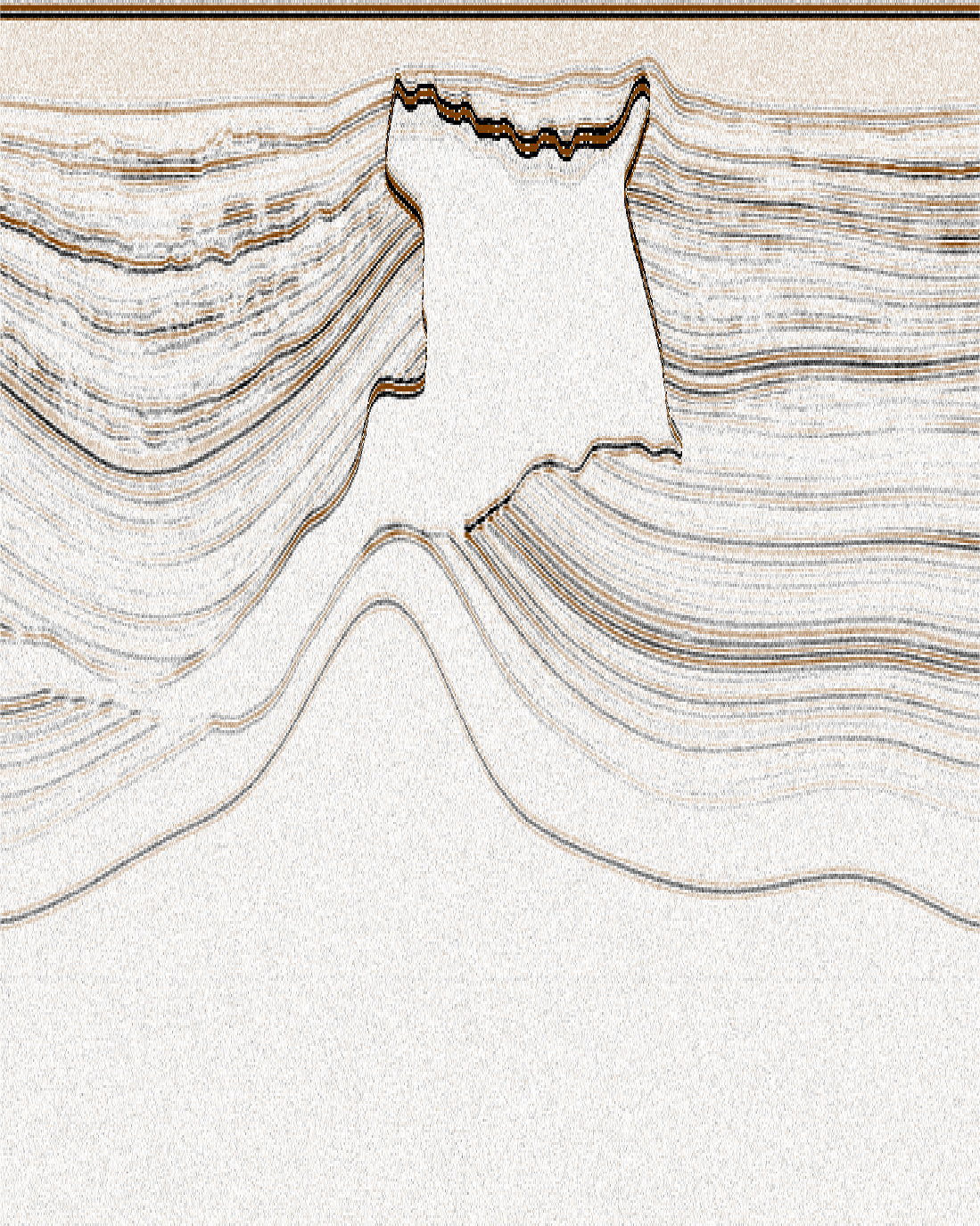}
    }\\
    \subfloat[]{
    \includegraphics[height=0.25\textwidth]{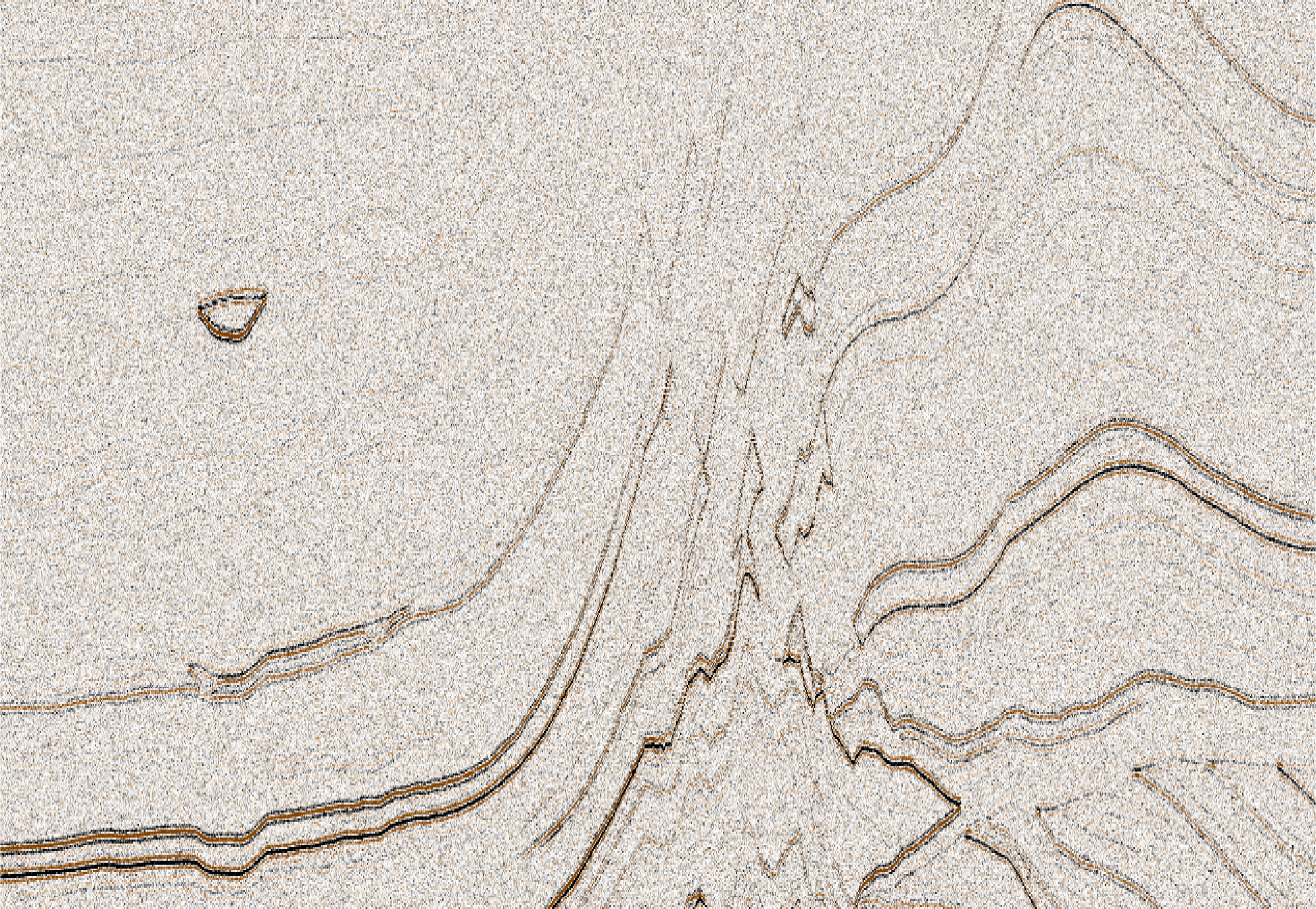}
    } & &
    \subfloat[]{
    \includegraphics[height=0.25\textwidth]{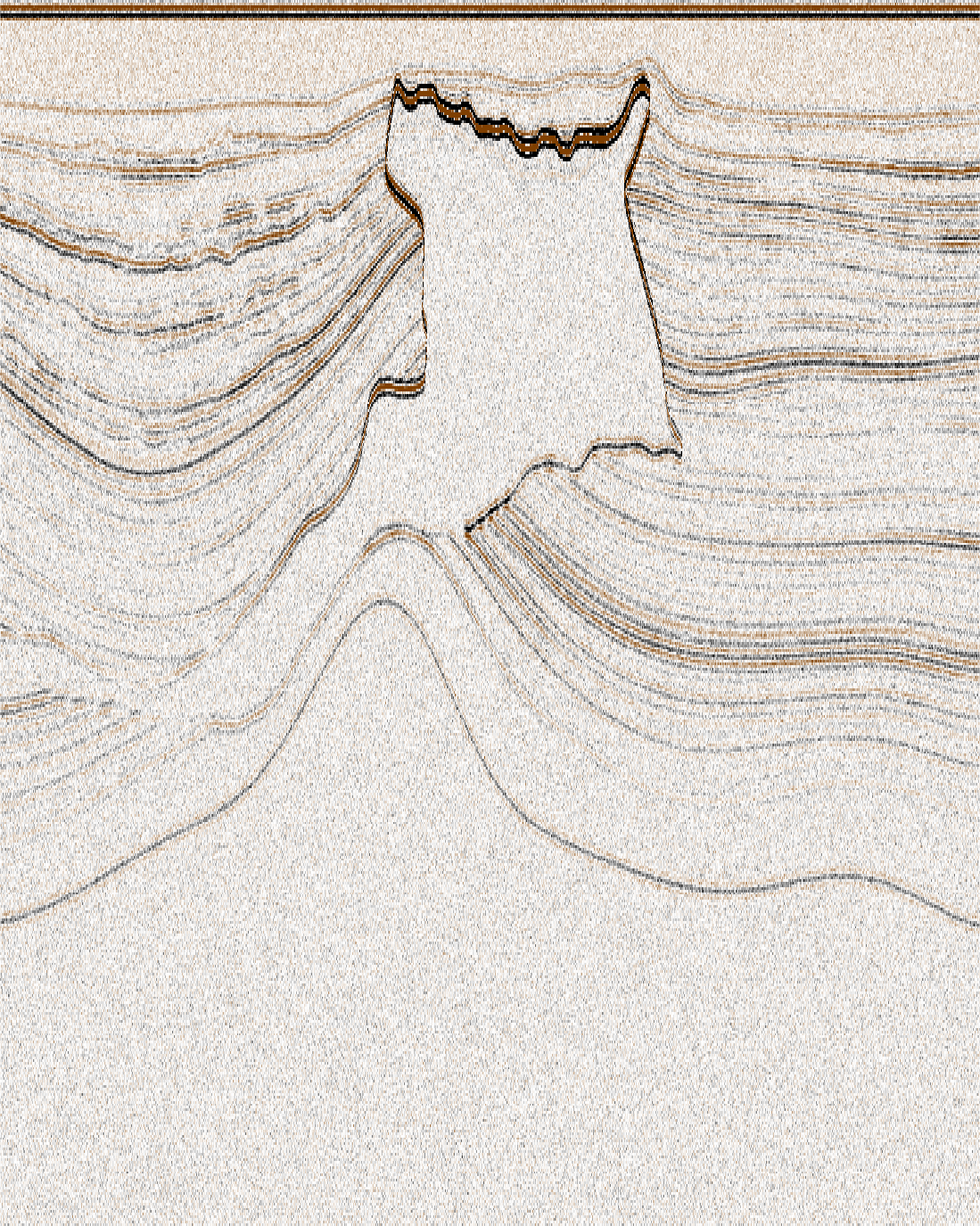}
    }
    \end{tabular}
	\caption{data from the Marmousi2 (left) and SEAM model (right): impedance profile (a,b) and seismic data without noise (c,d), medium noise PSNR $\approx33$ (e,f), high noise PSNR $\approx27$ (g,h). The colormap of all seismic images is the same and cuts of extreme values for a better contrast.}
    \label{fig:data}
\end{figure}

Since the neural network based methods require normalized data (i.e., mean zero and standard deviation one), both the impedance and seismic data sets have been normalized for this experiment. Figure~\ref{fig:data} shows the impedance profile of both models as well as the seismic data with no noise, medium noise, and high noise level.

\begin{figure}[h!tb]
	\centering
    \subfloat[]{\includegraphics[width=0.25\textwidth]{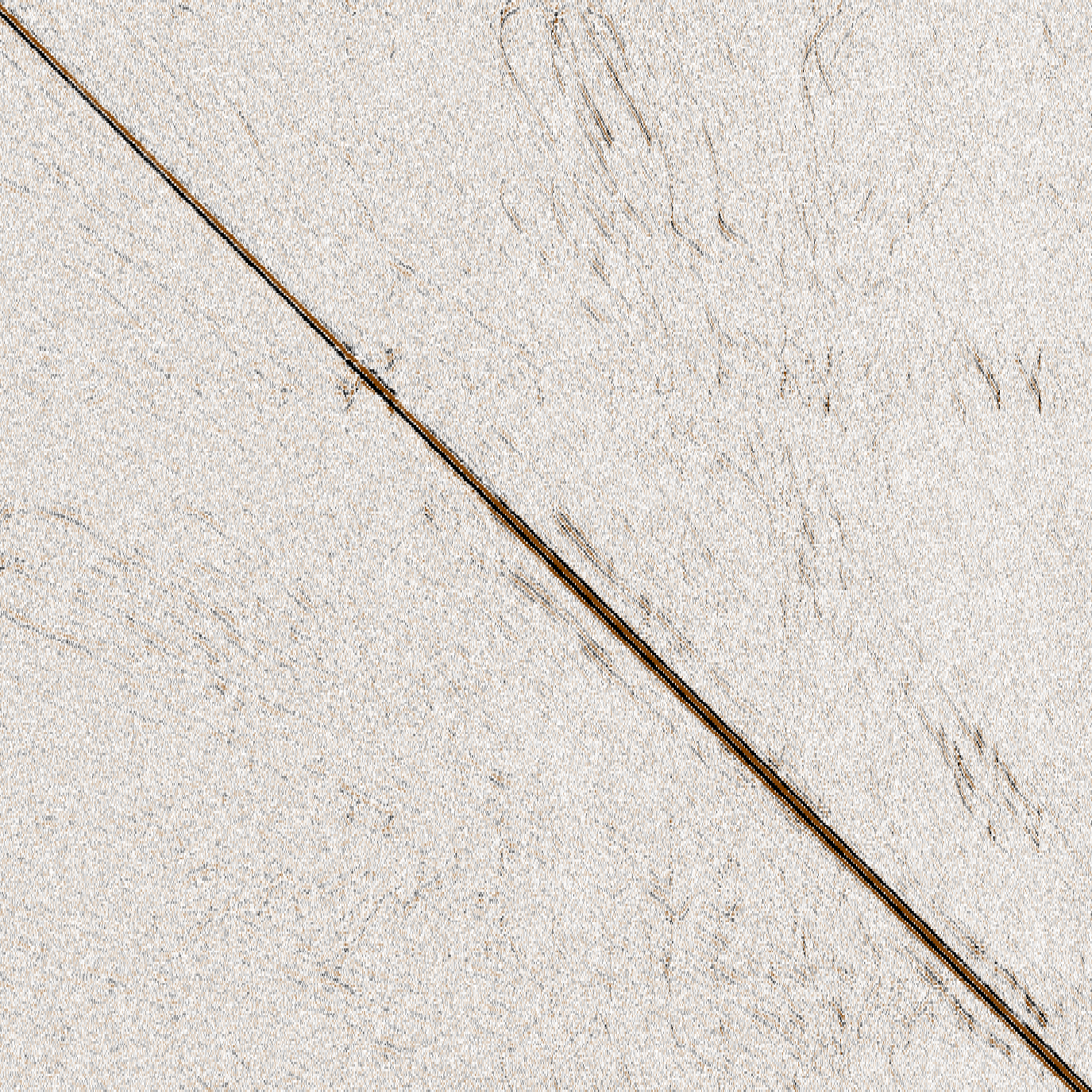}}\ \ 
    \subfloat[]{\includegraphics[width=0.25\textwidth]{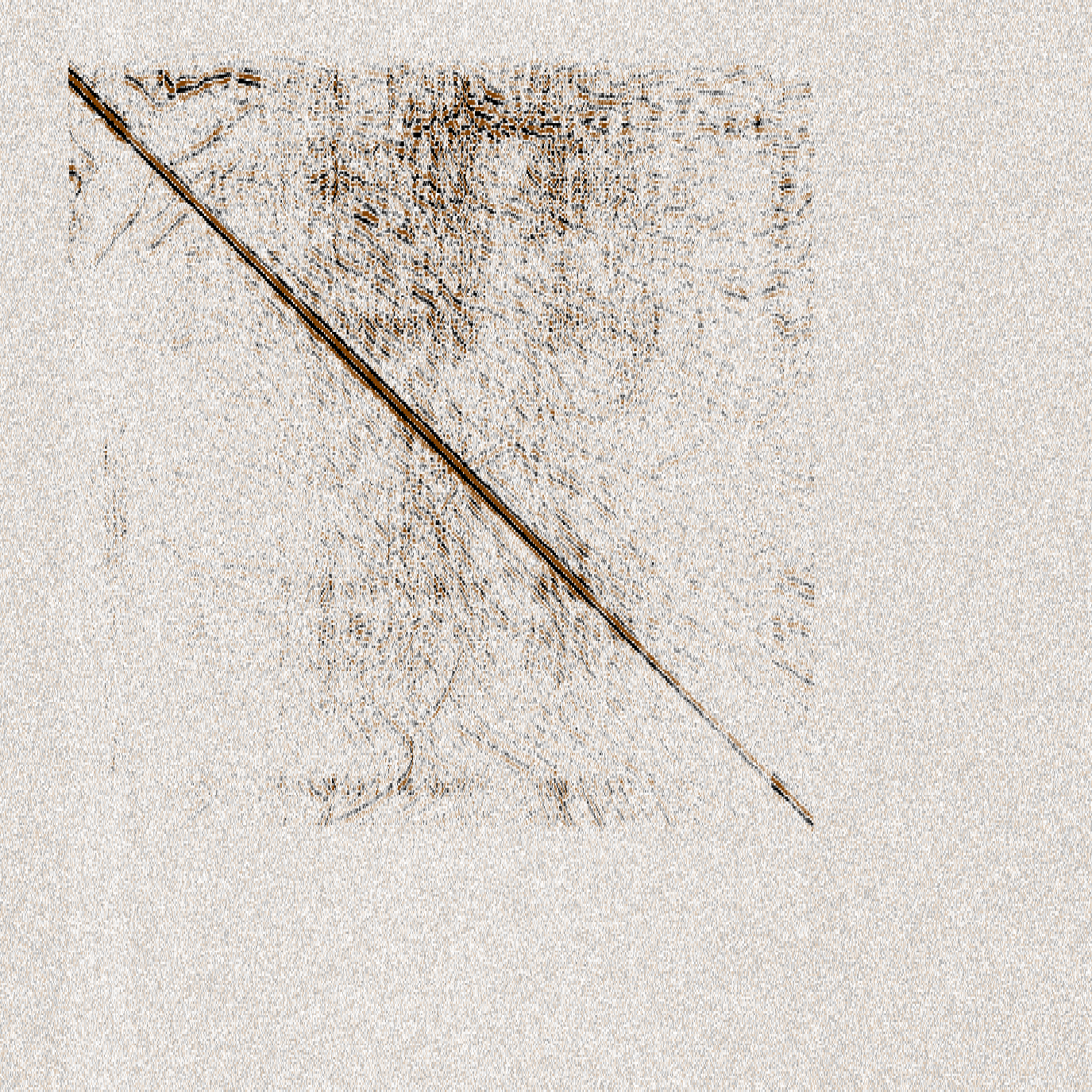}}\ \ 
    \colorbar{seis}{-0.25}{0.25}{0.25\textwidth}
	\caption{Linear forward operator $K$ learned for Marmousi2 (a) and SEAM  (b).}
    \label{fig:K}
\end{figure}

For this experiment we assume that we have perfect information on the forward model. We train the neural network initialization methods with $N_p=20$ traces of the true impedance. Furthermore, we learn the linear forward operator $K$ from the true impedance and noiseless seismic data as one layer fully connected linear neural network. The resulting operator matrix for both datasets is shown in Figure \ref{fig:K}. We see that both operators are concentrated along the main diagonal. This matches with the linear model $K=W\nabla$ mentioned in the previous subsection, where $W$ is a convolution matrix, $\nabla$ a difference operator and both matrices ($W$ and $\nabla$) are concentrated along the main diagonal. For the SEAM model $K$ has nearly vanishing values all around the boundary. This is due to areas of constant impedance at the top and bottom, which we added to have the same number of time samples in both models.

\begin{figure}[h!tb]
	\centering
    \subfloat[]{\includegraphics[height=0.22\textwidth]{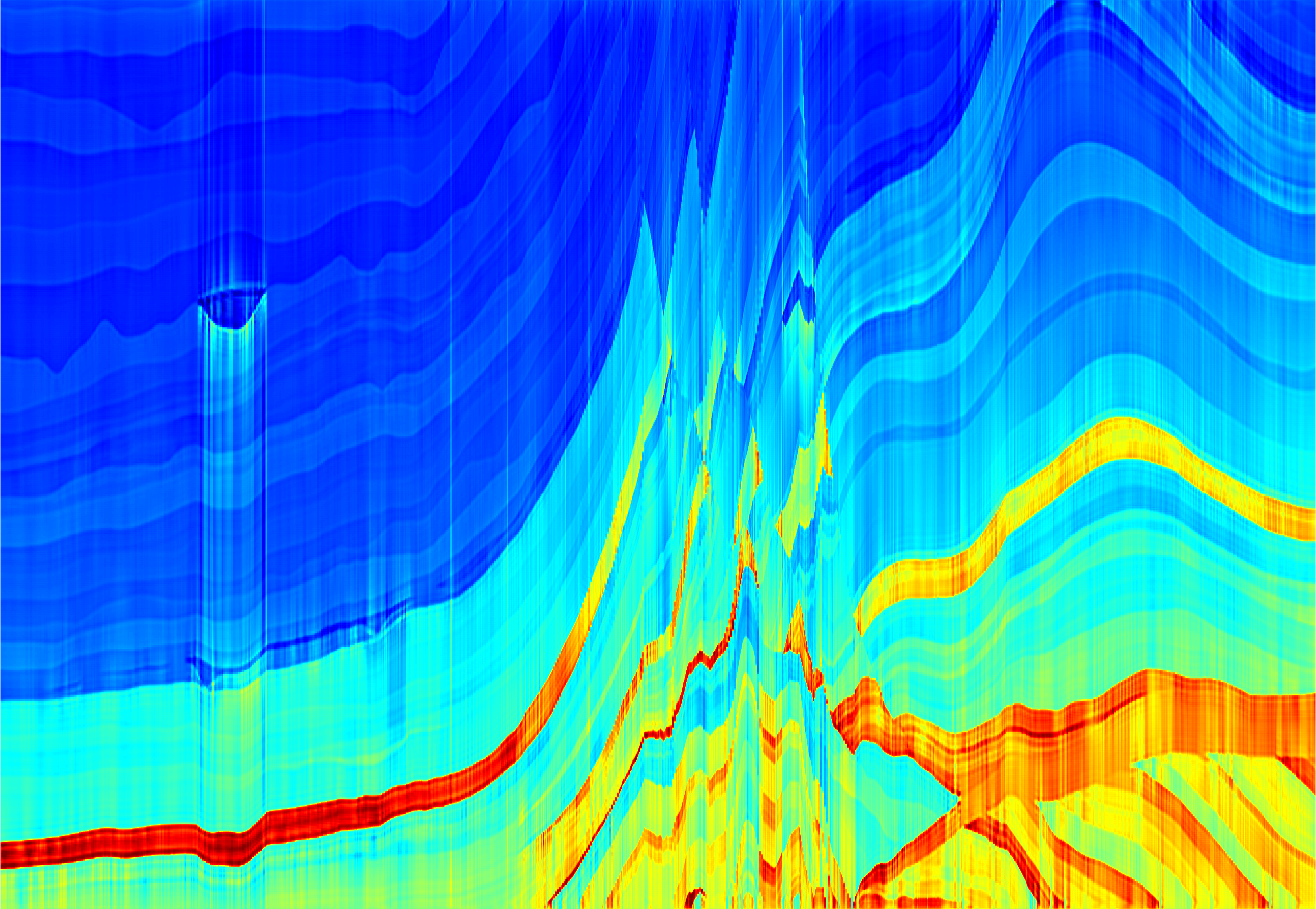}}\ \ 
    \subfloat[]{\includegraphics[height=0.22\textwidth]{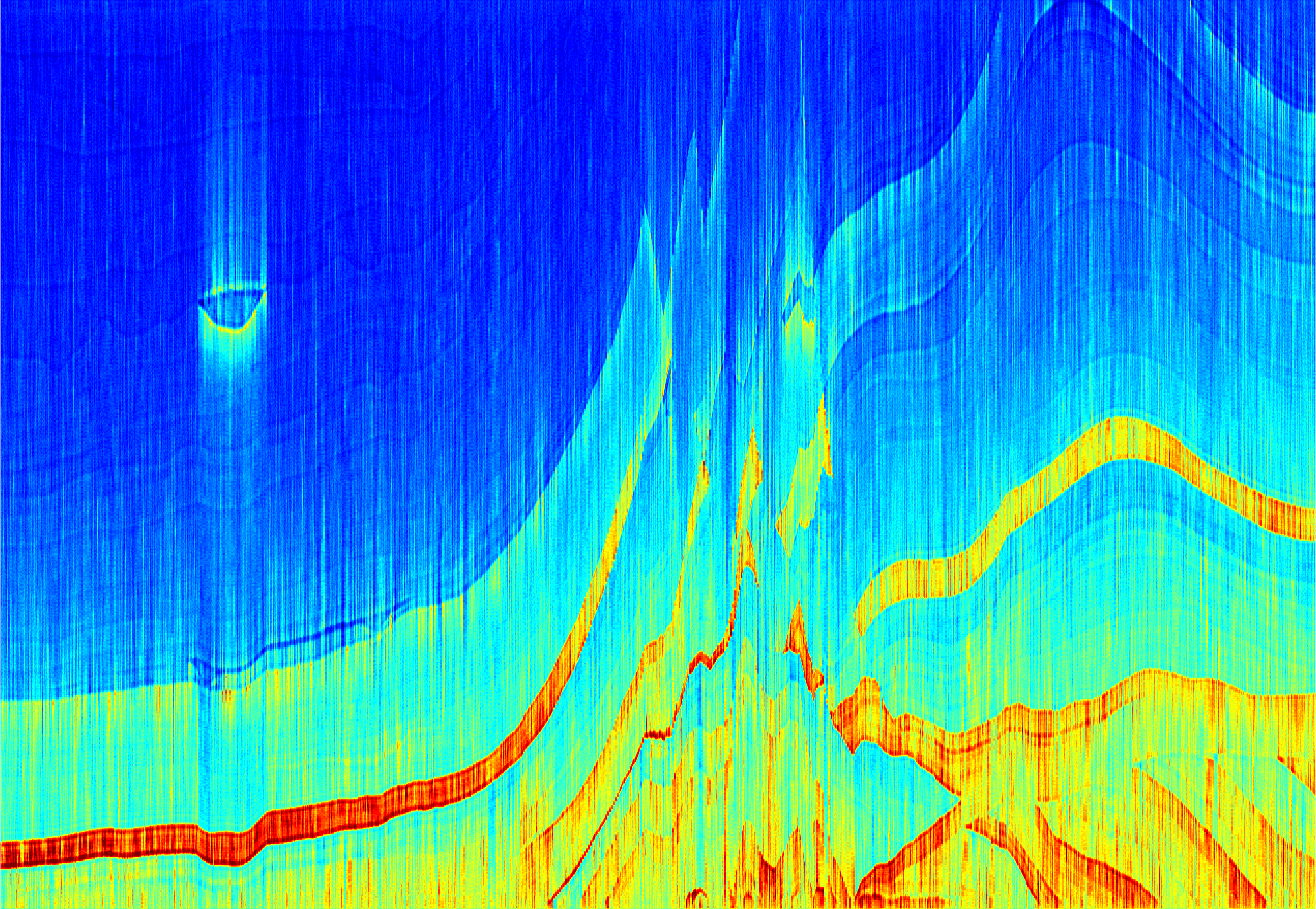}}\ \ 
    \subfloat[]{\includegraphics[height=0.22\textwidth]{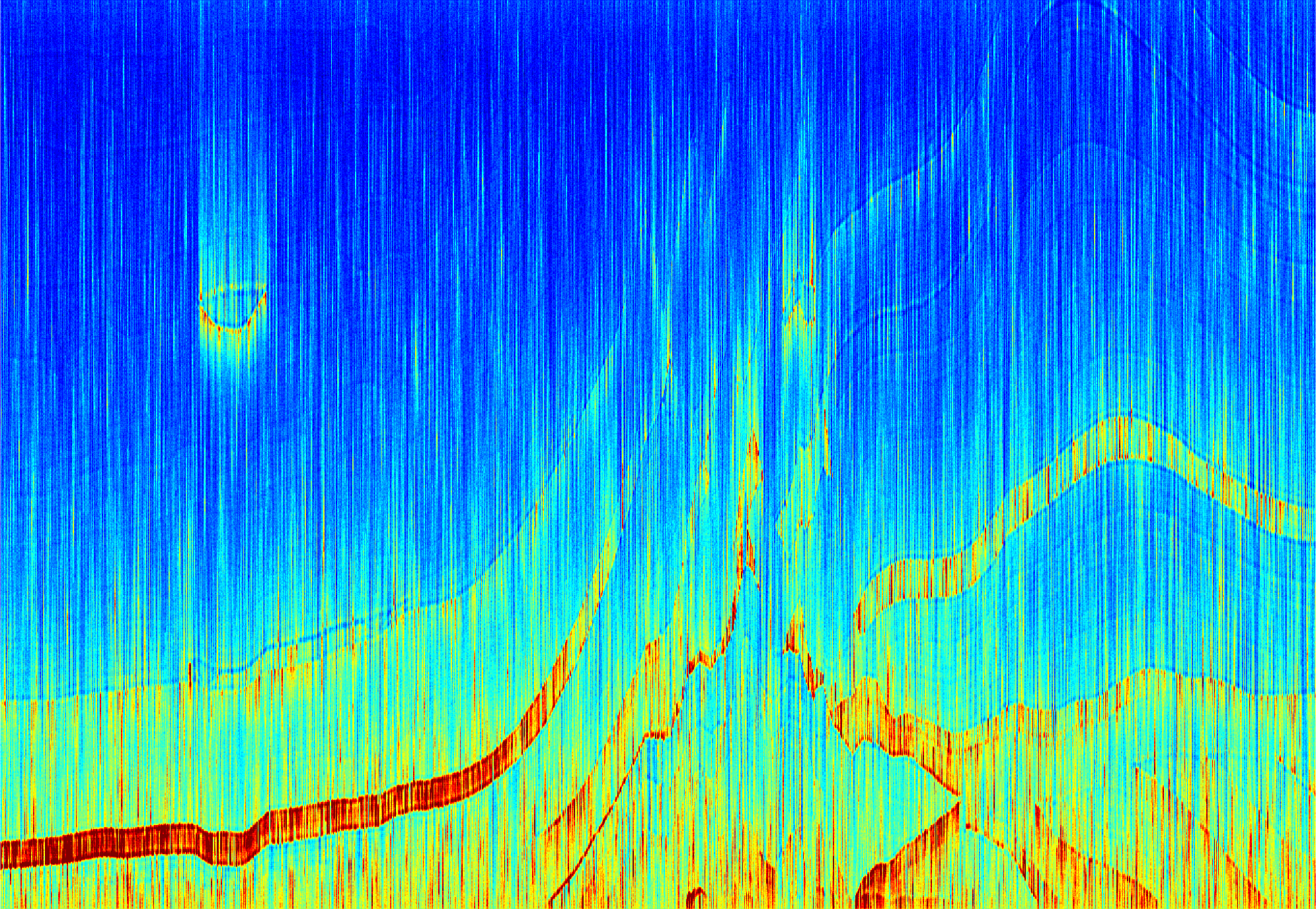}}\\
    \subfloat[]{\includegraphics[height=0.22\textwidth]{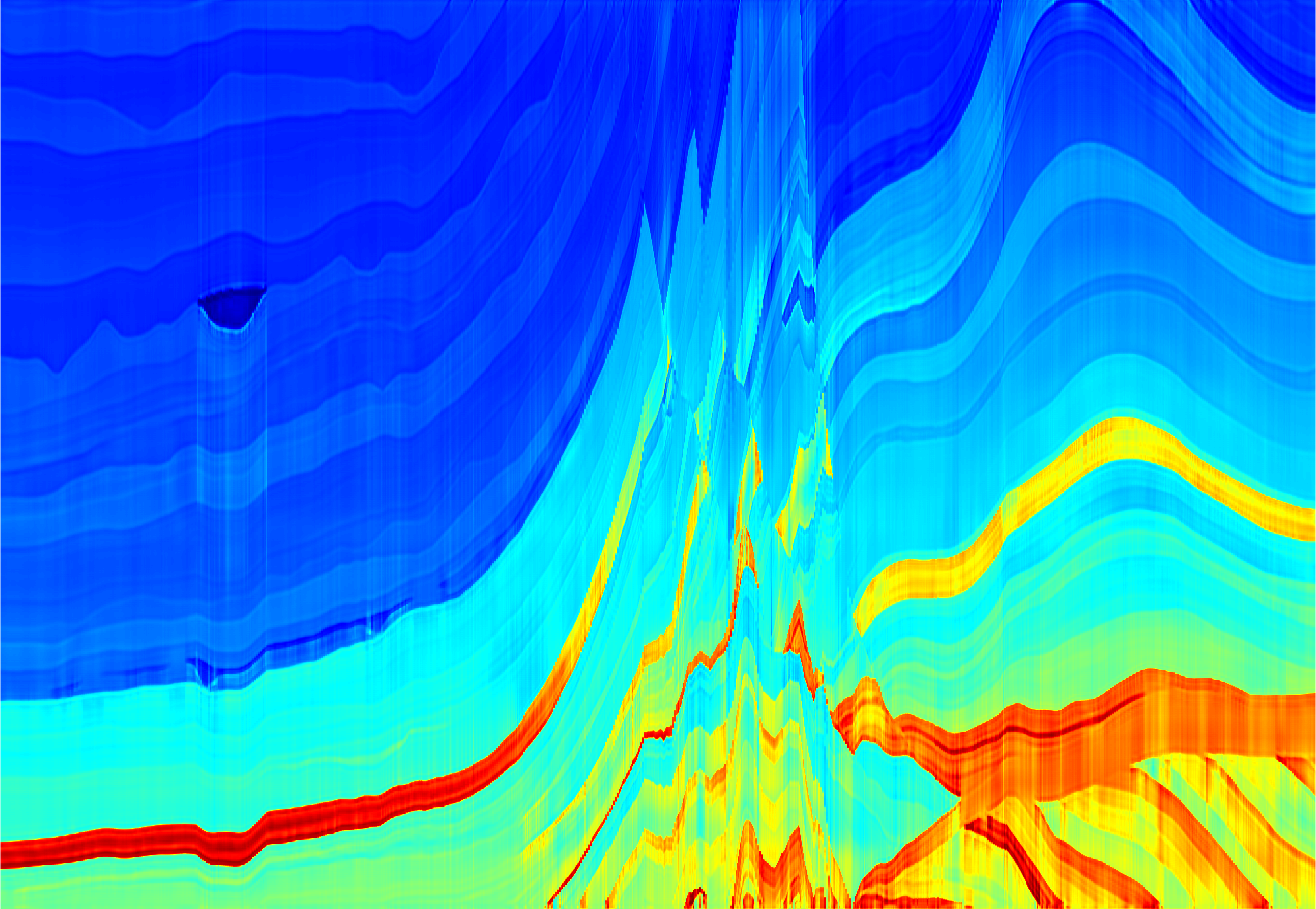}}\ \ 
    \subfloat[]{\includegraphics[height=0.22\textwidth]{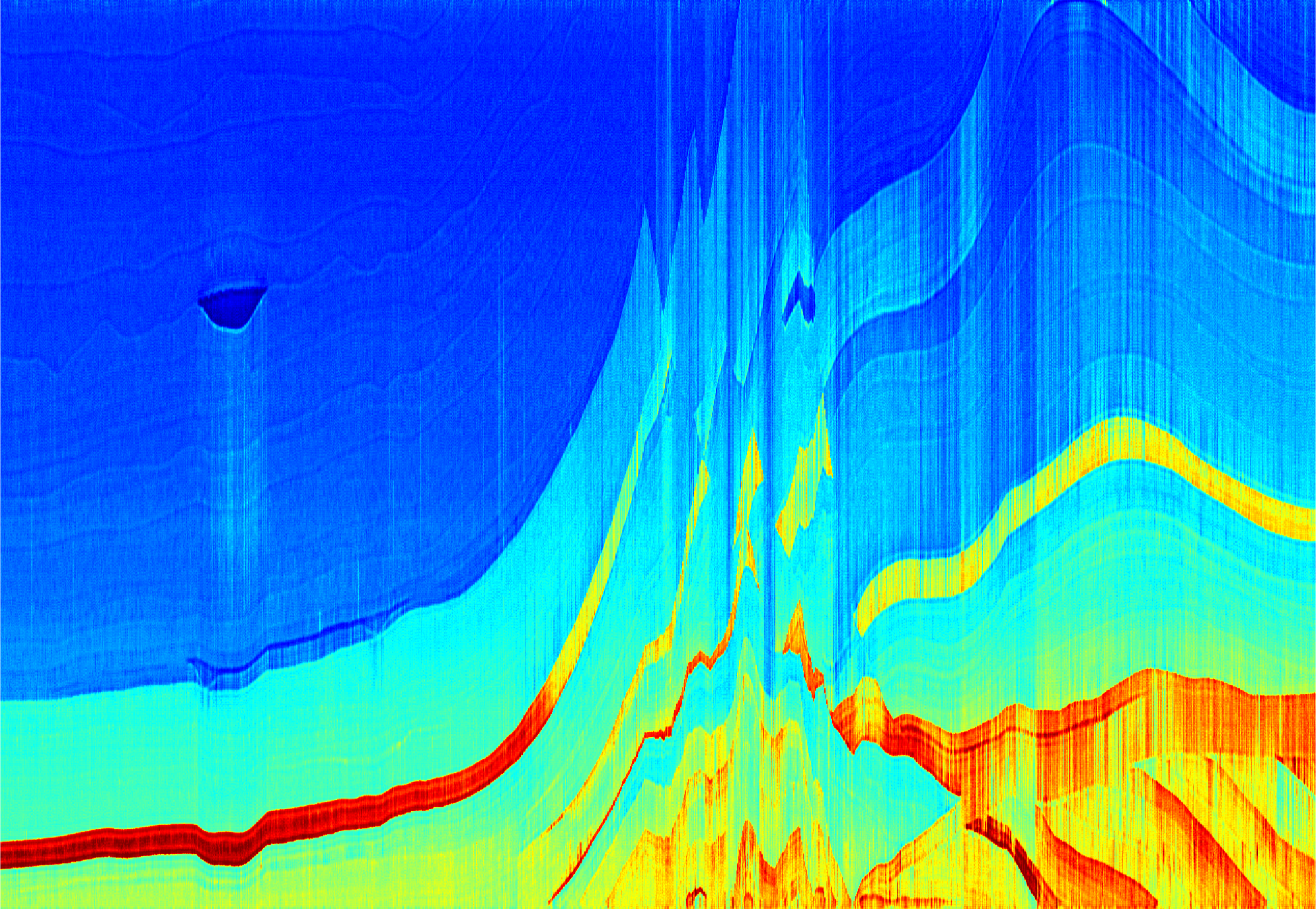}}\ \ 
    \subfloat[]{\includegraphics[height=0.22\textwidth]{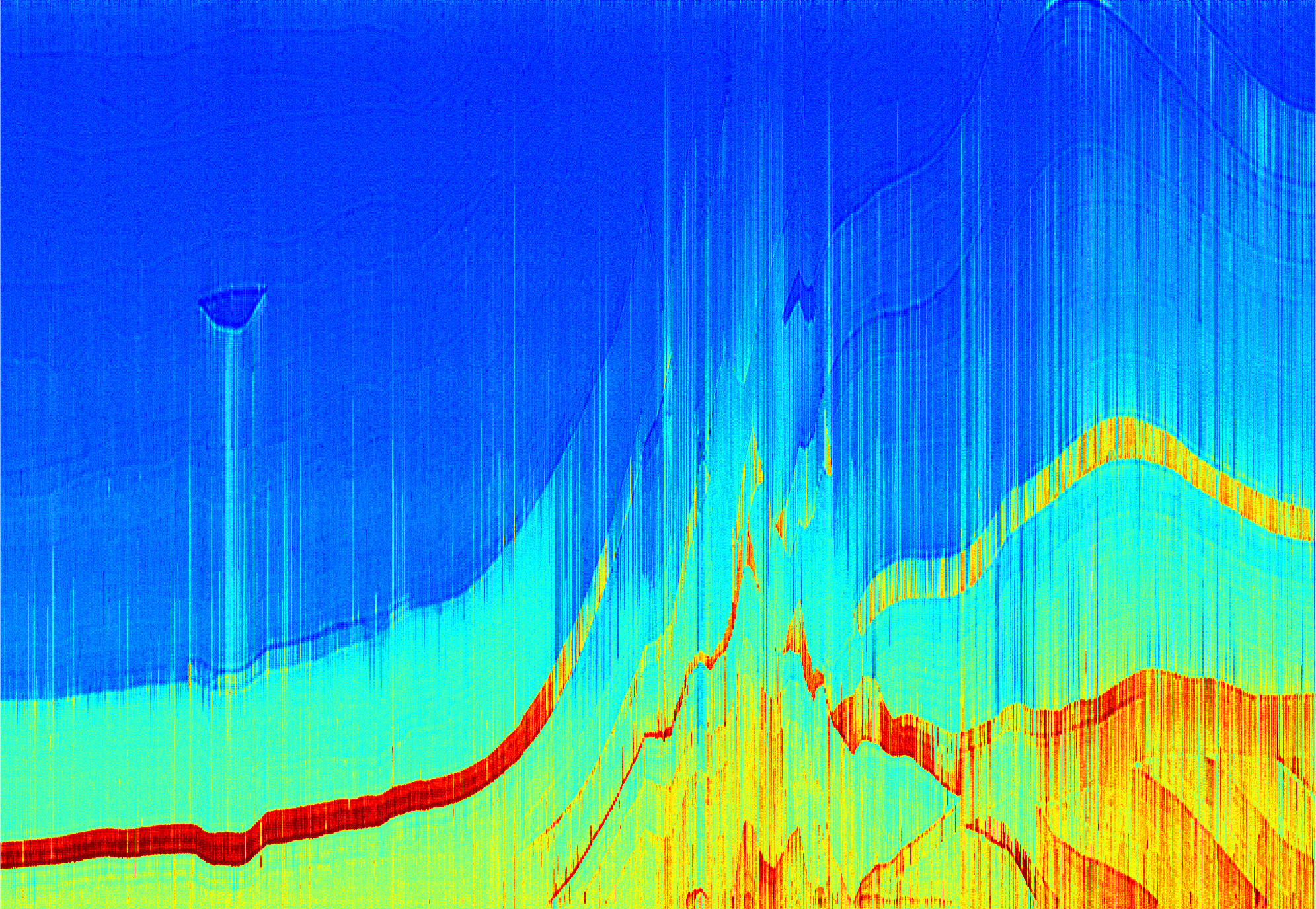}}\\
    \subfloat[]{\includegraphics[height=0.22\textwidth]{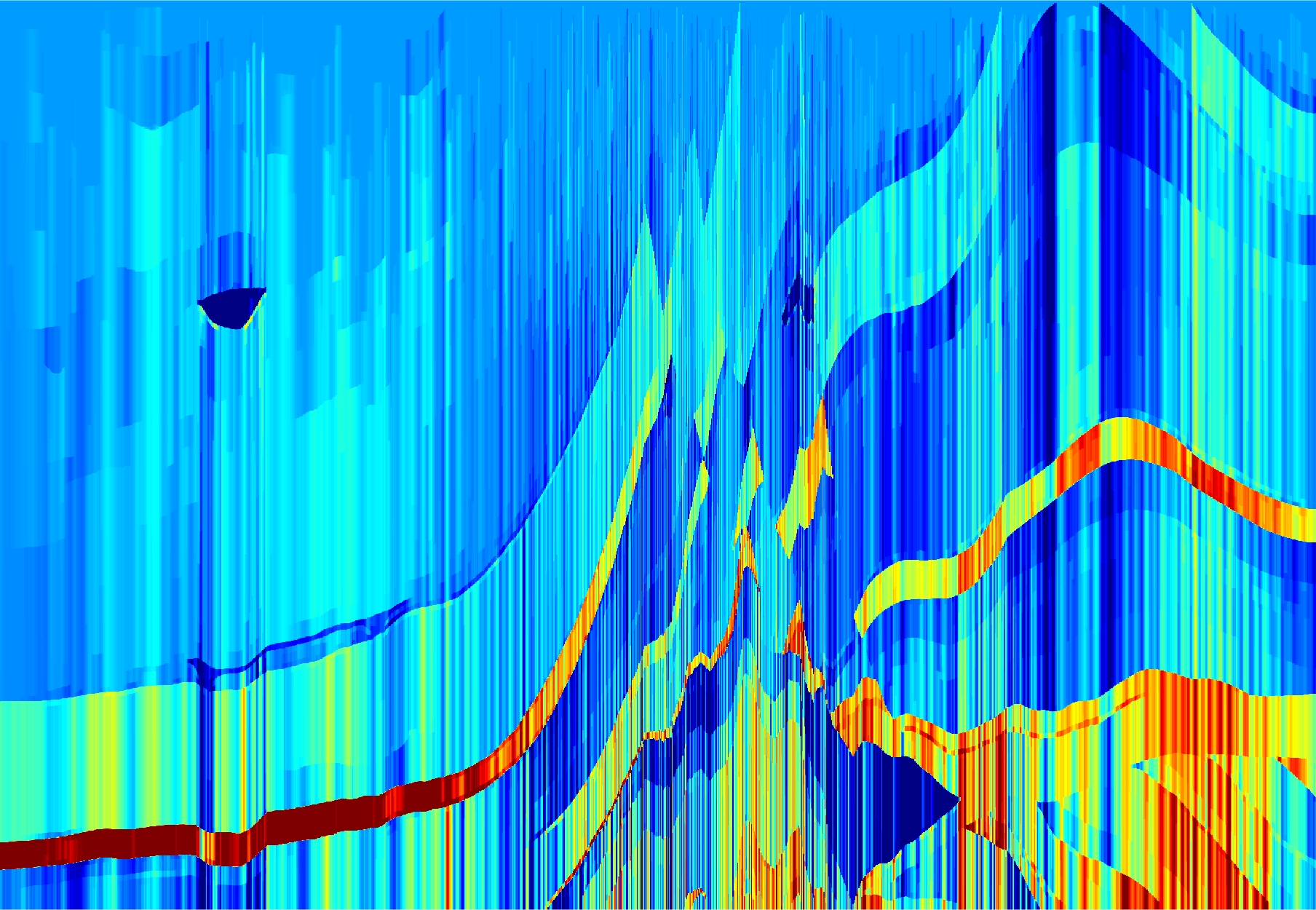}}\ \ 
    \subfloat[]{\includegraphics[height=0.22\textwidth]{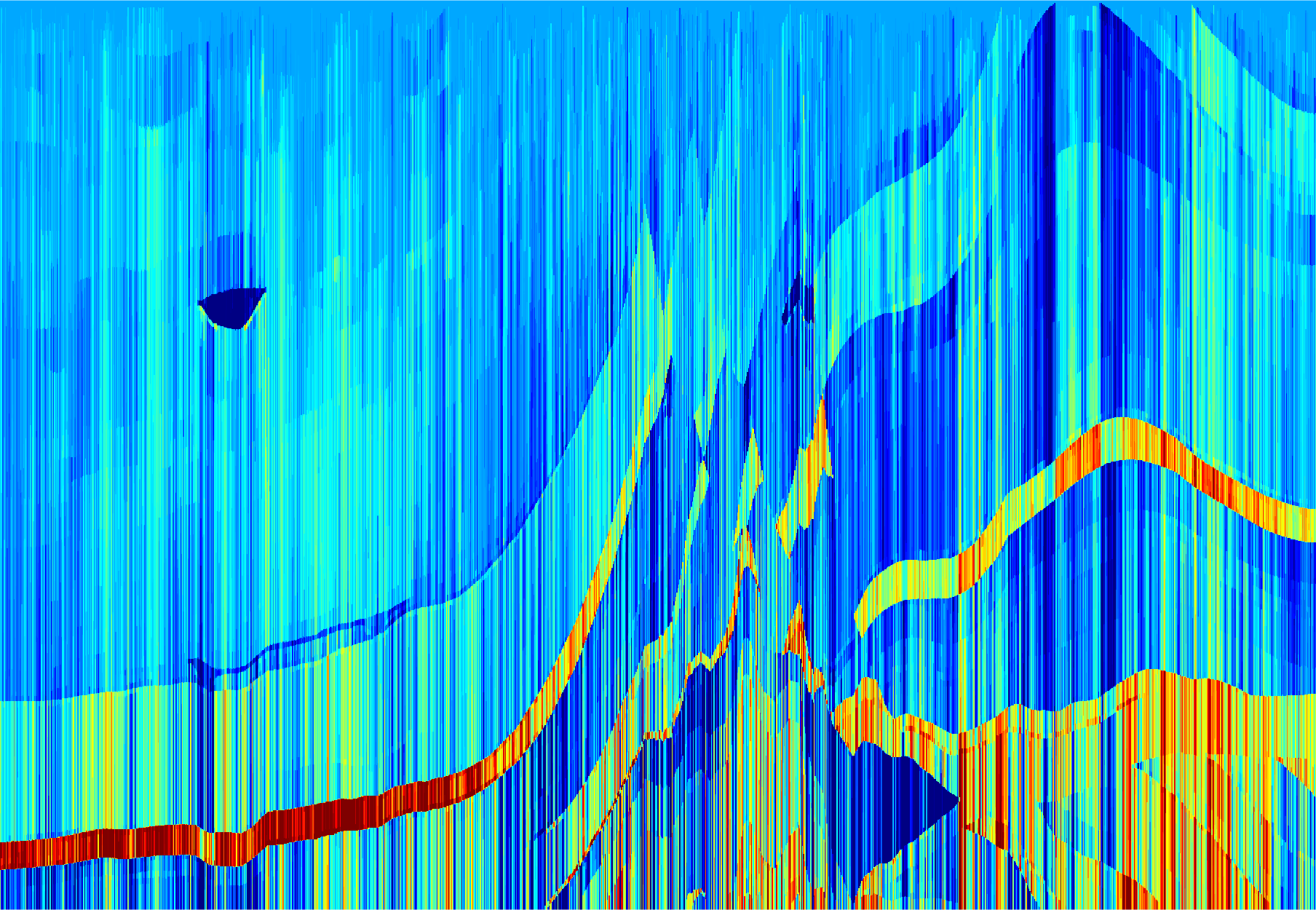}}\ \ 
    \subfloat[]{\includegraphics[height=0.22\textwidth]{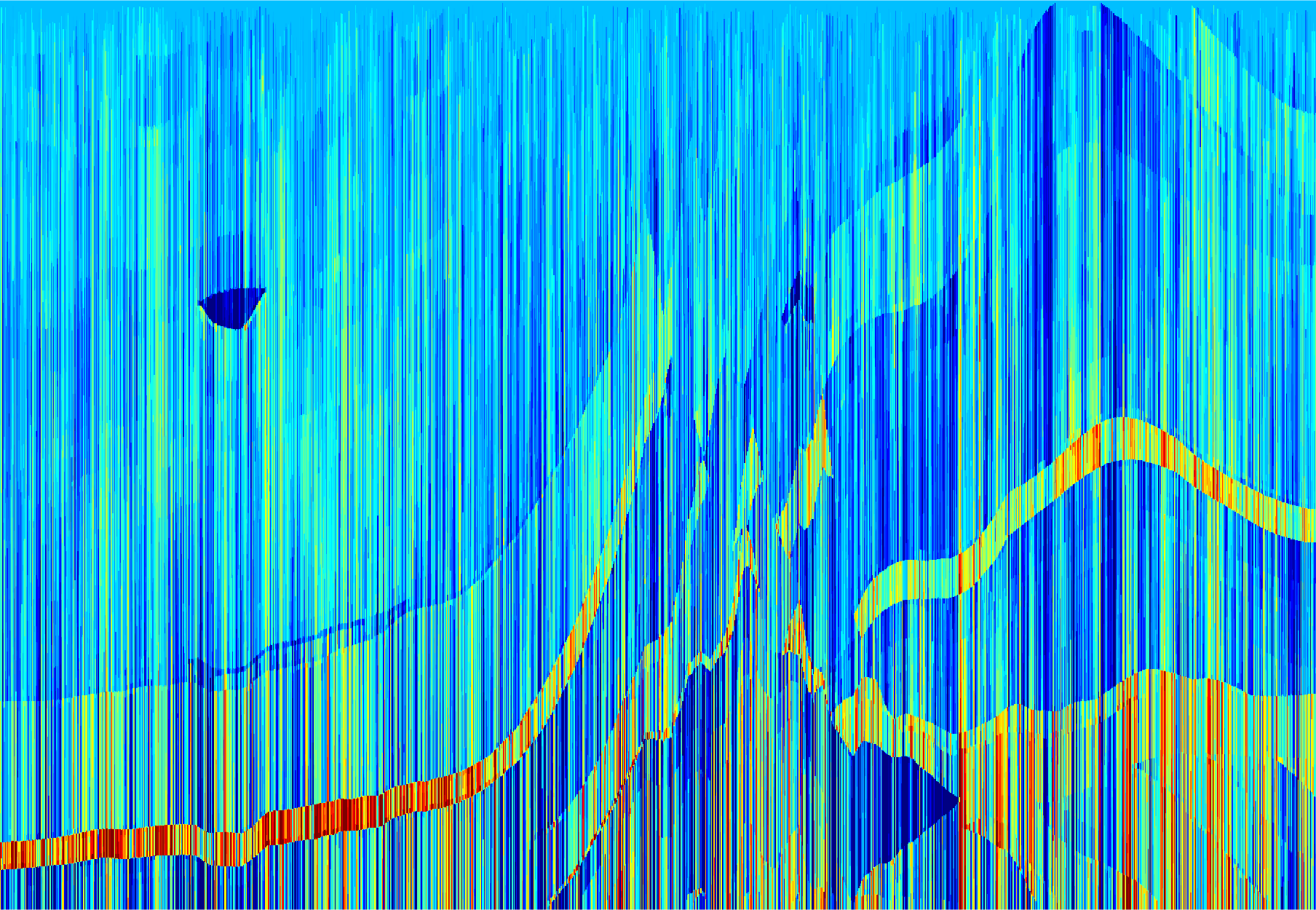}}\\
    \subfloat[]{\includegraphics[height=0.22\textwidth]{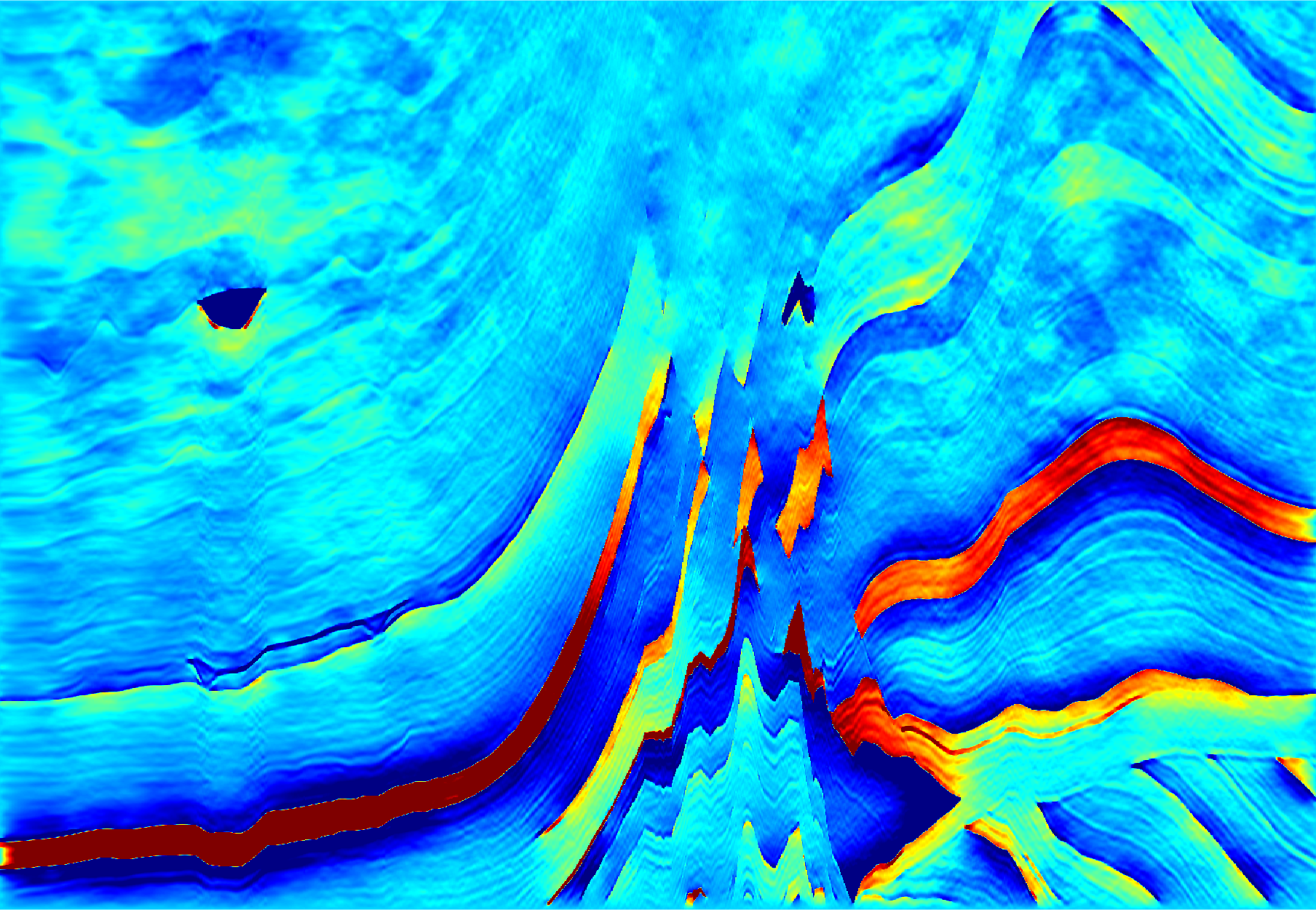}}\ \ 
    \subfloat[]{\includegraphics[height=0.22\textwidth]{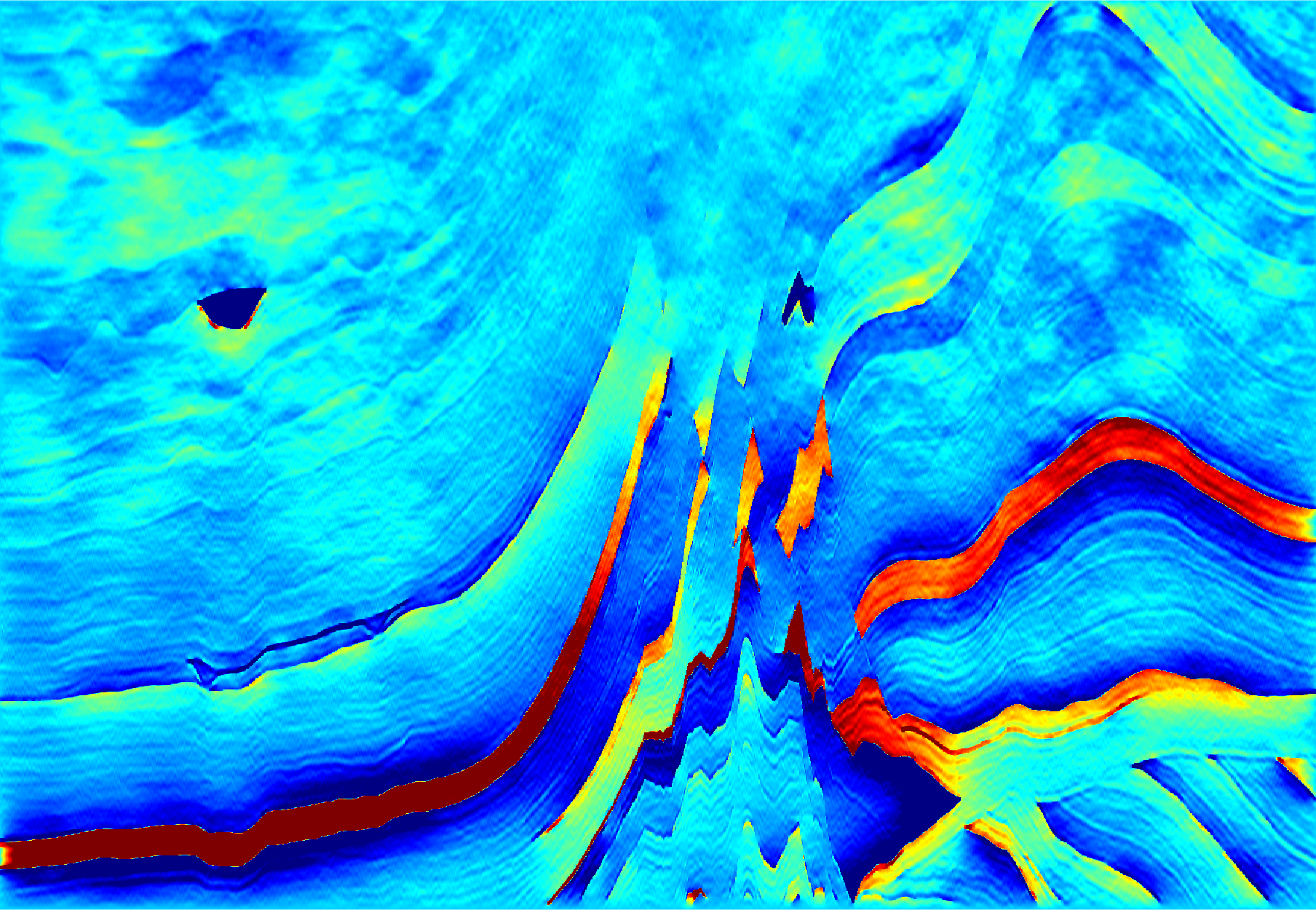}}\ \ 
    \subfloat[]{\includegraphics[height=0.22\textwidth]{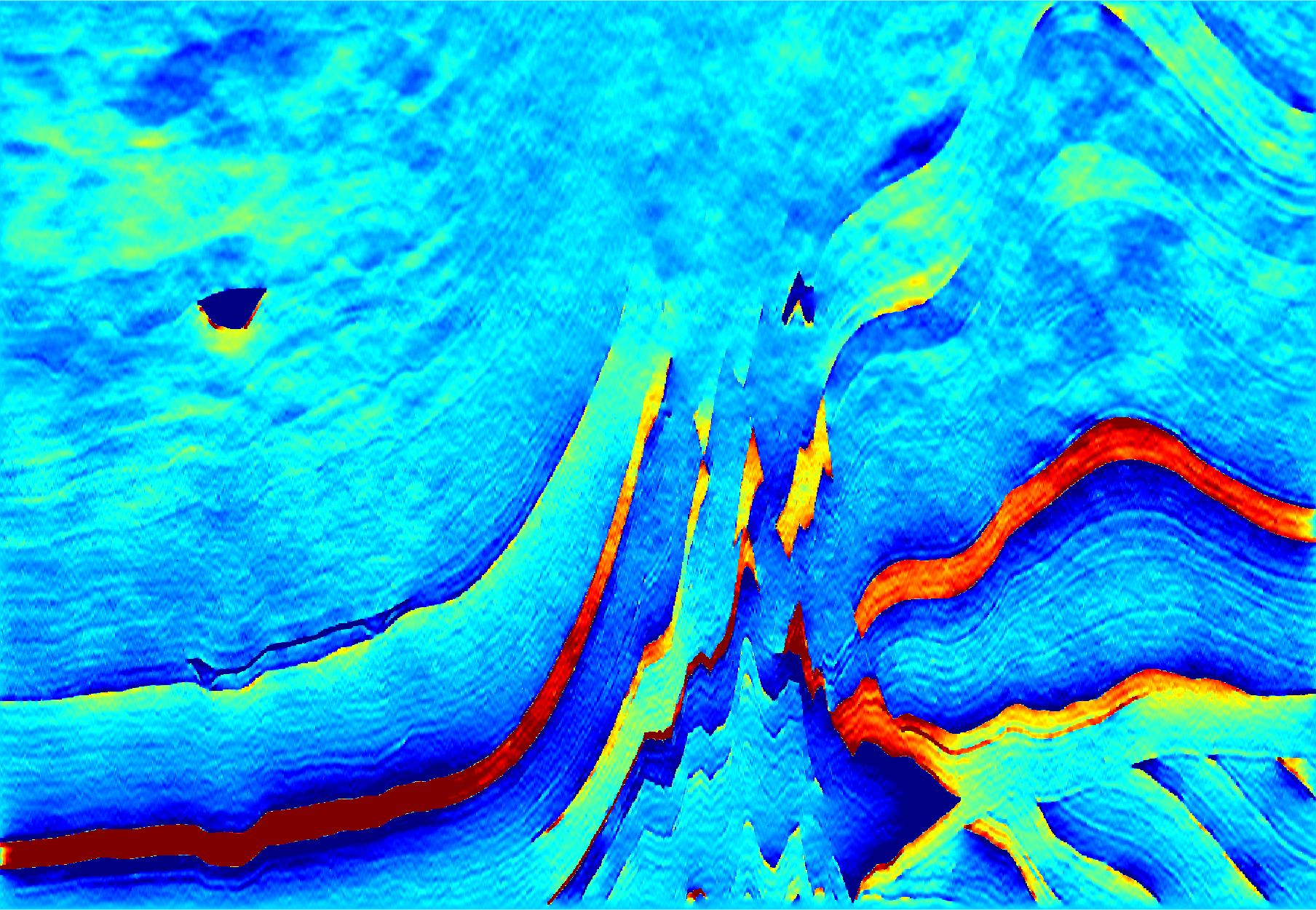}}\\
    \marmousiCBH{0.33\textwidth}
	\caption{Reconstructed impedance profiles for Marmousi2 from different initialization methods and different noise levels: AA (a,b,c), Liu (d,e,f), SSI (g,h,i), SB (j,k,l); without noise (left column: a,d,g,j), medium noise (middle column: b,e,h,k), high noise (right column: c,f,i,l).}
    \label{fig:Marmousi_init}
\end{figure}

\begin{figure}[h!tb]
	\centering
    \subfloat[]{\includegraphics[height=0.22\textwidth]{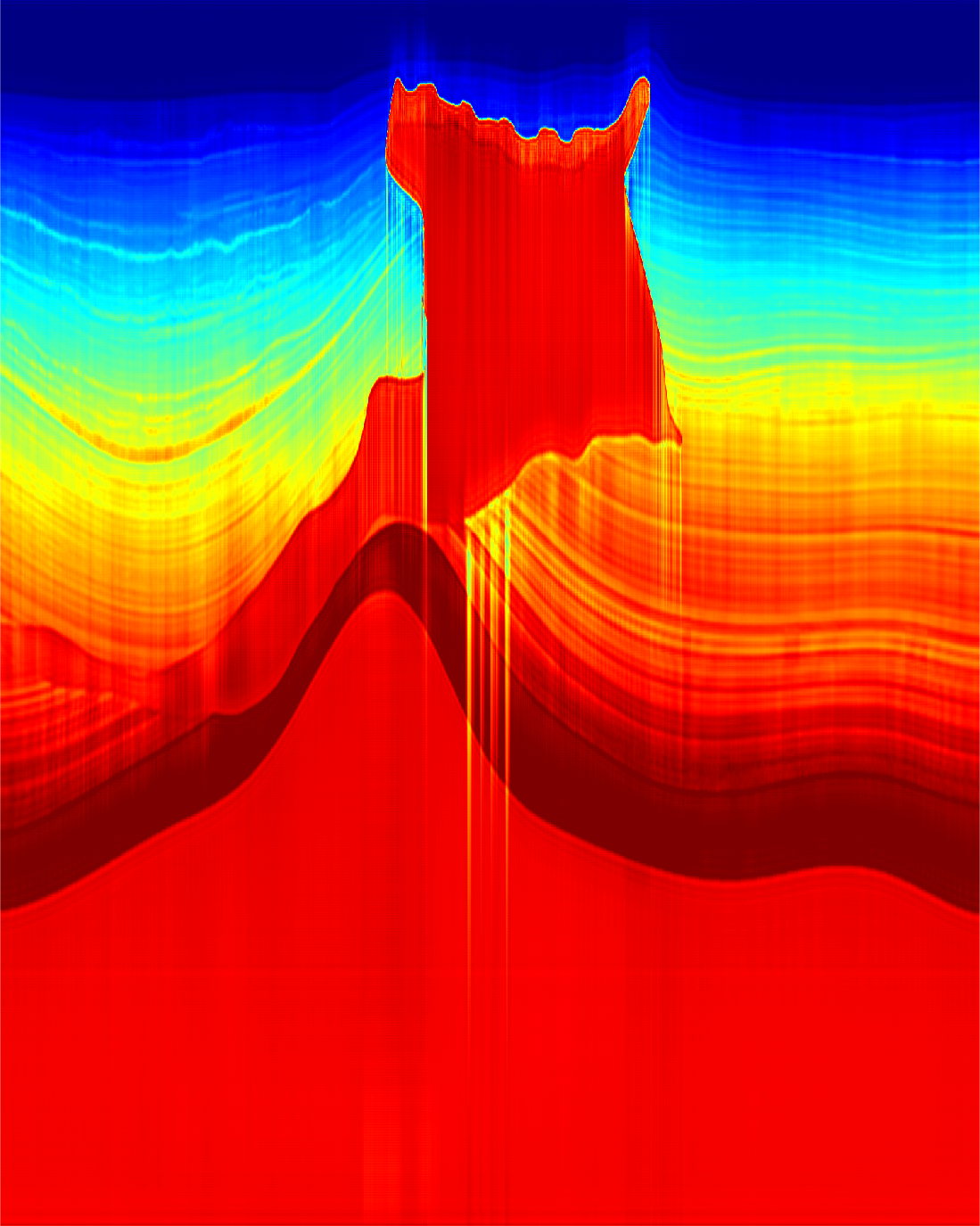}}\ \ 
    \subfloat[]{\includegraphics[height=0.22\textwidth]{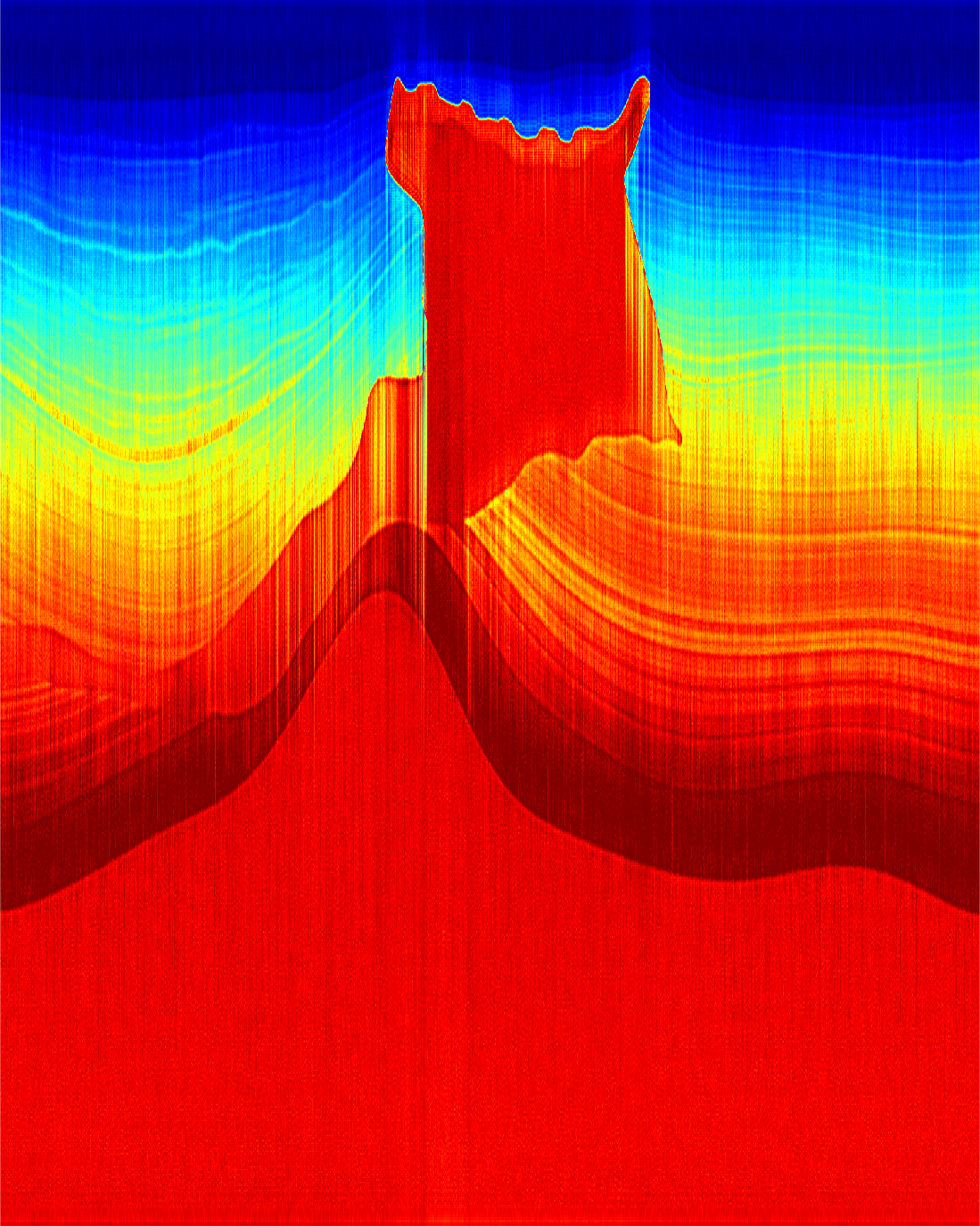}}\ \ 
    \subfloat[]{\includegraphics[height=0.22\textwidth]{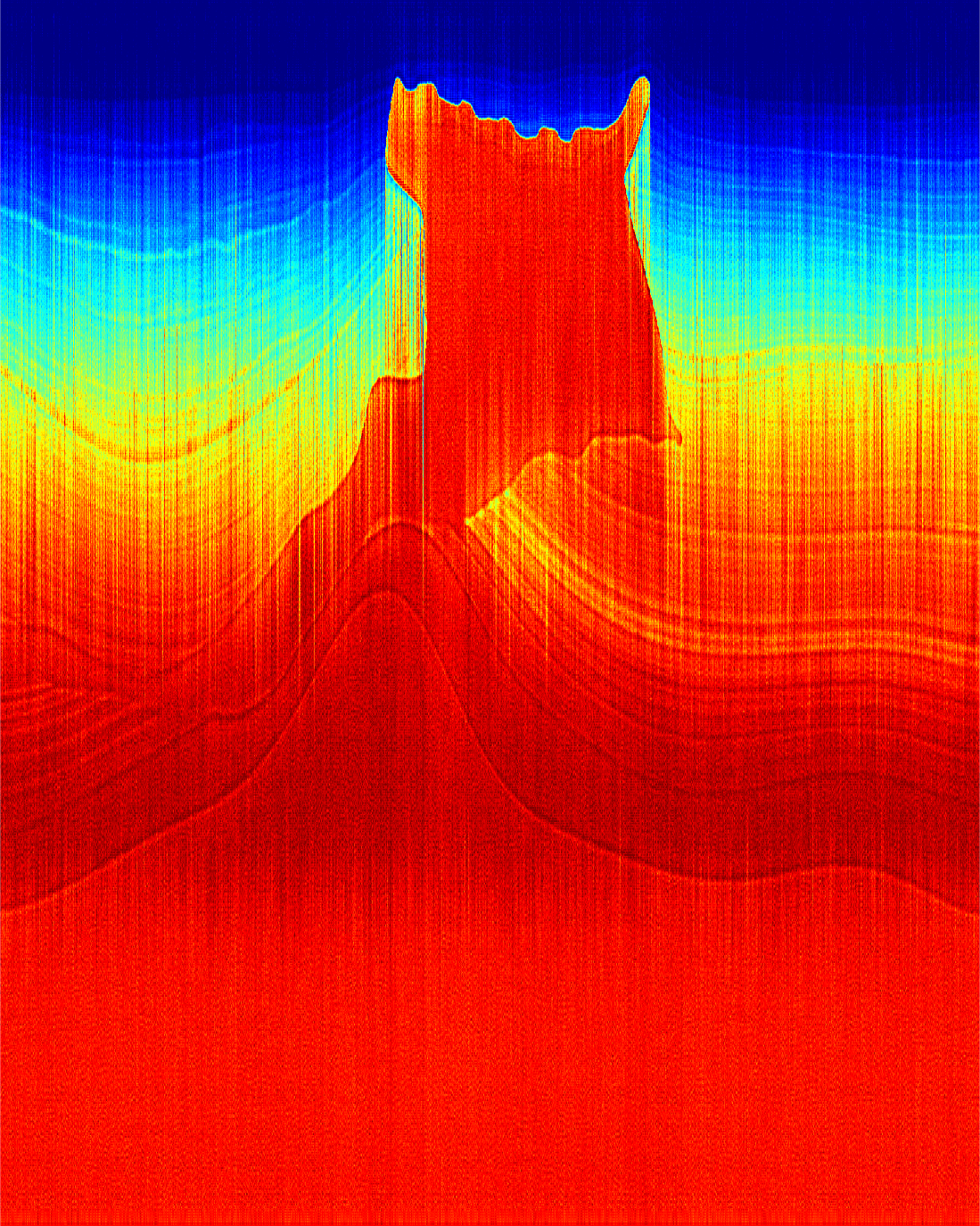}}\\
    \subfloat[]{\includegraphics[height=0.22\textwidth]{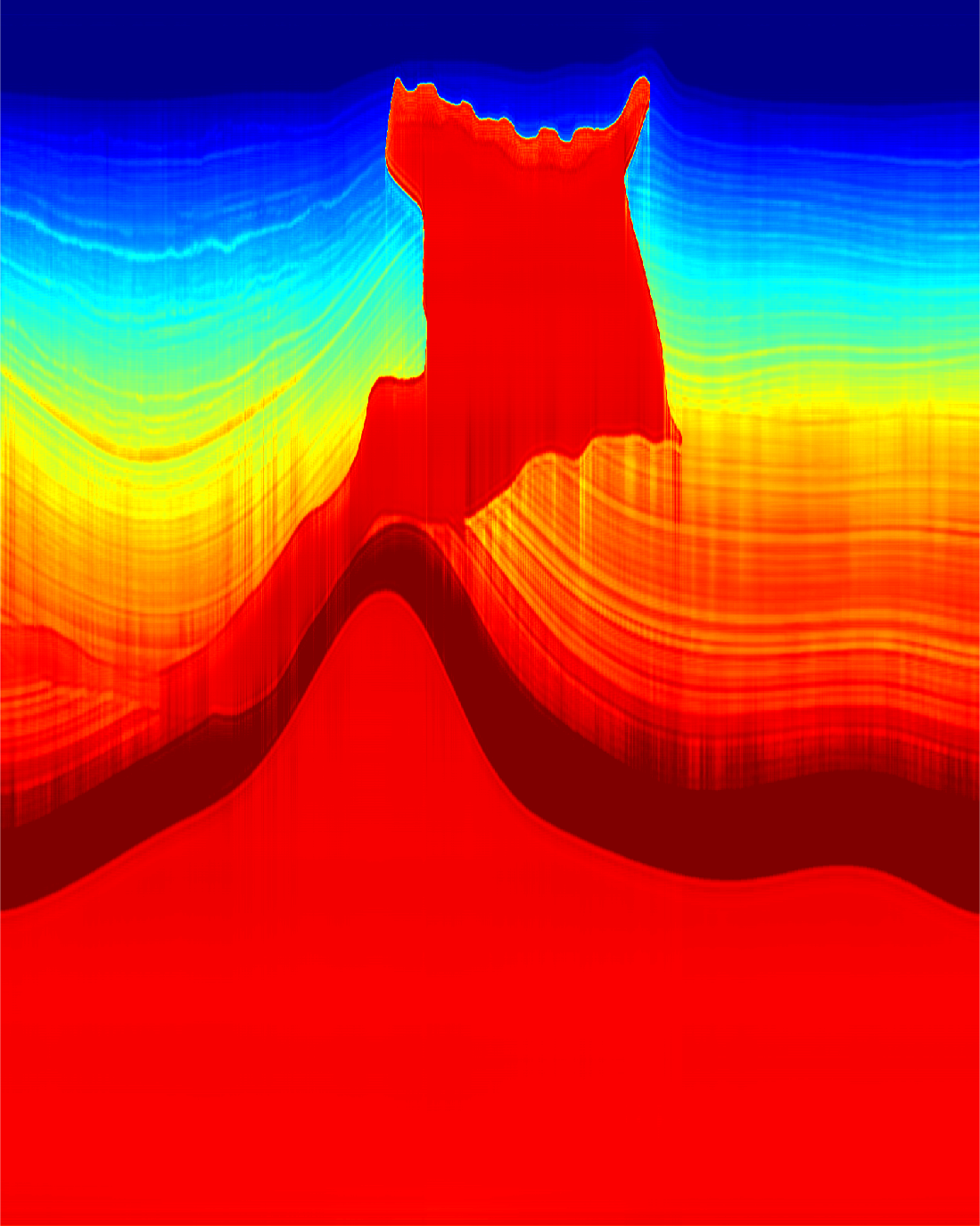}}\ \ 
    \subfloat[]{\includegraphics[height=0.22\textwidth]{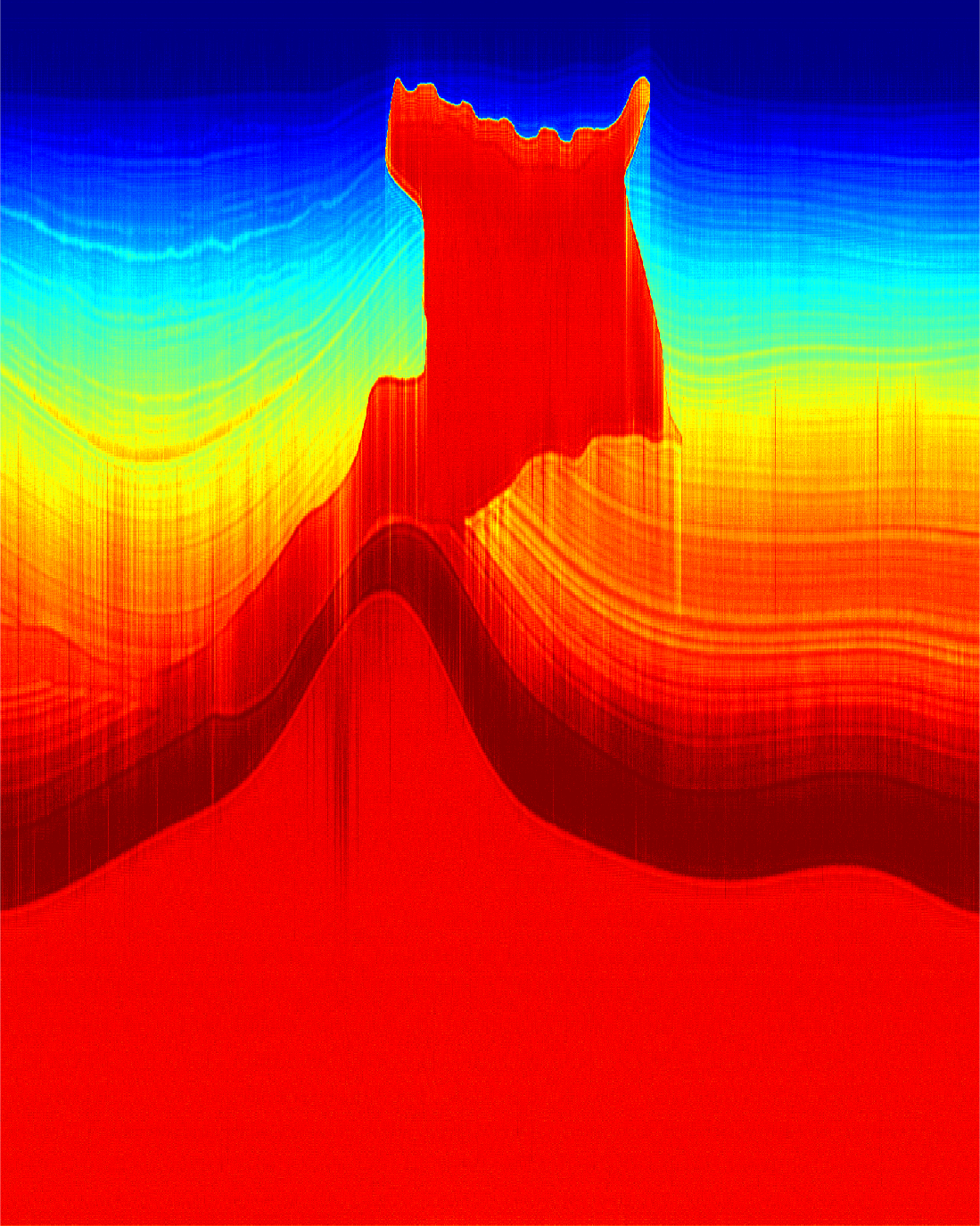}}\ \ 
    \subfloat[]{\includegraphics[height=0.22\textwidth]{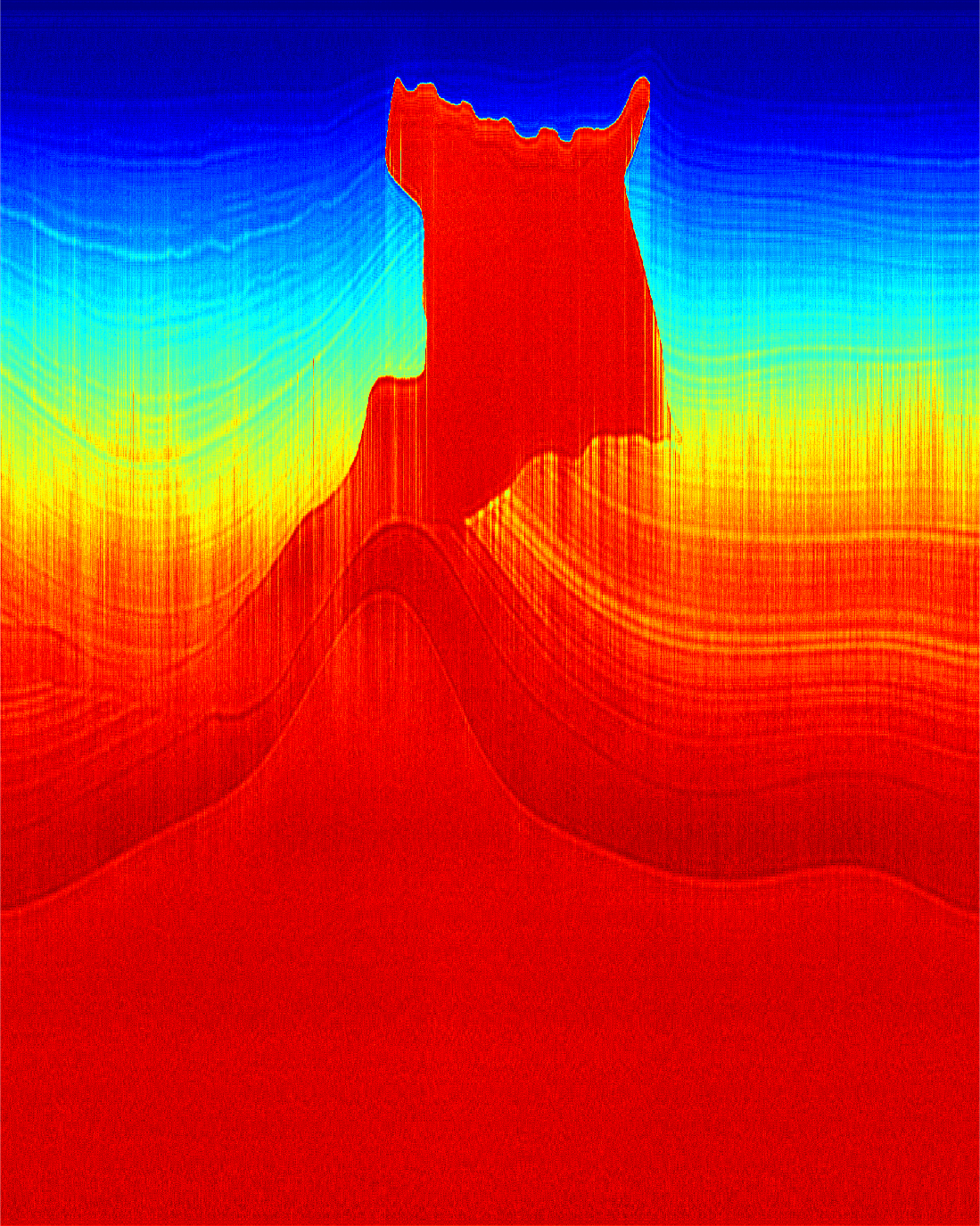}}\\
    \subfloat[]{\includegraphics[height=0.22\textwidth]{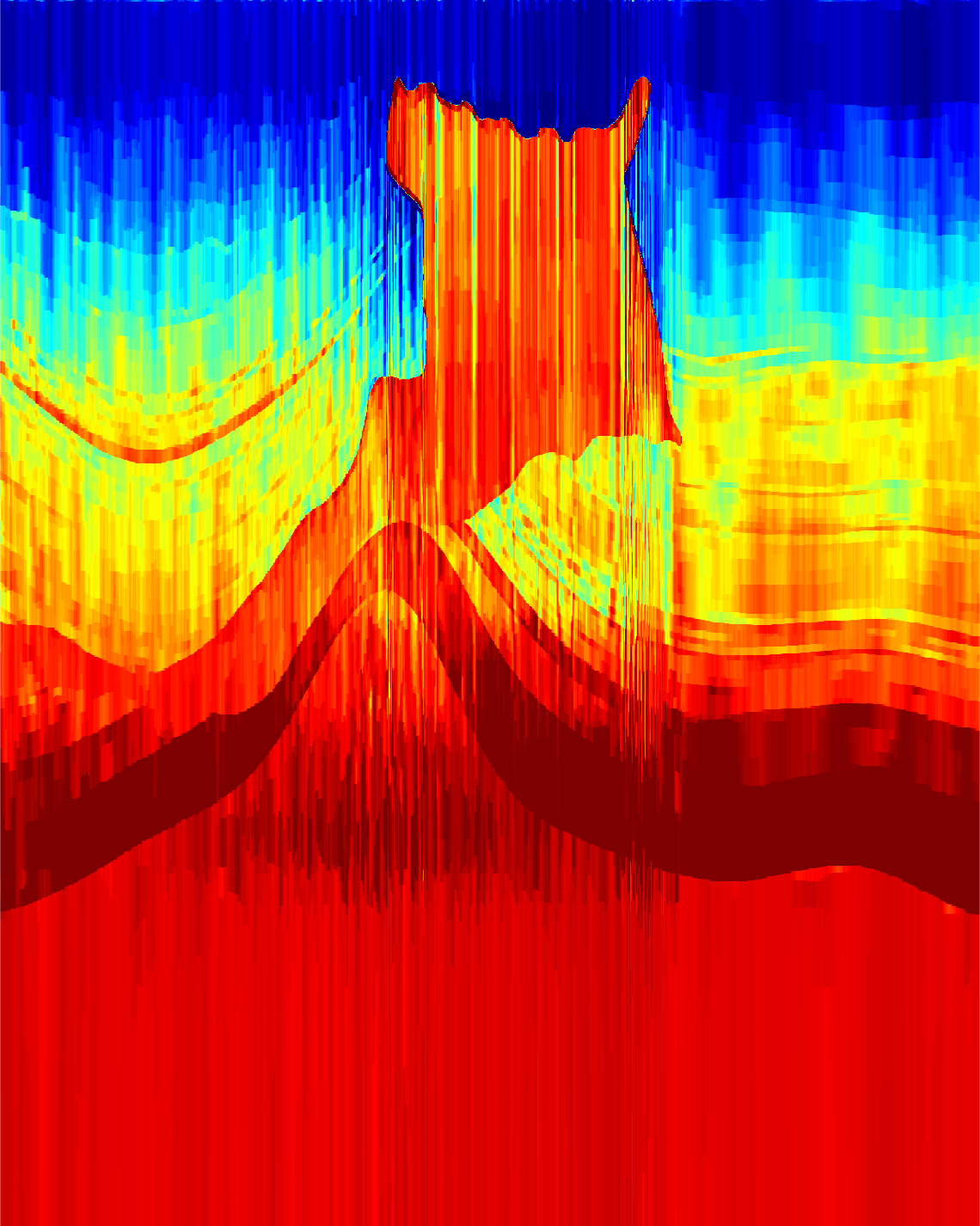}}\ \ 
    \subfloat[]{\includegraphics[height=0.22\textwidth]{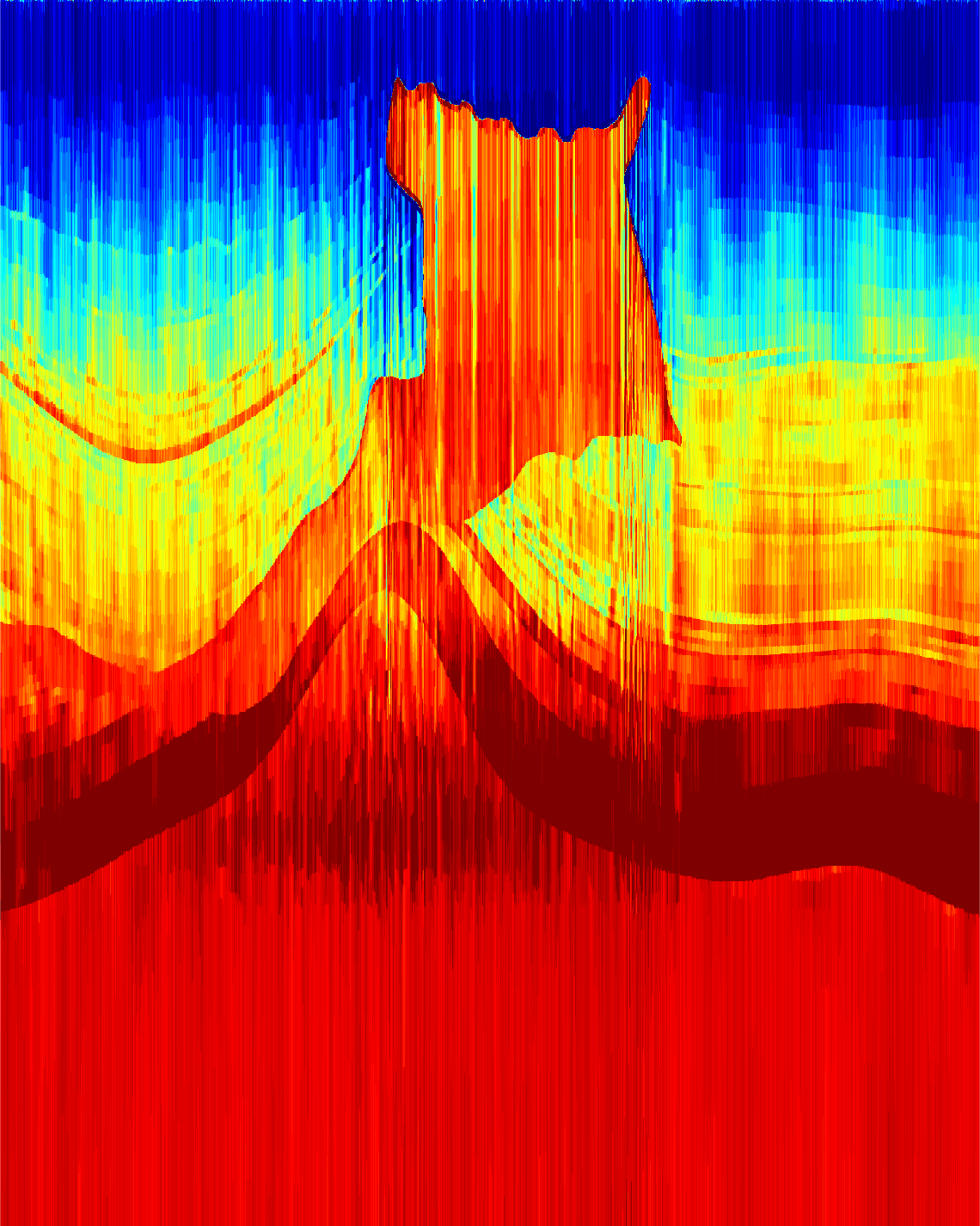}}\ \ 
    \subfloat[]{\includegraphics[height=0.22\textwidth]{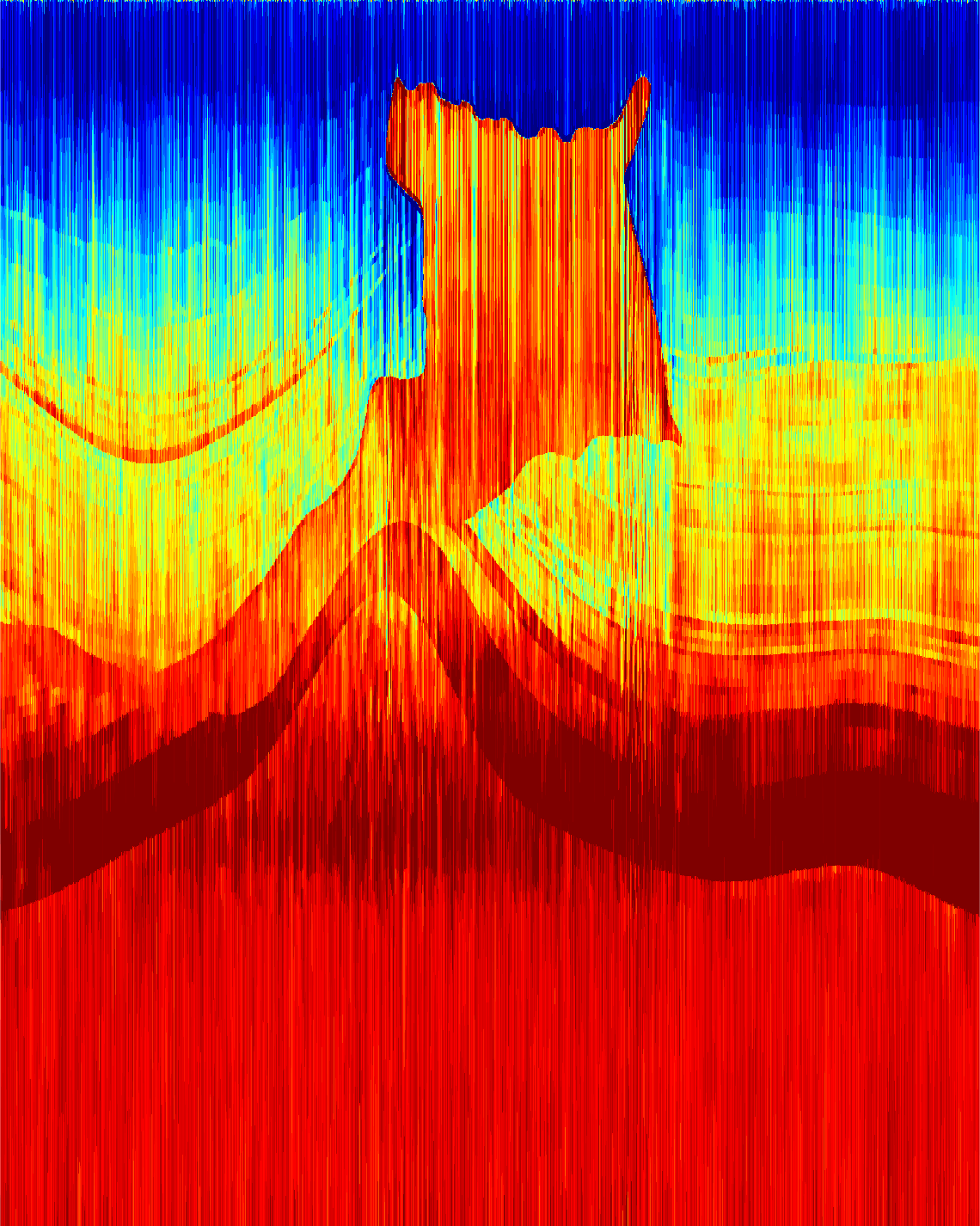}}\\
    \subfloat[]{\includegraphics[height=0.22\textwidth]{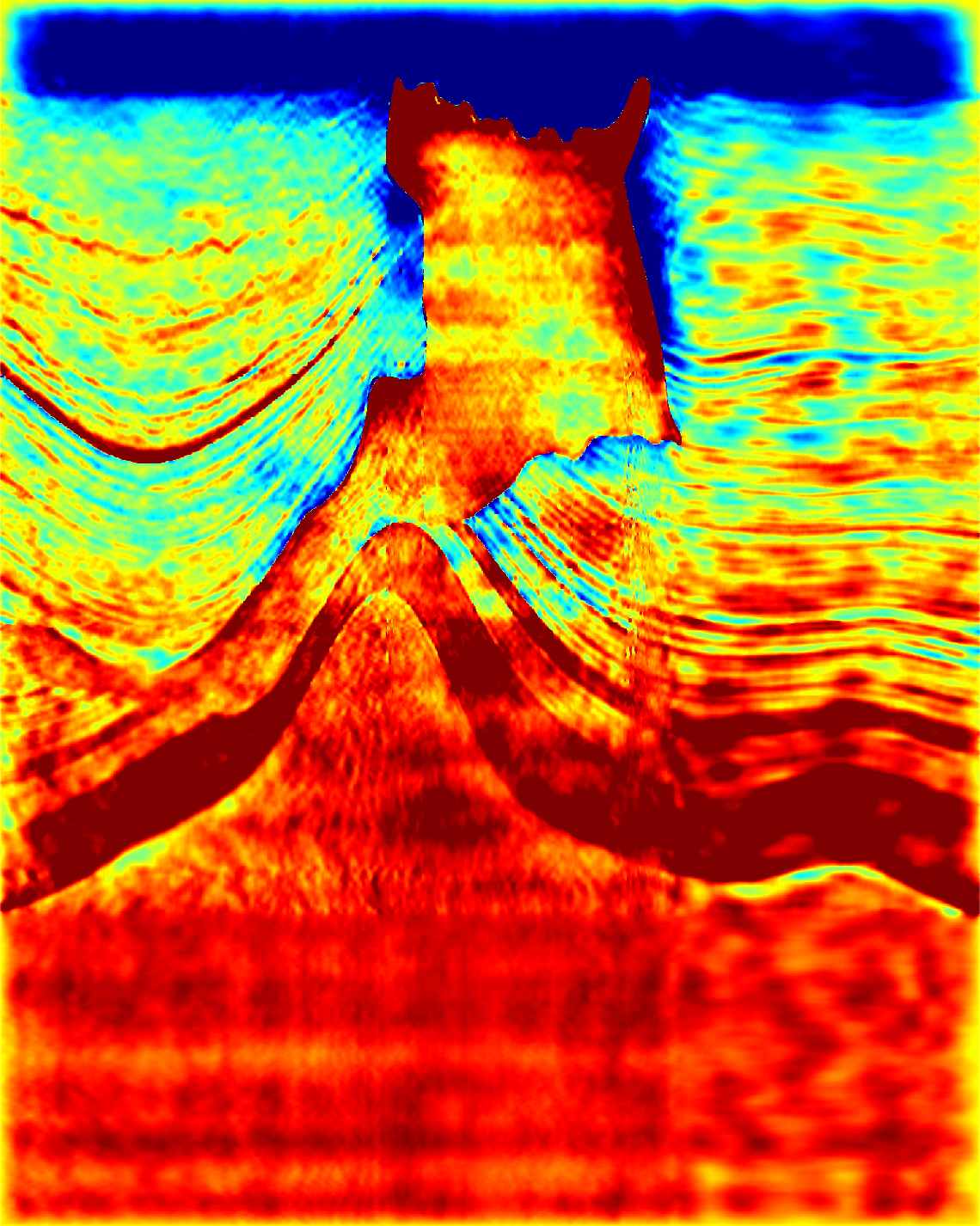}}\ \ 
    \subfloat[]{\includegraphics[height=0.22\textwidth]{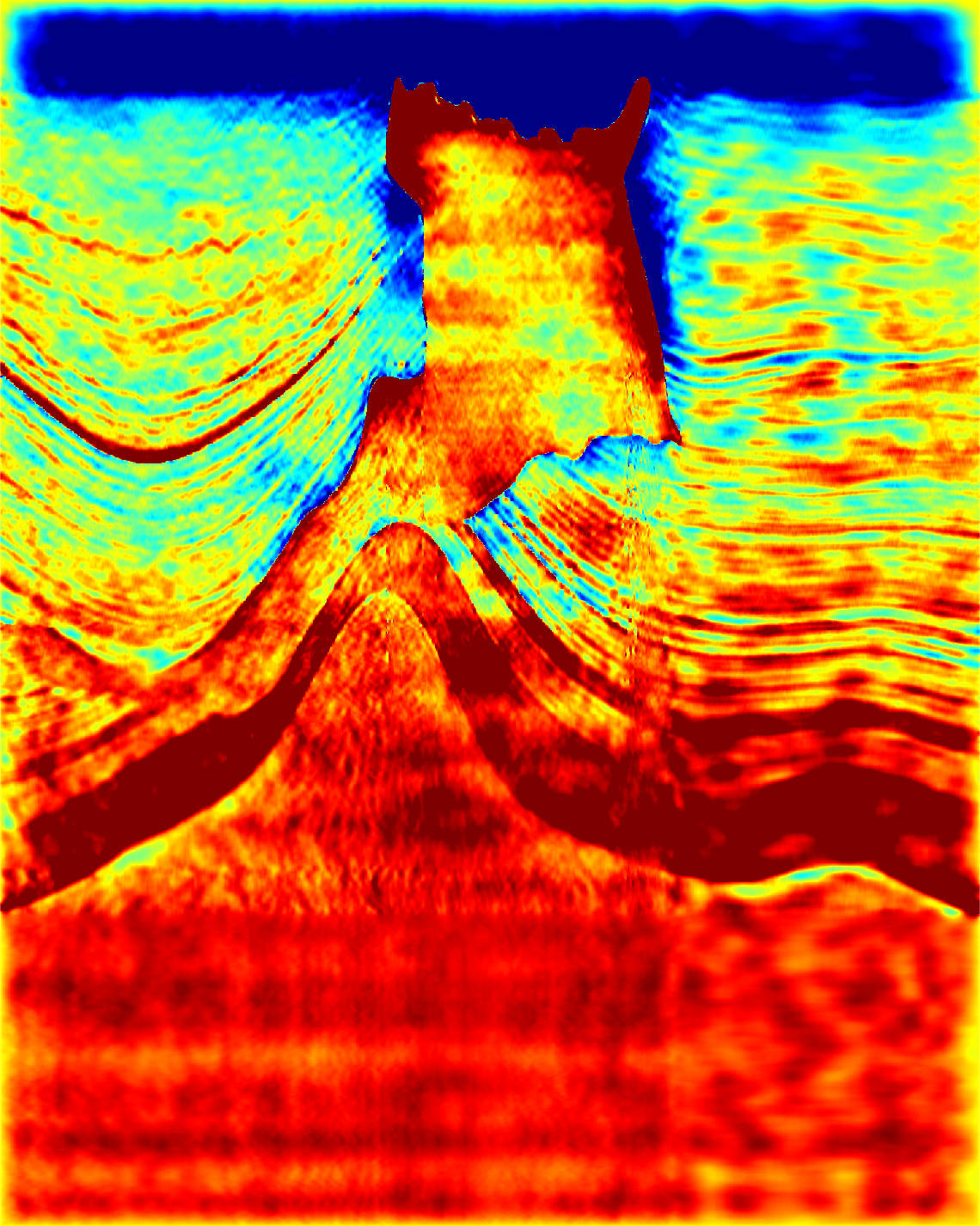}}\ \ 
    \subfloat[]{\includegraphics[height=0.22\textwidth]{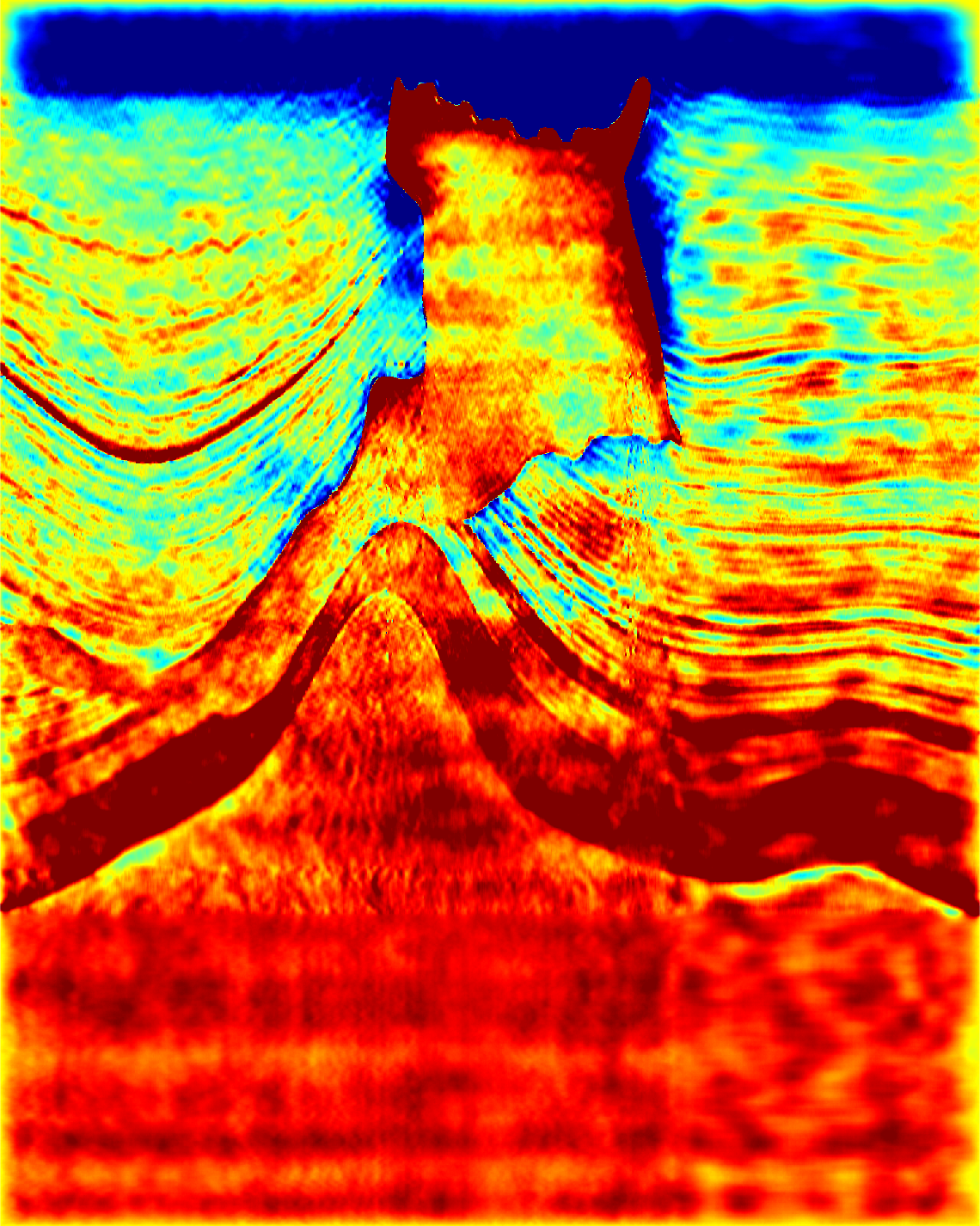}}\\
    \seamCBH{0.33\textwidth}
	\caption{Reconstructed impedance profiles for SEAM from different initialization methods and different noise levels: AA (a,b,c), Liu (d,e,f), SSI (g,h,i), SB (j,k,l); without noise (left column: a,d,g,j), medium noise (middle column: b,e,h,k), high noise (right column: c,f,i,l).}
    \label{fig:SEAM_init}
\end{figure}

We use the four initialization methods discussed earlier to obtain a starting guess $\bx_0^\delta$. Because of the highly underdetermined system and the noise, we obtained better results using large regularization parameters in equation \ref{eq:classical_model}. For SSI we set $\alpha=15$ for Marmousi2 and $\alpha=7.5$ for SEAM; SB worked best with $\alpha=\beta=40$ for both models. The reconstruction results for all methods on different noise levels are shown in Figure~\ref{fig:Marmousi_init} (Marmousi2) and Figure~\ref{fig:SEAM_init} (SEAM). For a better comparison we use the exact same color-coding as in Figure \ref{fig:data} and normalize the reconstructed impedance profiles to fit the colormap. The trace-wise reconstruction methods AA and SSI are more prone to noise. In addition, both classical methods tend to over-smooth the reconstruction due to the large regularization parameters.

\begin{figure}[h!tb]
	\centering
    \begin{tabular}{lll}
    \subfloat[]{\includegraphics[height=0.22\textwidth]{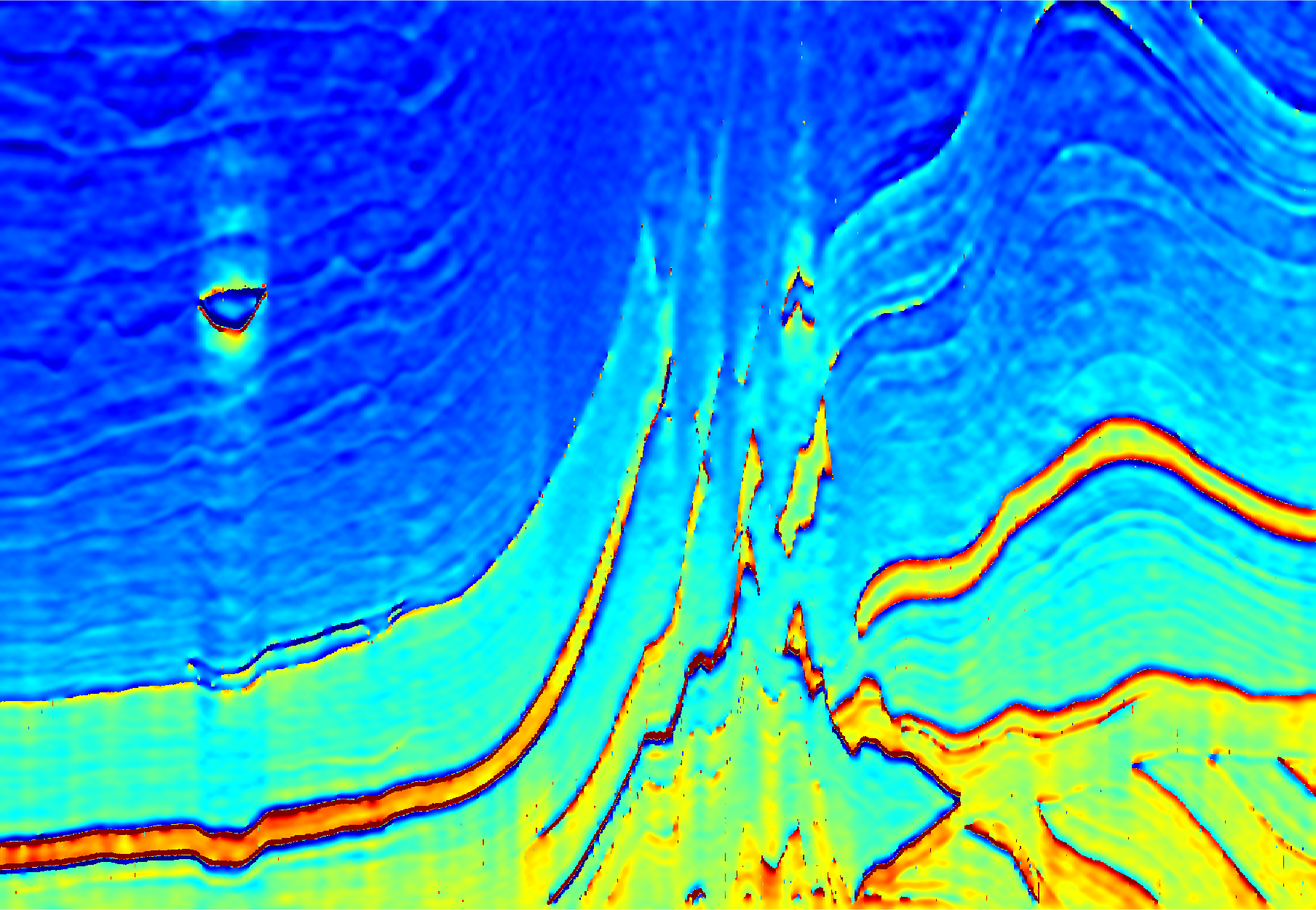}} &
    \subfloat[]{\includegraphics[height=0.22\textwidth]{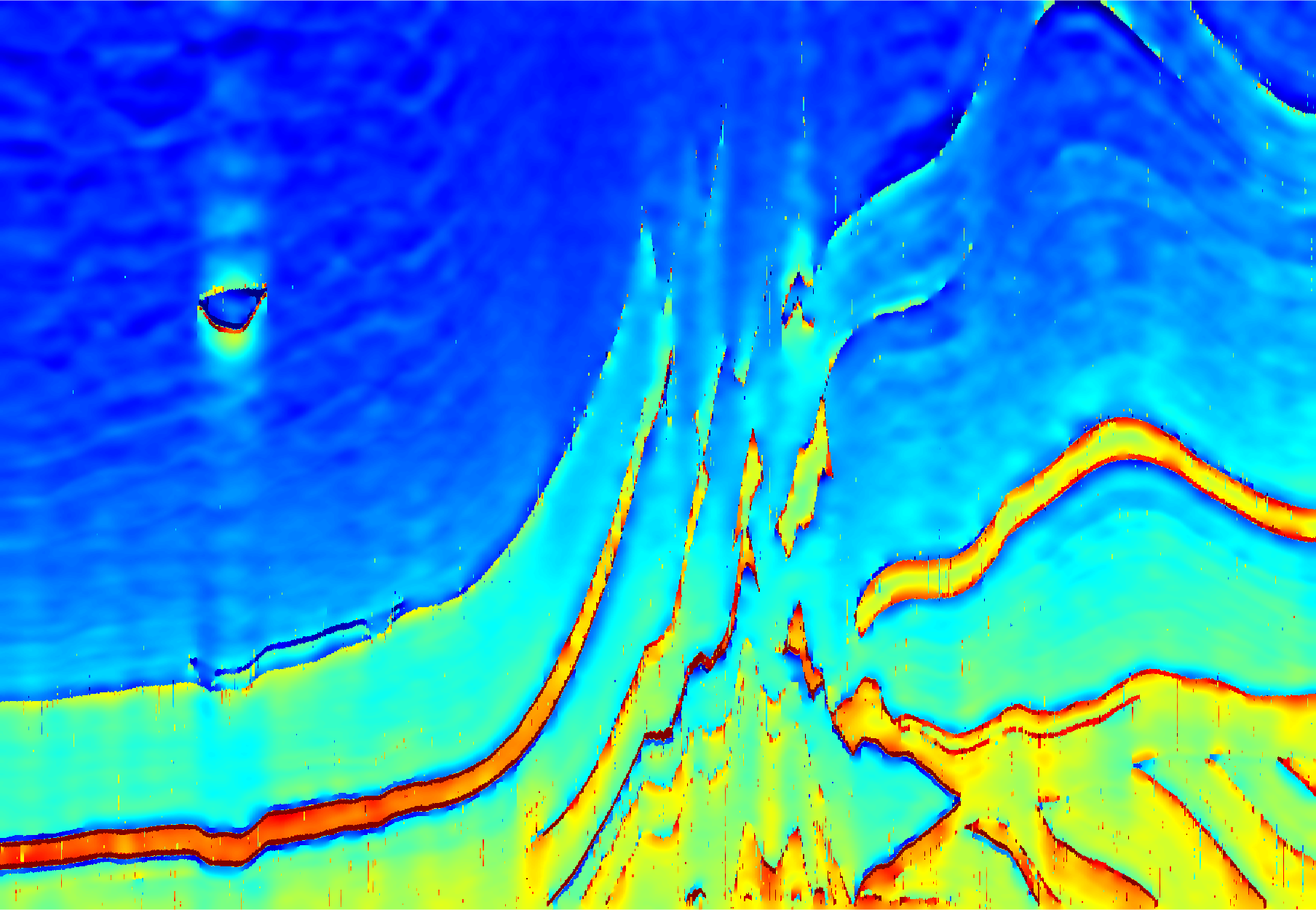}} &
    \marmousiCB{0.22\textwidth} \\
    \subfloat[]{\includegraphics[height=0.22\textwidth]{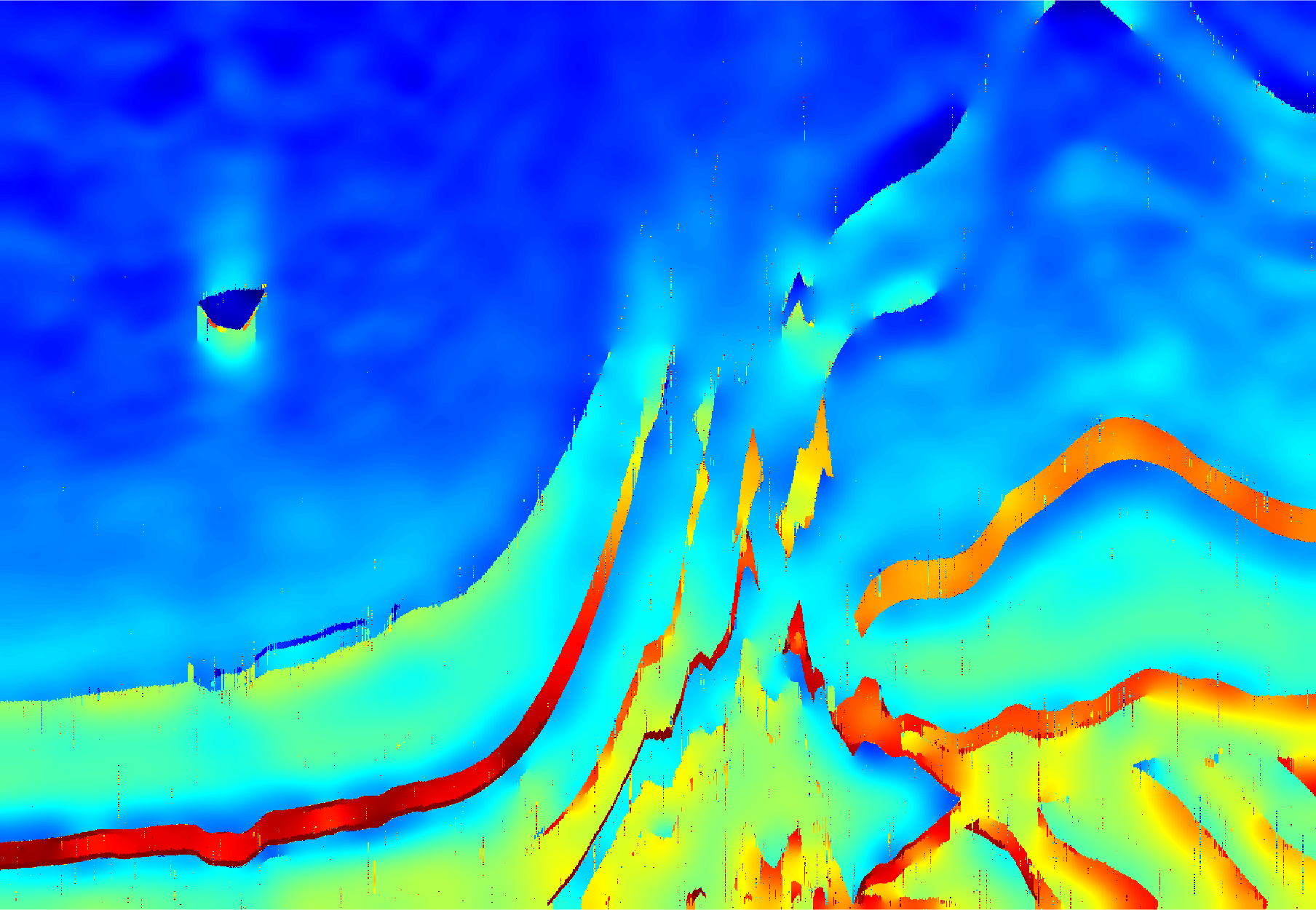}} & 
    \subfloat[]{\includegraphics[height=0.22\textwidth]{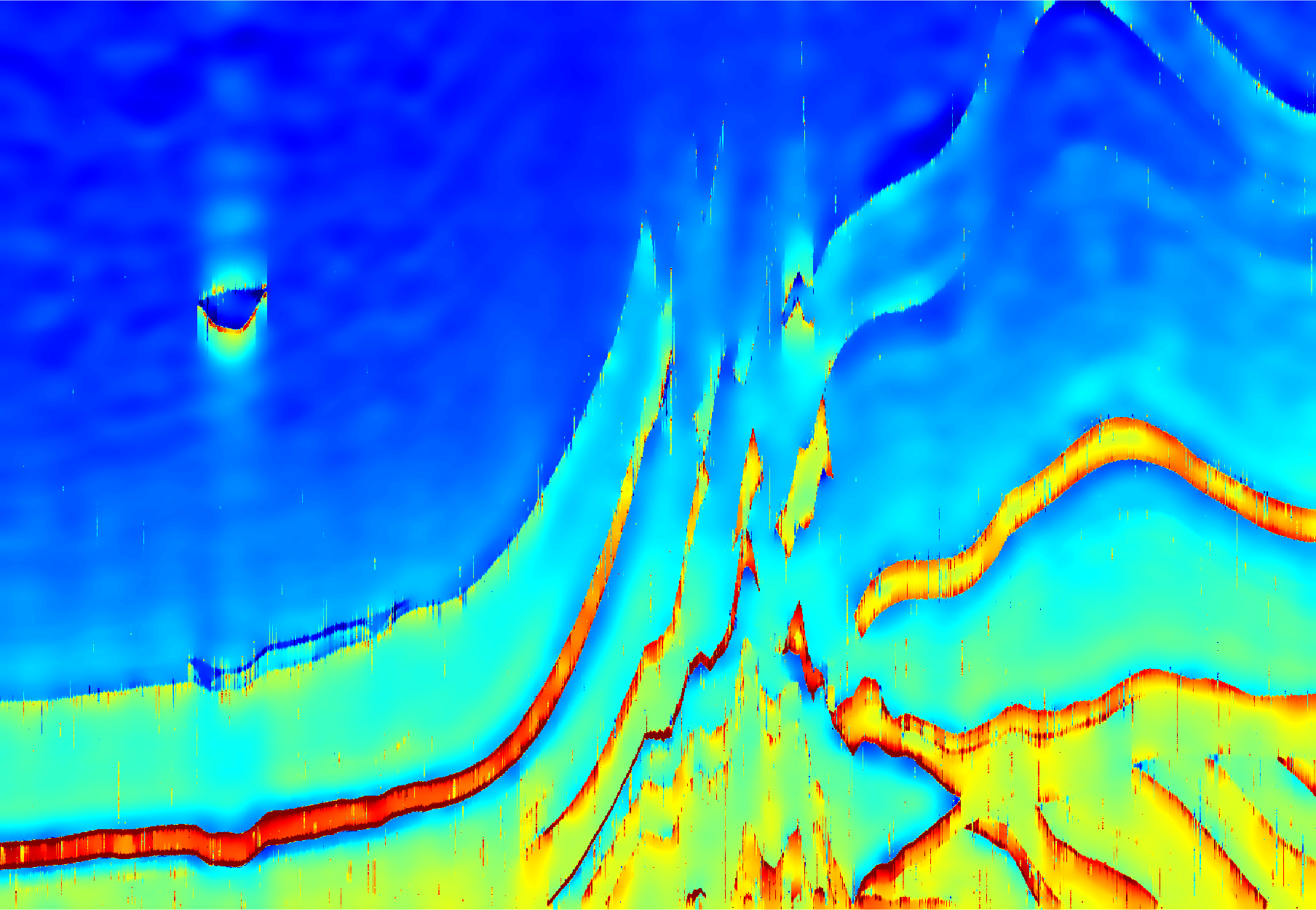}}
    \end{tabular}
	\caption{Reconstructed impedance profiles for Marmousi2 with medium noise, AA initialization, and different radii $R$ for \ref{iterated_graphLaNet}: $R=2$ (a), $R=3$ (b), $R=7$ (c), and $R=7$ after 2 iterations (d). Larger radii converge faster but tend to remove smaller structures and are prone to overregularization.}
    \label{fig:RcomparisonRecs}
\end{figure}

To start the iterative process \ref{iterated_graphLaNet} we first need to find a suitable parameters $\sigma$ and $R$. We performed multiple experiments to test the influence of $\sigma$ on the reconstruction. As it turns out, this parameter is not too relevant for the result and a choice of $\sigma\in[0.1,1]$ performs reasonable well as long as the initial guess $\bx^\delta_0$ is normalized for the calculation of the graph Laplacian matrix. Hence, for the following experiment we choose the optimal $\sigma$ from a list of only four values $\sigma=1,0.5,0.25,0.1$. The choice of the radius $R$ is more influential on the result. Generally, a larger radius tends to return better results after only a few iterations. However, it is prone to oversmoothing effects and can completely remove small structures and details if iterated too long.  In Figure~\ref{fig:Rcomparison}, we show the approximation error as a function of the number of iterations for medium-noise-level Marmousi data, using AA as the initialization method for different radii $R = 2, 3, 7$ with corresponding values $\sigma = 0.25, 1, 1$. For $R=7$ the best SSIM value is obtained after 2 iterations but the method becomes unstable afterwards. $R=2$ shows the slowest convergence in the first steps but proves to be the most stable setup in the end. Hence, we recommend using larger radii $R$ only if a reliable stopping criterion is known. As illustration, we show all obtained reconstructions after 10 iterations as well as the reconstruction for $R=7$ after two iterations in Figure \ref{fig:RcomparisonRecs} (also compare to Figure \ref{fig:Marmousi_init} (b) which shows the initial reconstruction used). For a better visual comparison the data was normalized and color-coded exactly the same.

\begin{figure}[h!tb]
	\centering
    \subfloat[]{
    \begin{tikzpicture}
		\begin{axis}[enlargelimits=false,
            x tick label style = {font = \footnotesize},
            y tick label style = {font = \footnotesize},
            xlabel = {Iterations}, ylabel = {D-MSE},
			axis on top,
            width=0.45\textwidth, height=0.45\textwidth]
			\addplot[blue,line width=2pt] coordinates {
            (0,0.00917361) (1,0.00458998) (2,0.00448352) (3,0.00469454) (4,0.00491143) (5,0.00507608) (6,0.00518885) (7,0.00527658) (8,0.0053424) (9,0.00540029) (10,0.00544803)
            };
            \addplot[red,dashed,line width=2pt] coordinates {
            (0,0.00917361) (1,0.0044101) (2,0.00525068) (3,0.00761923) (4,0.0100452) (5,0.0118751) (6,0.0131511) (7,0.0140283) (8,0.0146721) (9,0.0151562) (10,0.0154973)
            };
            \addplot[green!50!black,dotted,line width=2pt] coordinates {
            (0,0.00917361) (1,0.00513312) (2,0.00467862) (3,0.00461918) (4,0.0047393) (5,0.00486511) (6,0.00494926) (7,0.00499905) (8,0.00503687) (9,0.00506377) (10,0.00508606)
            };
		\end{axis}
	\end{tikzpicture}
    }\ \ 
    \subfloat[]{
    \begin{tikzpicture}
		\begin{axis}[enlargelimits=false,
            x tick label style = {font = \footnotesize},
            y tick label style = {font = \footnotesize},
            xlabel = {Iterations}, ylabel = {SSIM},
			axis on top,
            width=0.45\textwidth, height=0.45\textwidth]
			\addplot[blue,line width=2pt] coordinates {
            (0,0.379016) (1,0.670252) (2,0.719356) (3,0.731151) (4,0.735411) (5,0.737351) (6,0.738415) (7,0.739052) (8,0.739419) (9,0.739575) (10,0.739556)
            };
            \addplot[red,dashed,line width=2pt] coordinates {
            (0,0.379016) (1,0.704362) (2,0.725262) (3,0.724326) (4,0.720861) (5,0.717068) (6,0.713433) (7,0.709347) (8,0.704731) (9,0.699895) (10,0.695067)
            };
            \addplot[green!50!black,dotted,line width=2pt] coordinates {
            (0,0.379016) (1,0.630062) (2,0.693433) (3,0.714836) (4,0.725183) (5,0.731031) (6,0.734418) (7,0.736479) (8,0.737654) (9,0.738382) (10,0.738769)
            };
		\end{axis}
	\end{tikzpicture}
    }
	\caption{D-MSE (a) and SSIM (b) for Marmousi2 with medium noise and AA initialization against number of \ref{iterated_graphLaNet} iterations with different radius $R$: $R=7$ (red, dashed), $R=3$ (blue), and $R=2$ (green, dotted). A larger radius leads to a more optimal but also unstable reconstruction with overregularization after some iterations.}
    \label{fig:Rcomparison}
\end{figure}

Another important consideration for the choice of $R$ is the runtime of the algorithm. From equation \ref{eq:edge-weight-function:applications} it follows that the number of non-zero entries in $\Delta_{\bx}$ strongly depends on the choice of $R$. The algorithm performance increases the sparser $\Delta_{\bx}$ becomes. Besides $R$ the data size has the second largest influence on the measured runtime. Calculations for the Marmousi2 model with about $5.1$ million pixels took considerably longer than for the SEAM model with about $2.8$ million pixels. The average runtime in seconds per iteration is shown in Table \ref{tab:runtimes}. The runtime was very consistent with only minimal variance. The experiments were performed in Matlab 2020a with a Intel Xeon Silver 4210R CPU (2.4GHZ, 2 processors) and 128GB RAM.

\begin{table}[]
    \centering
    \begin{tabular}{|c|c|c|c|}
        \hline
         Data & $R=2$ & $R=3$ & $R=7$ \\
         Marmousi2 & 58.67 & 71.29 & 239.16 \\
         SEAM & 32.54 & 39.90 & 111.26\\
         \hline
    \end{tabular}
    \caption{Algorithm runtime in seconds per iteration for both datasets and different choices of $R$.}
    \label{tab:runtimes}
\end{table}

\begin{figure}[h!tb]
	\centering
    \subfloat[]{\includegraphics[height=0.22\textwidth]{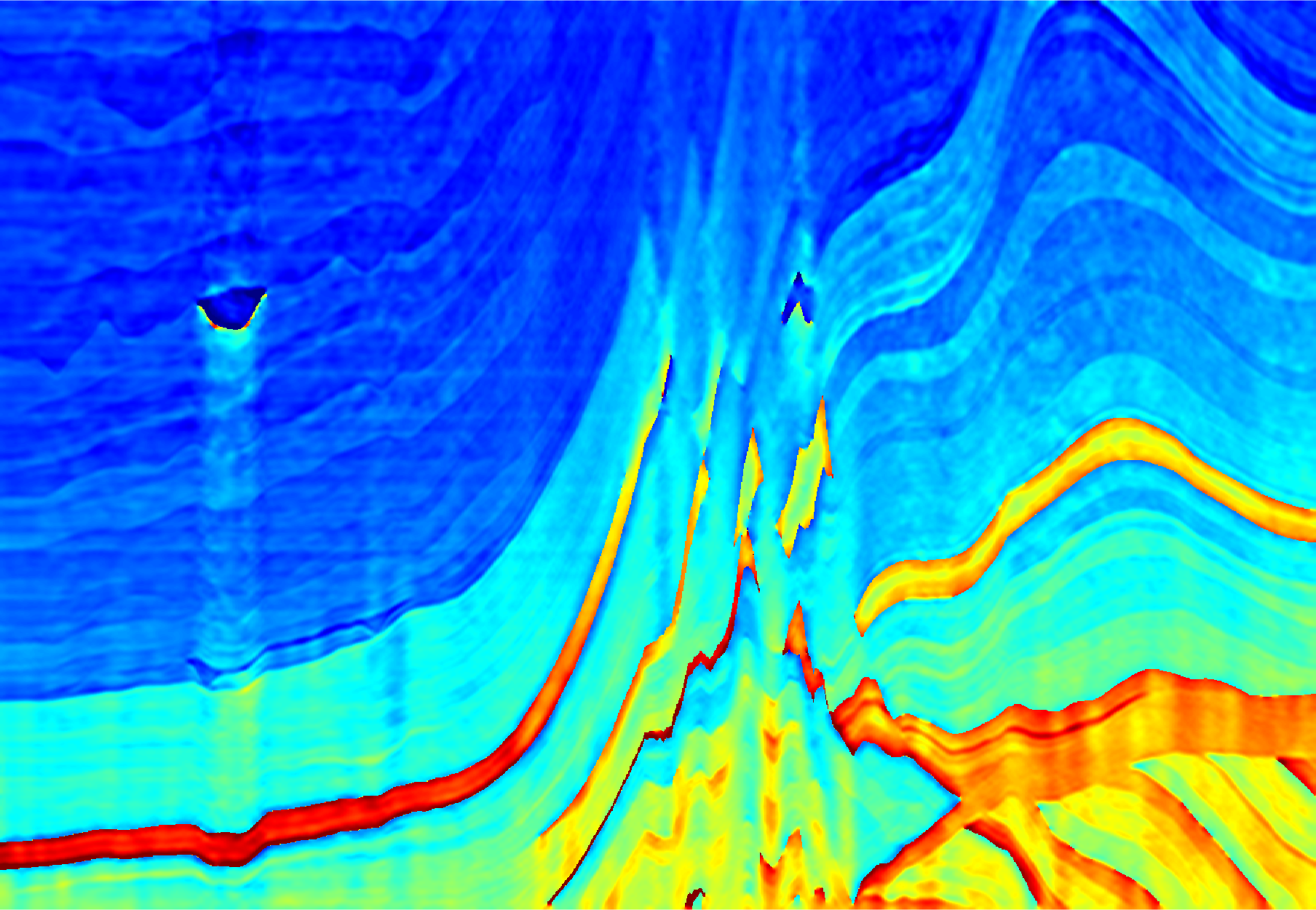}}\ \ 
    \subfloat[]{\includegraphics[height=0.22\textwidth]{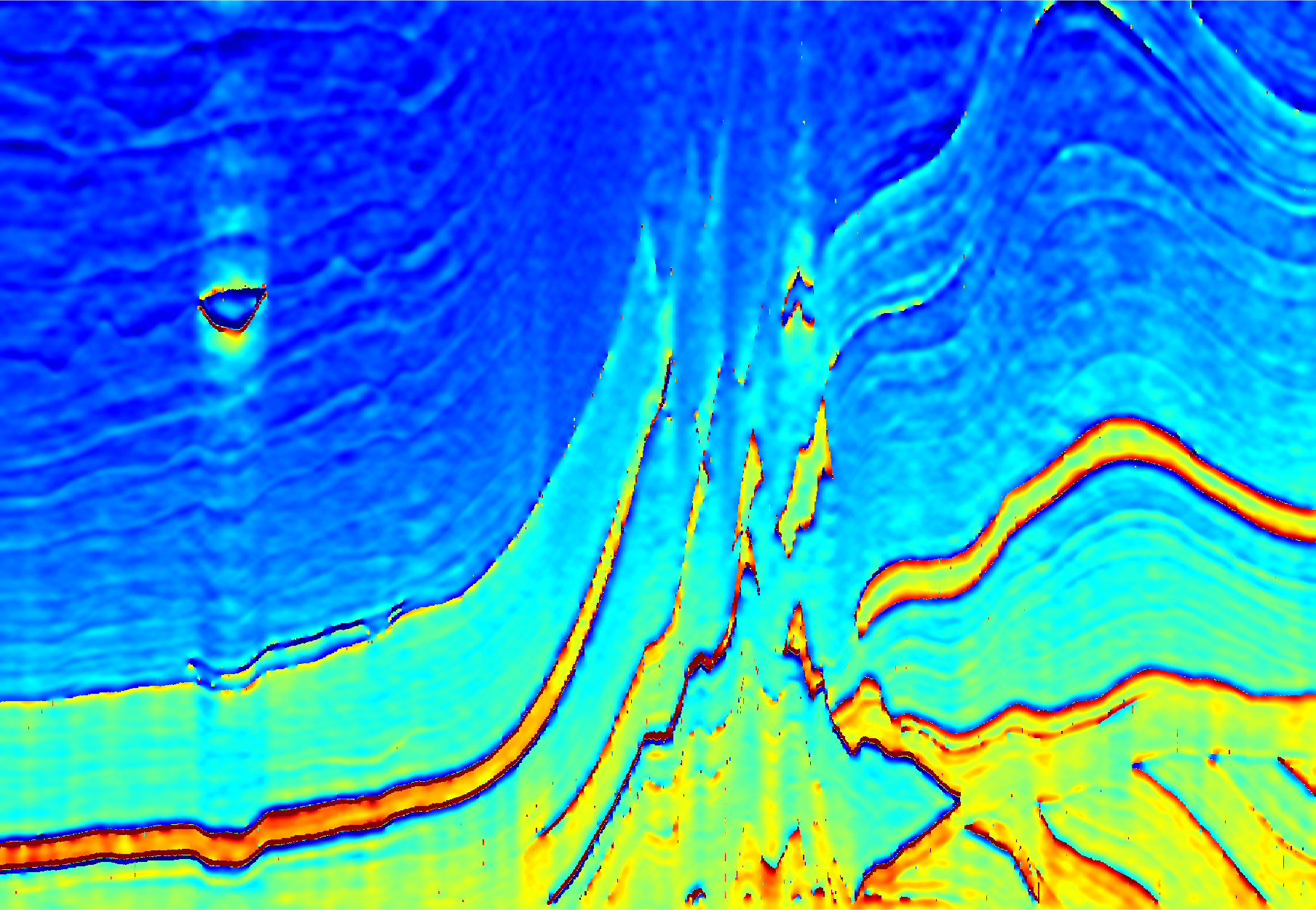}}\ \ 
    \subfloat[]{\includegraphics[height=0.22\textwidth]{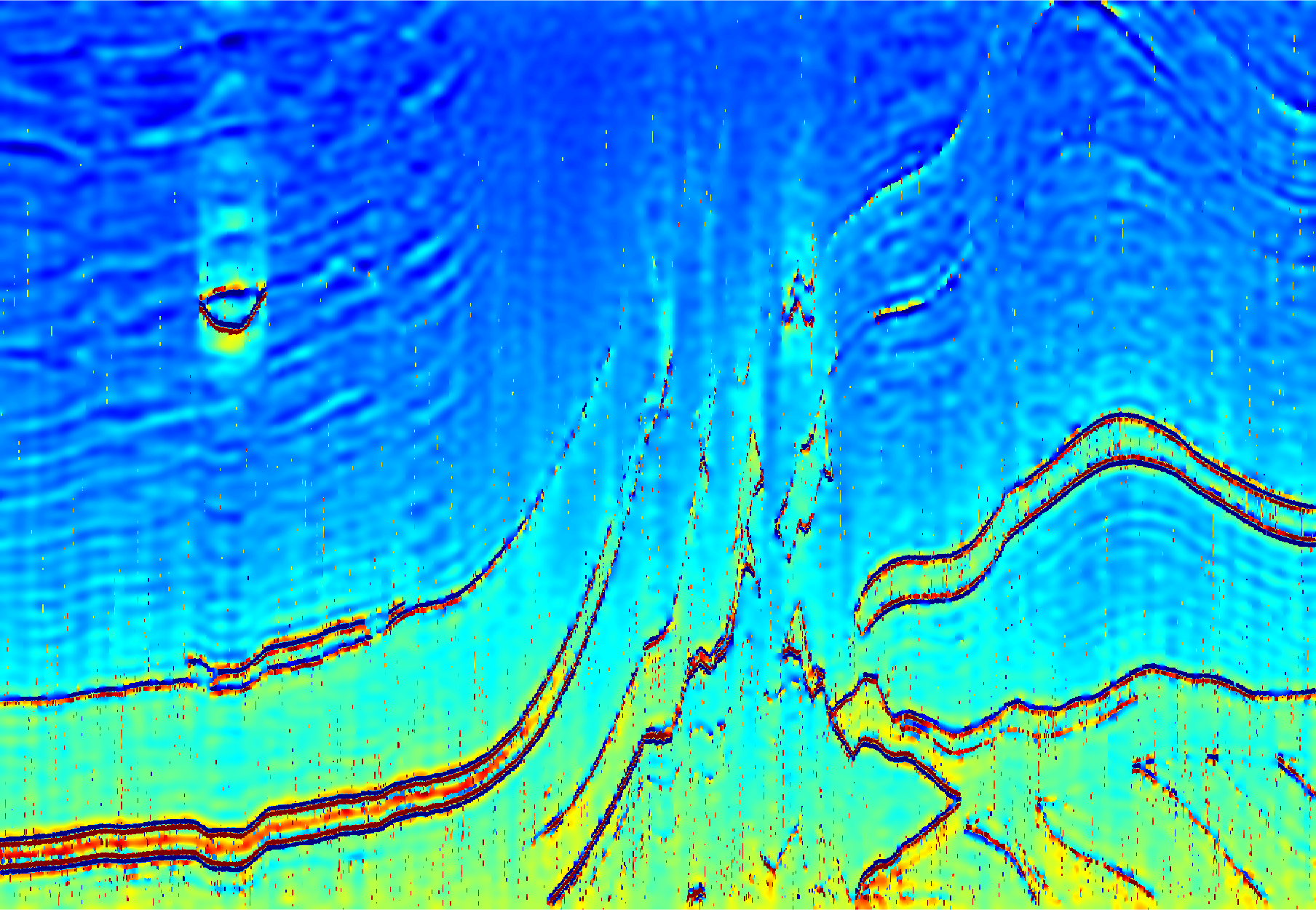}}\\
    \subfloat[]{\includegraphics[height=0.22\textwidth]{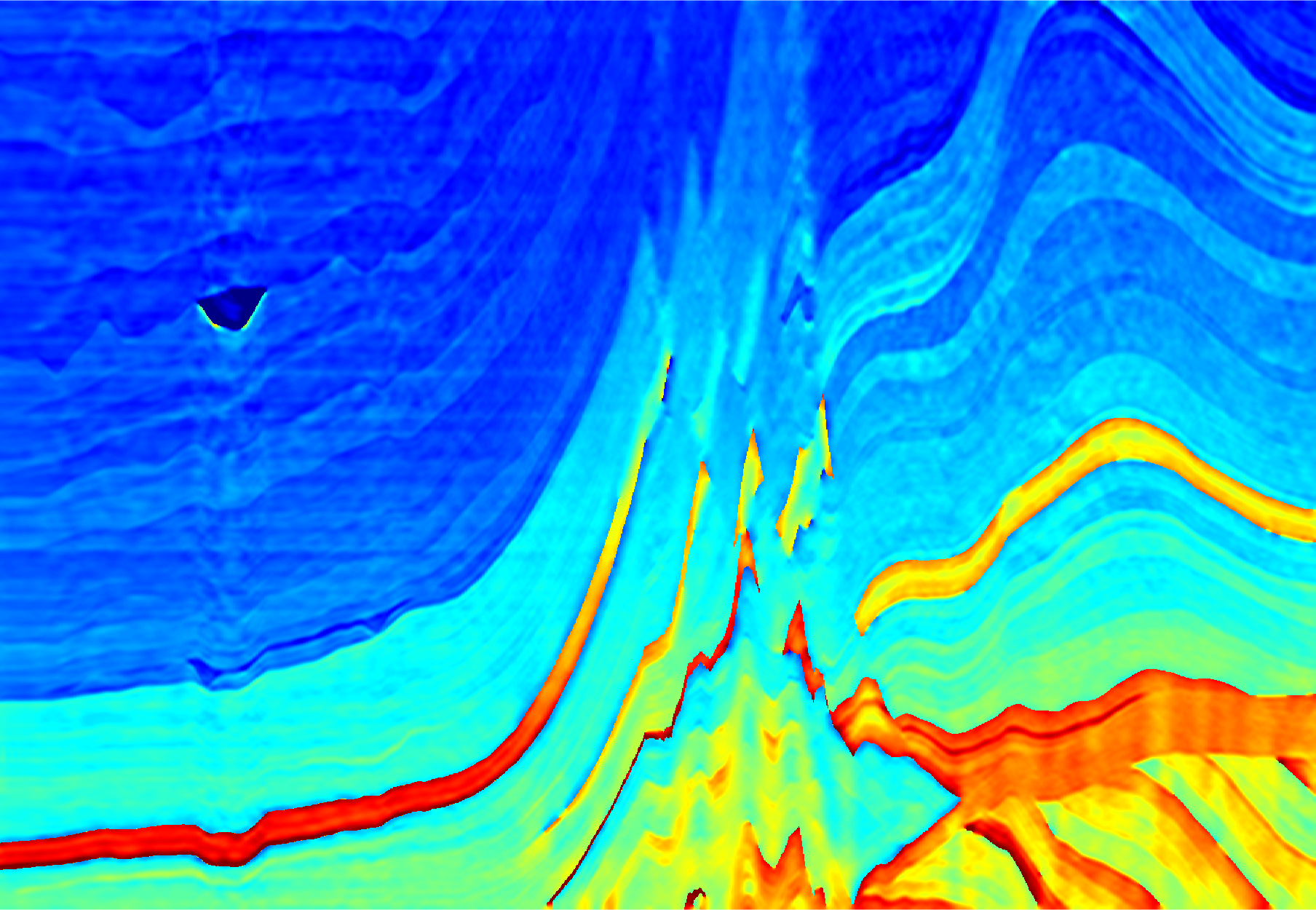}}\ \ 
    \subfloat[]{\includegraphics[height=0.22\textwidth]{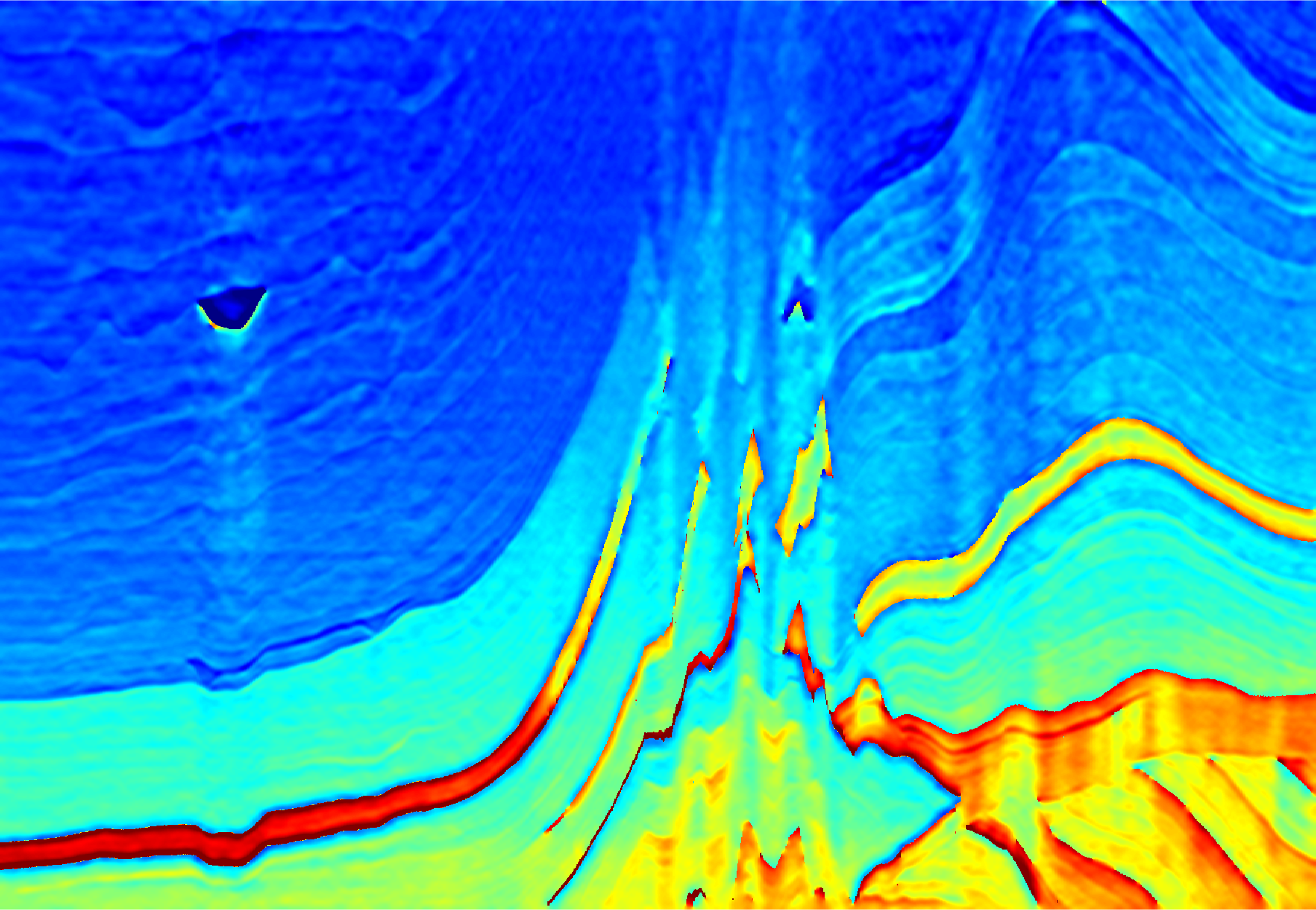}}\ \ 
    \subfloat[]{\includegraphics[height=0.22\textwidth]{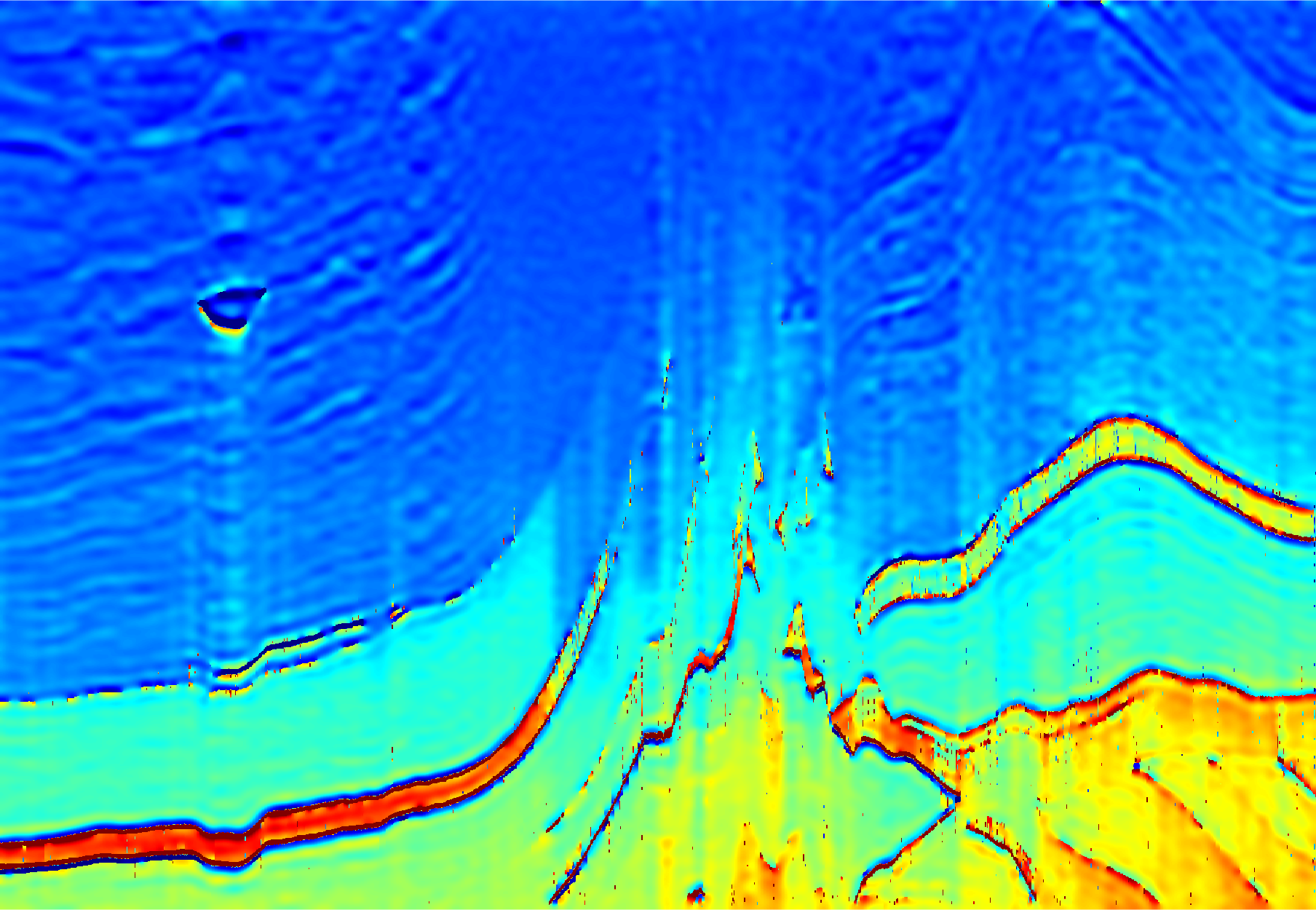}}\\
    \subfloat[]{\includegraphics[height=0.22\textwidth]{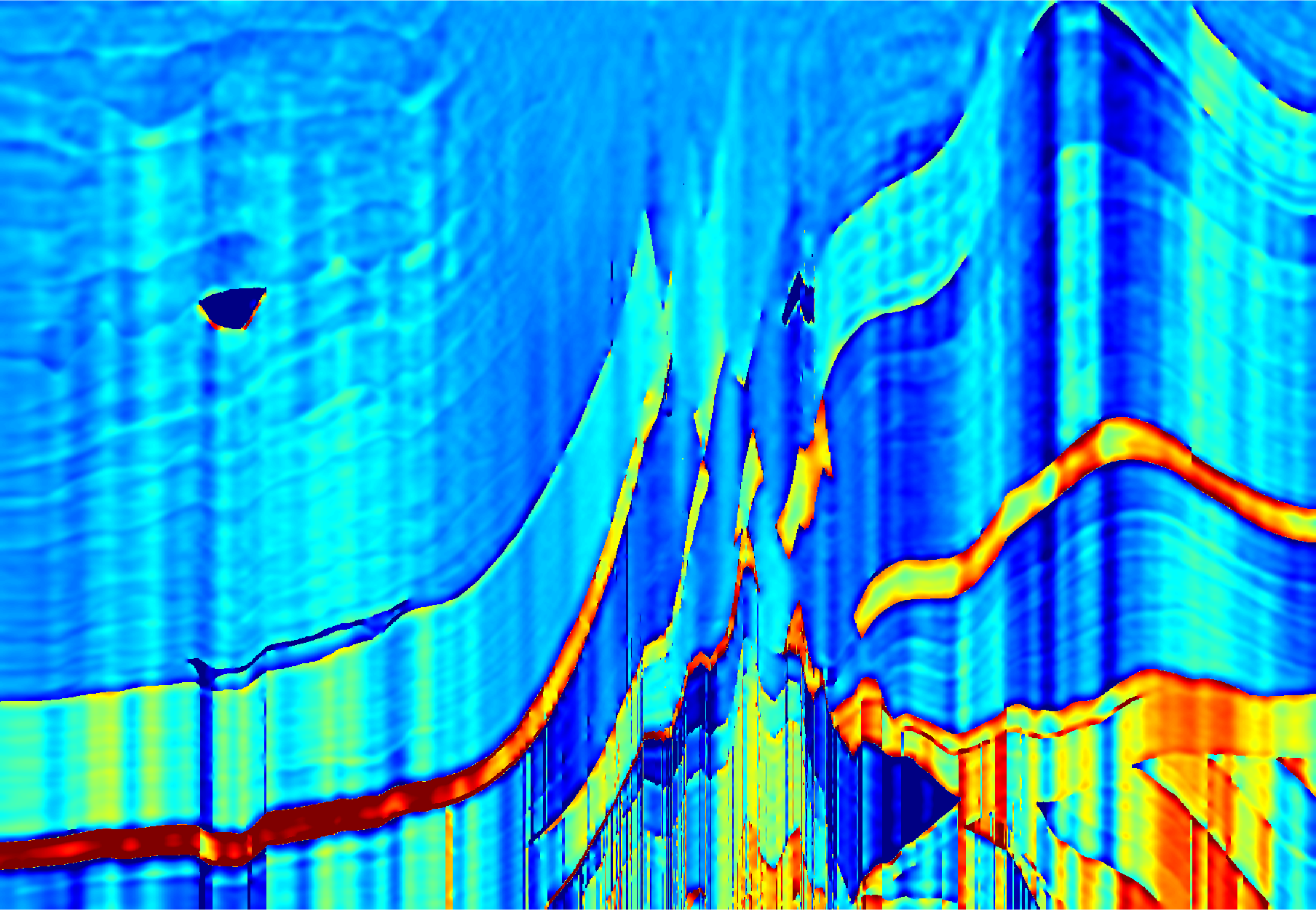}}\ \ 
    \subfloat[]{\includegraphics[height=0.22\textwidth]{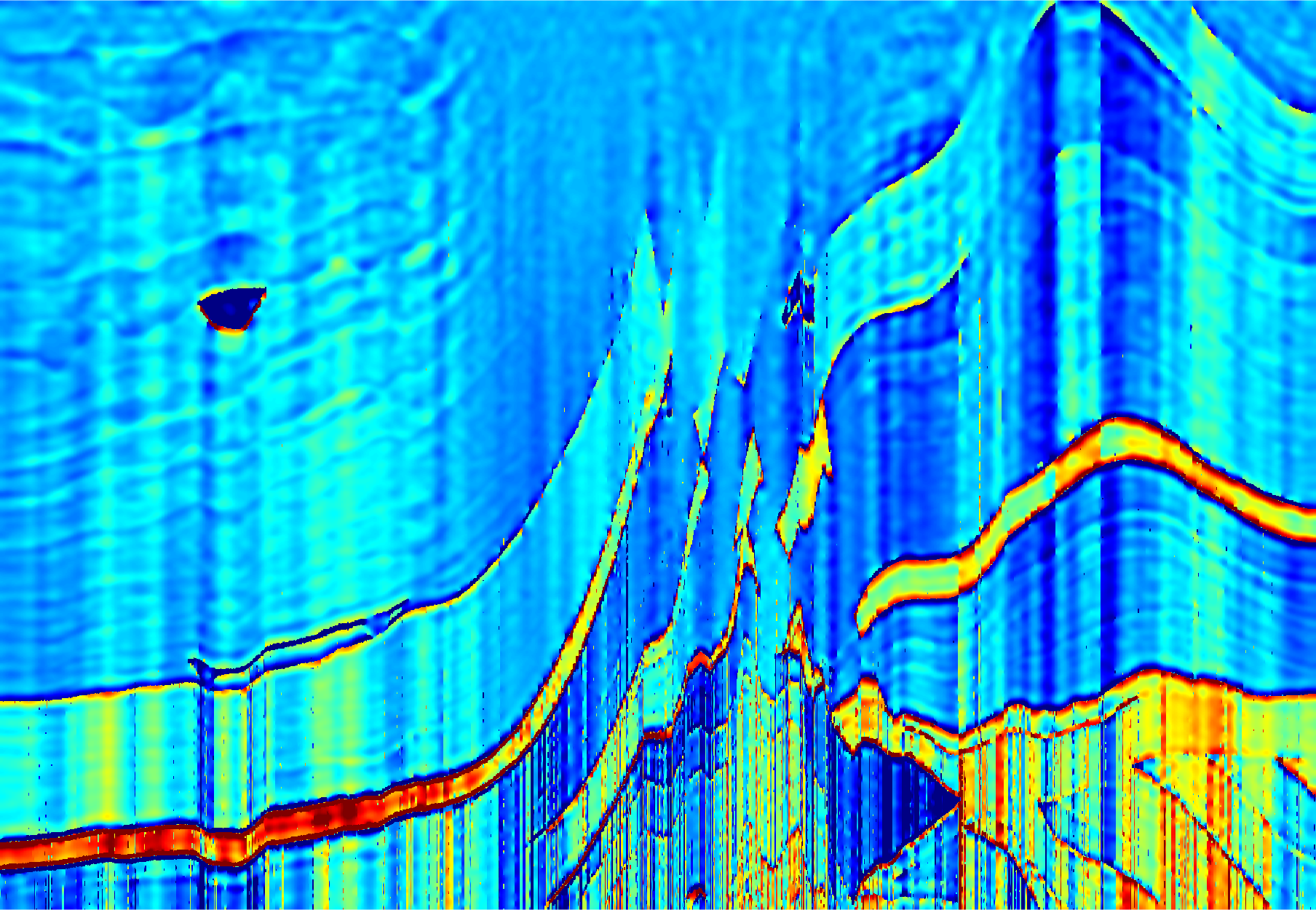}}\ \ 
    \subfloat[]{\includegraphics[height=0.22\textwidth]{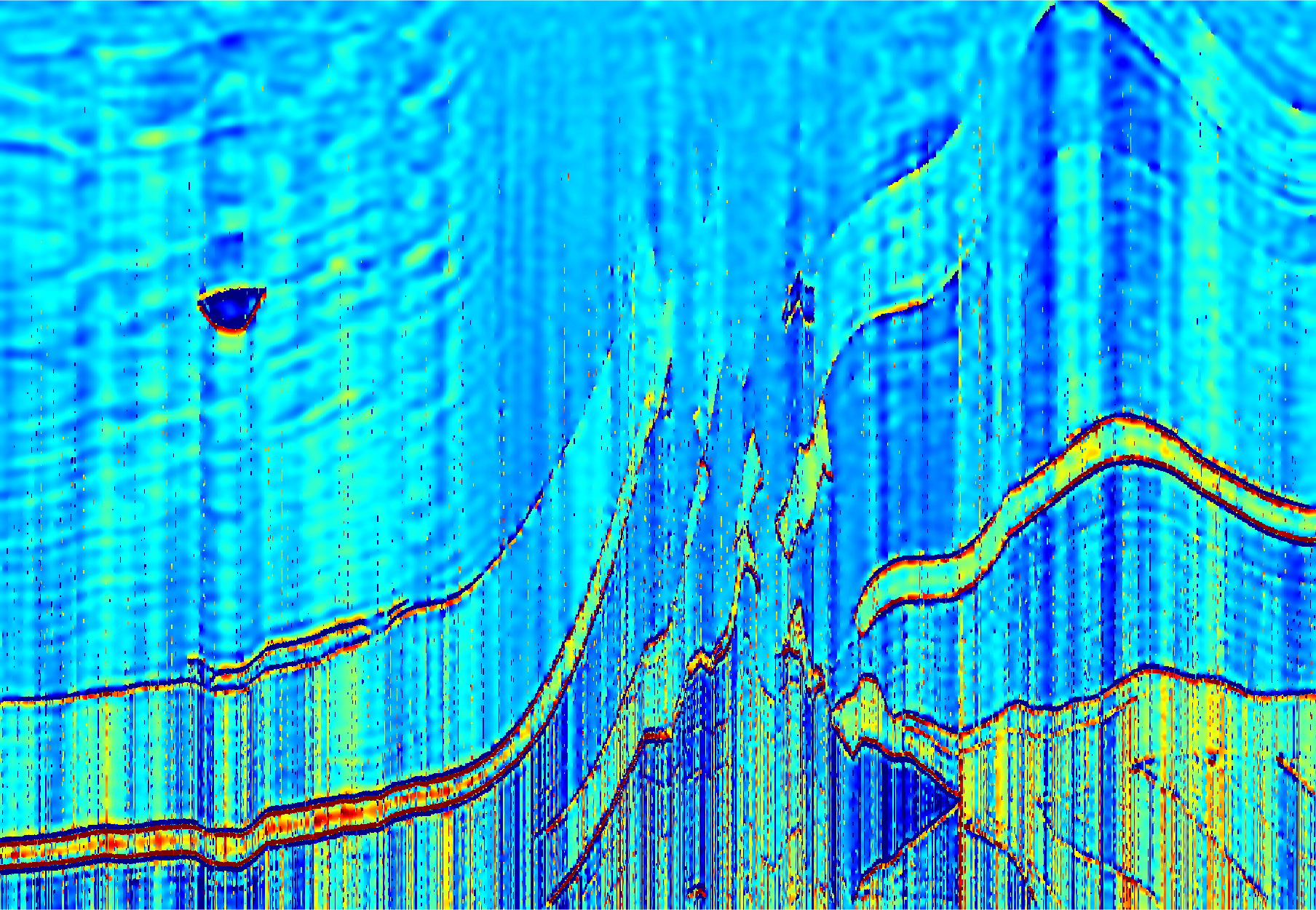}}\\
    \subfloat[]{\includegraphics[height=0.22\textwidth]{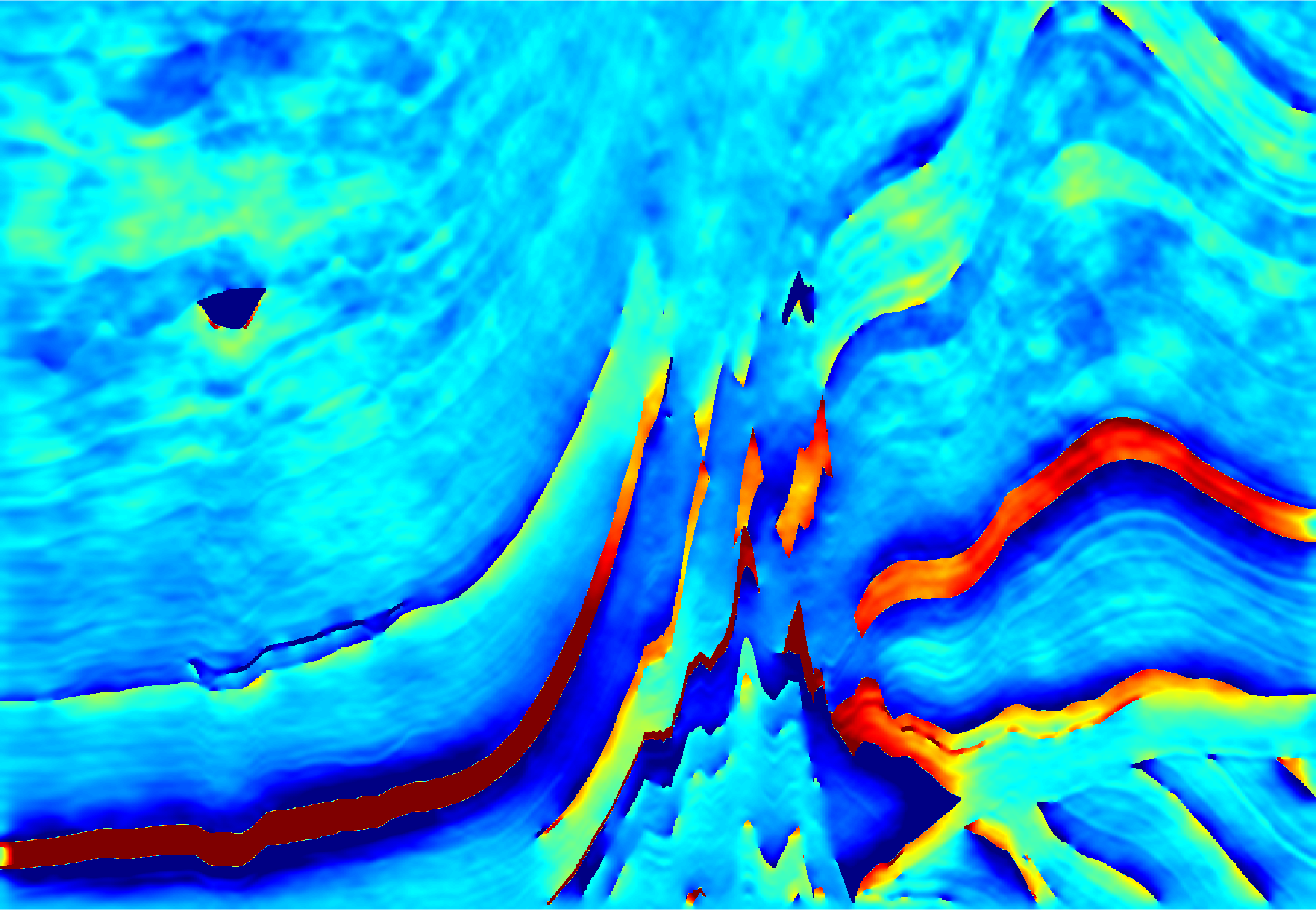}}\ \ 
    \subfloat[]{\includegraphics[height=0.22\textwidth]{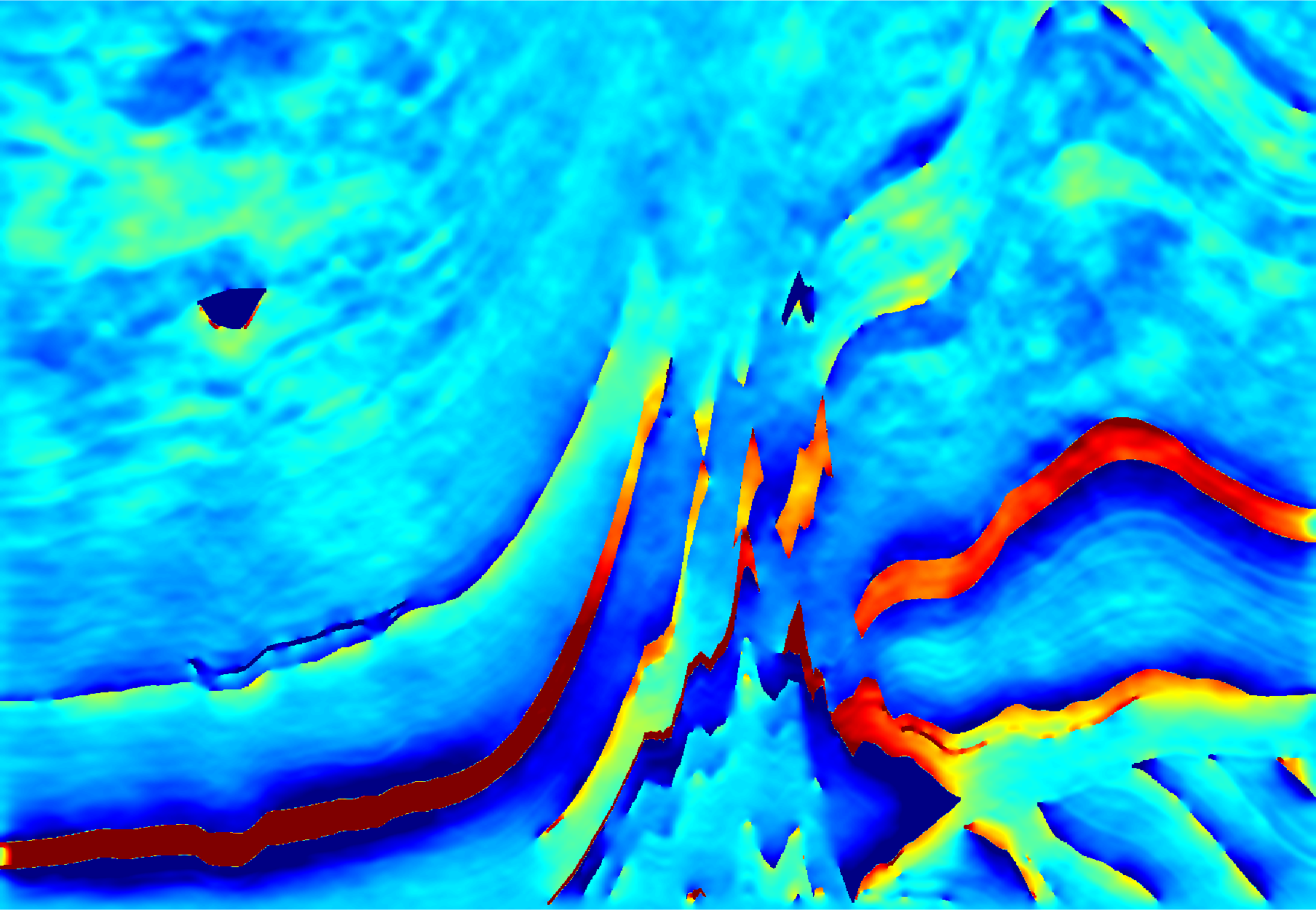}}\ \ 
    \subfloat[]{\includegraphics[height=0.22\textwidth]{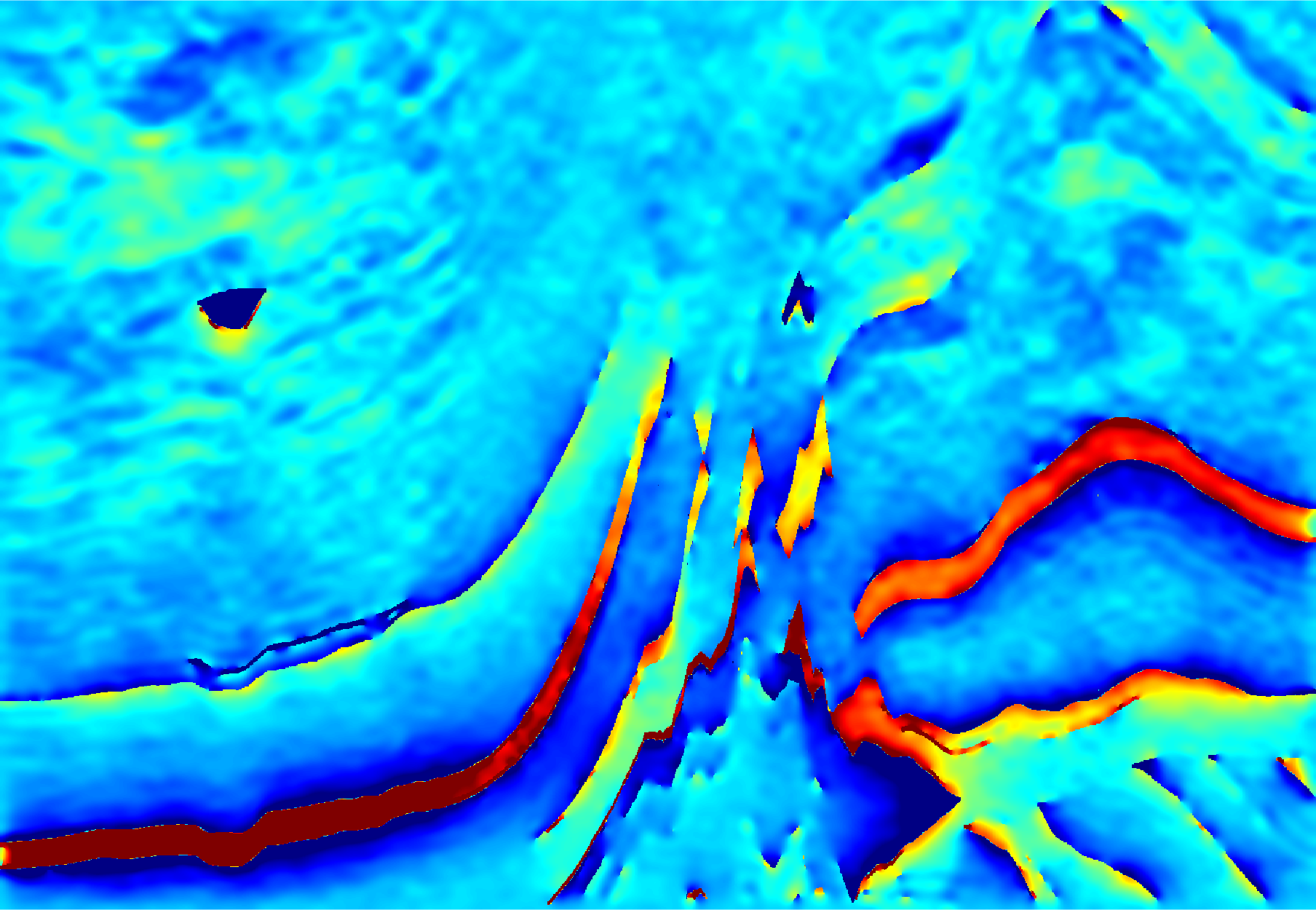}}\\
    \marmousiCBH{0.33\textwidth}
	\caption{Reconstructed impedance profiles for Marmousi2 from different initialization methods and different noise levels after 10 iterations of \ref{iterated_graphLaNet}: AA (a,b,c), Liu (d,e,f), SSI (g,h,i), SB (j,k,l); without noise (left column: a,d,g,j), medium noise (middle column: b,e,h,k), high noise (right column: c,f,i,l).}
    \label{fig:Marmousi_10iters}
\end{figure}

\begin{figure}[h!tb]
	\centering
    \subfloat[]{\includegraphics[height=0.22\textwidth]{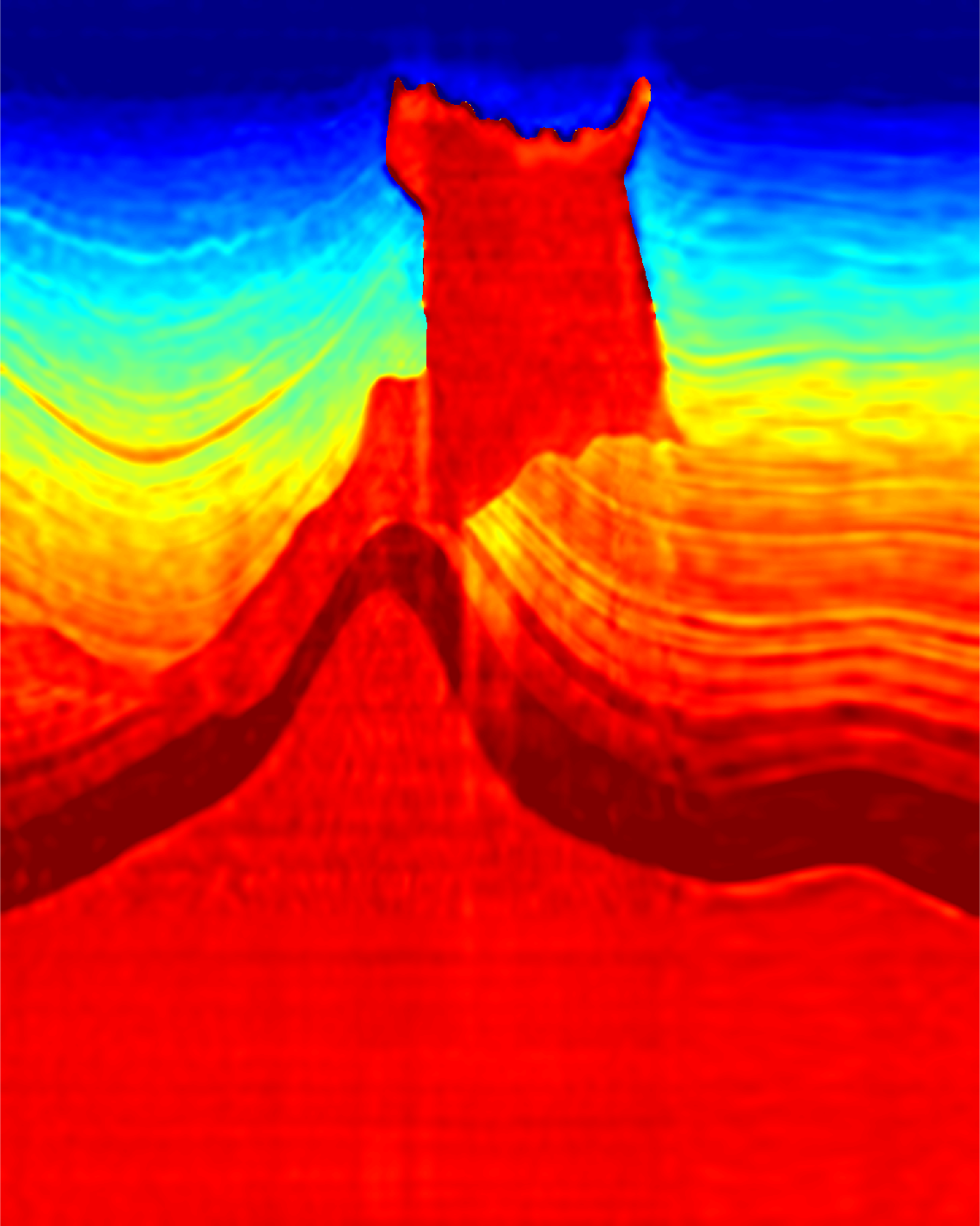}}\ \ 
    \subfloat[]{\includegraphics[height=0.22\textwidth]{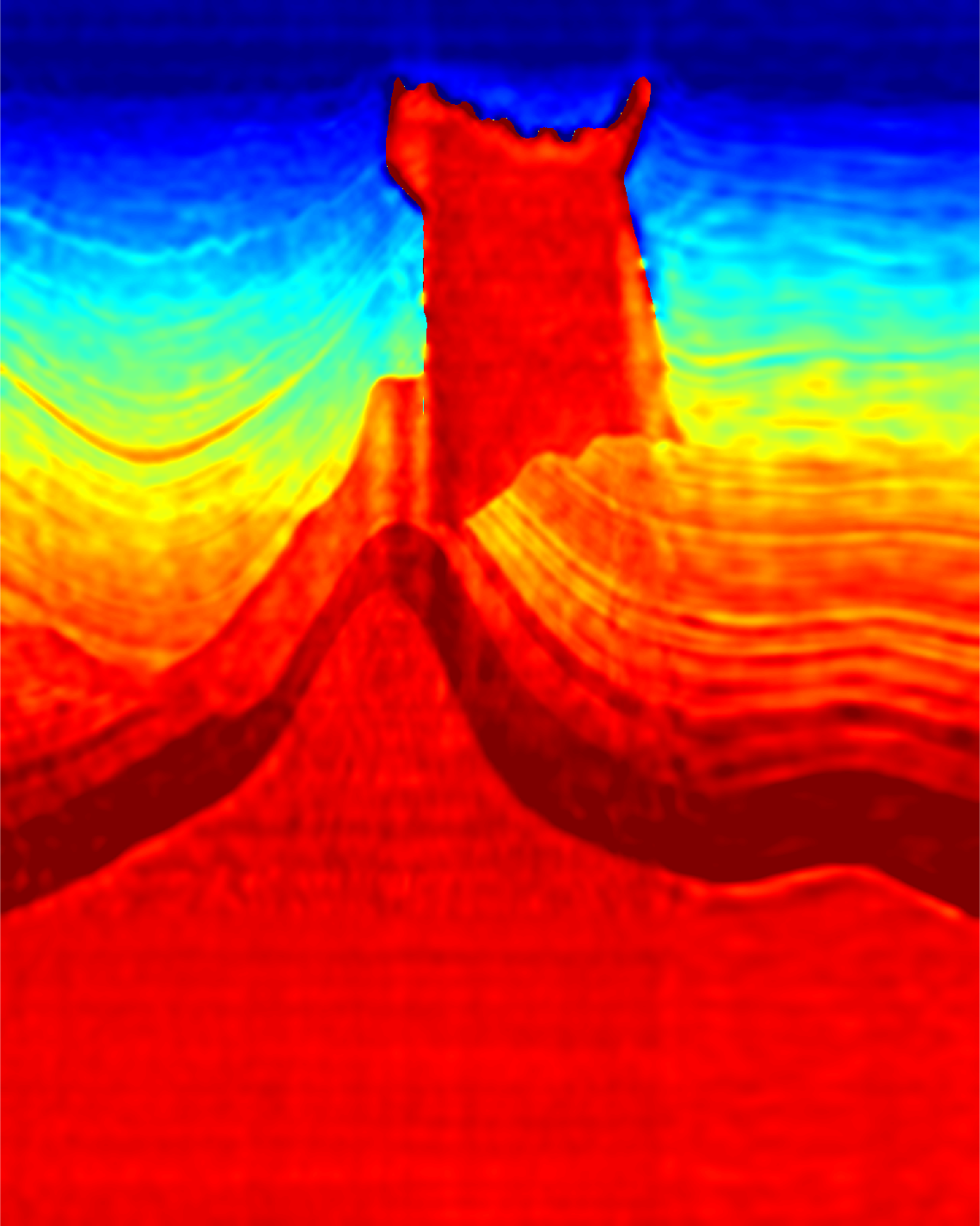}}\ \ 
    \subfloat[]{\includegraphics[height=0.22\textwidth]{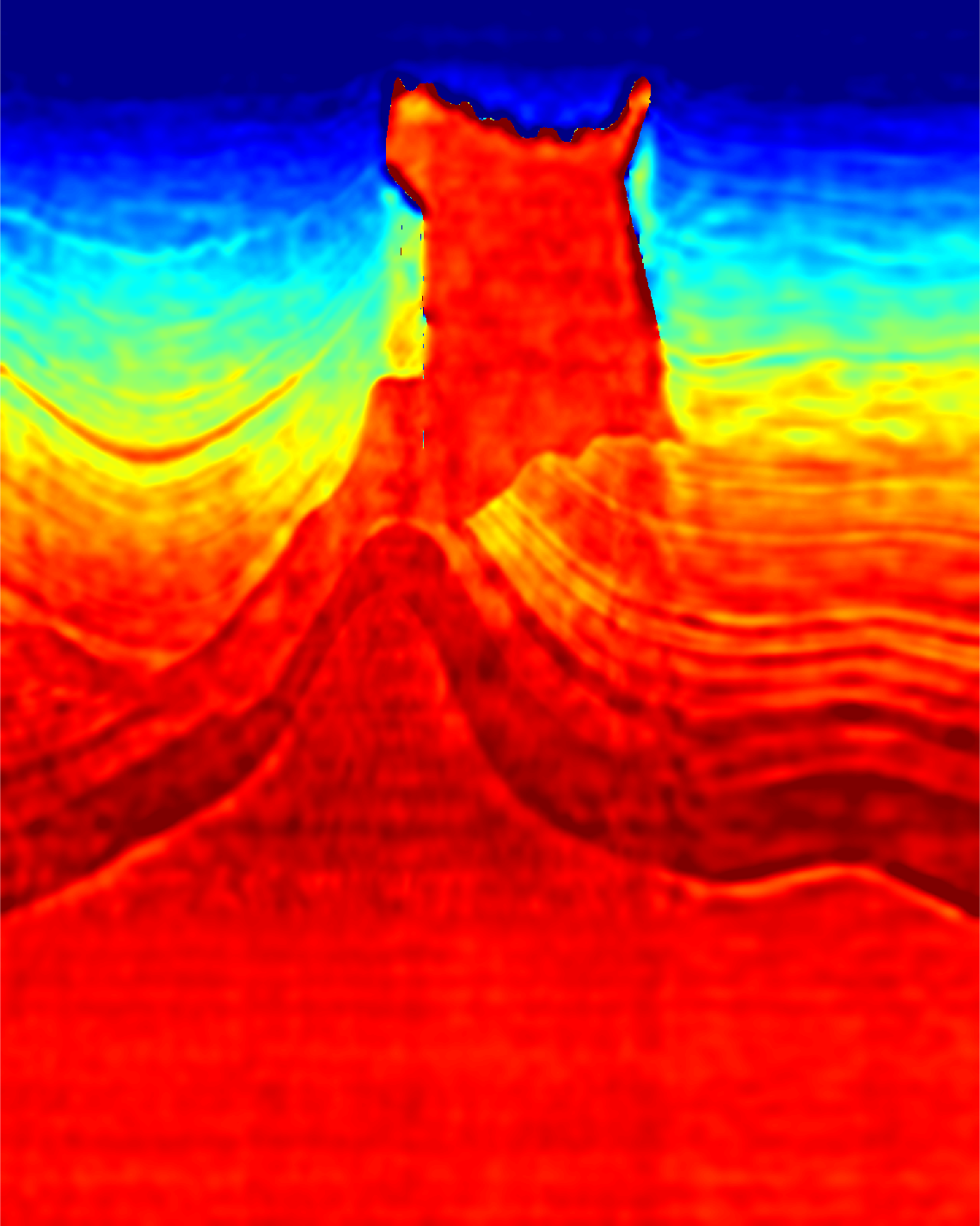}}\\
    \subfloat[]{\includegraphics[height=0.22\textwidth]{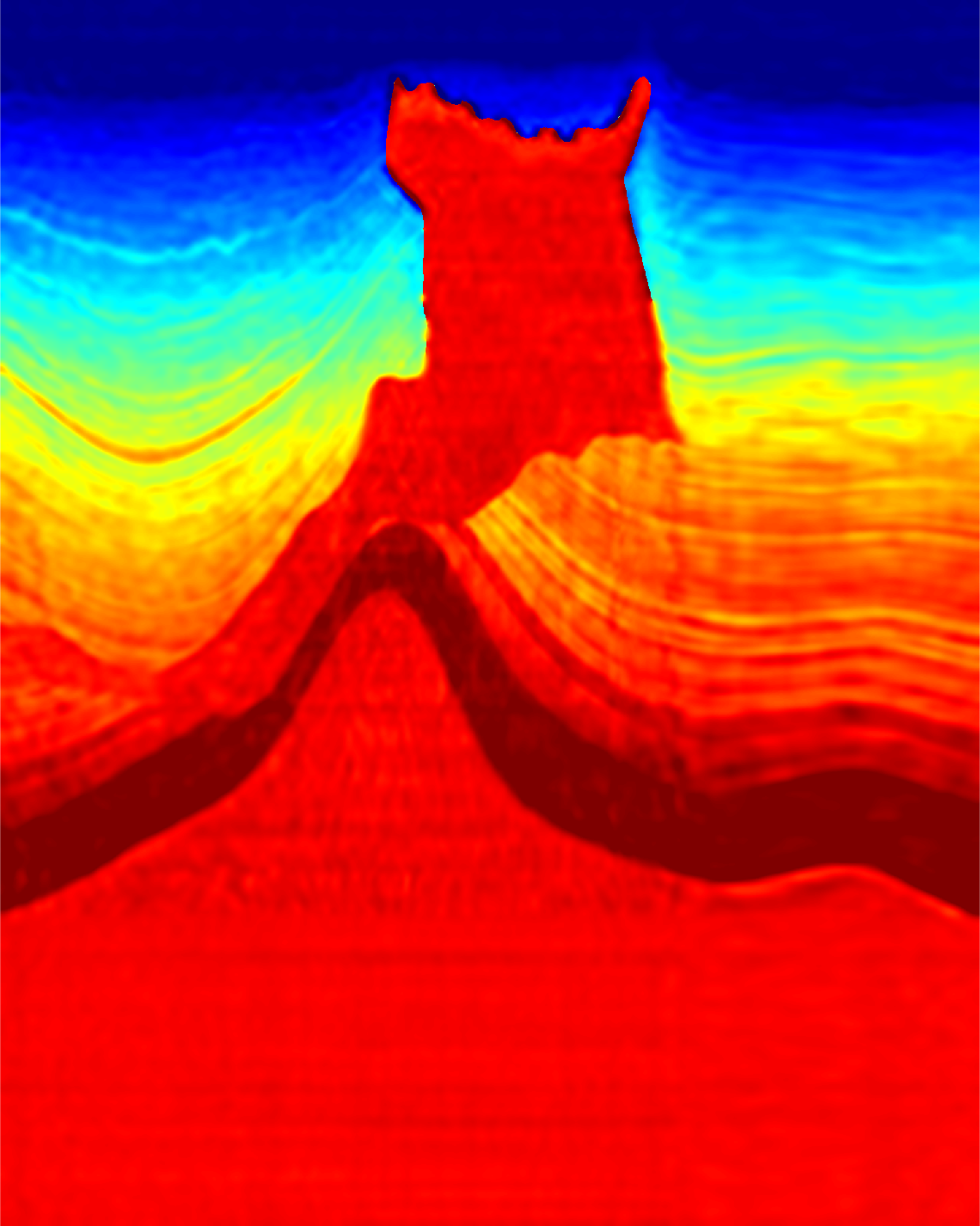}}\ \ 
    \subfloat[]{\includegraphics[height=0.22\textwidth]{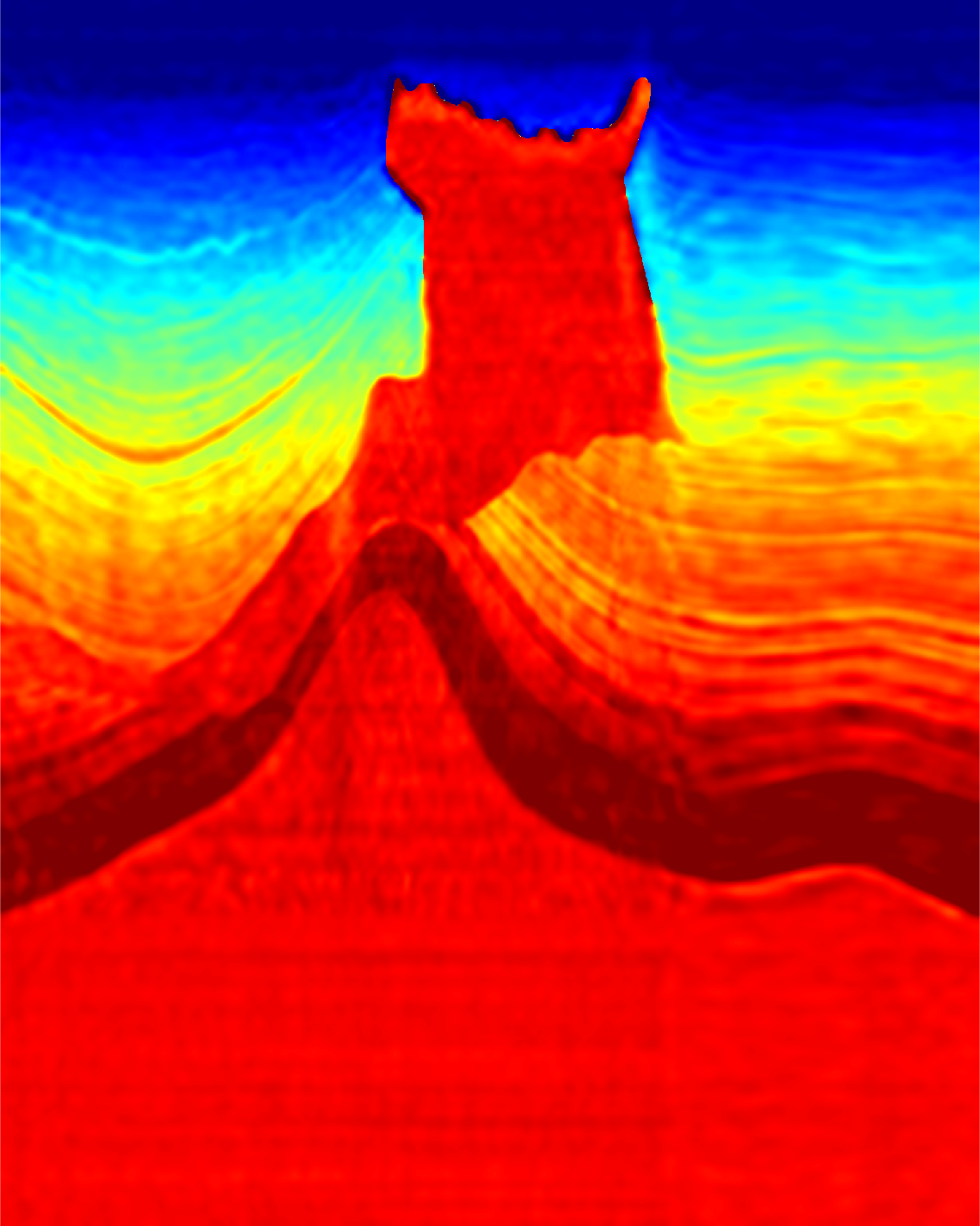}}\ \ 
    \subfloat[]{\includegraphics[height=0.22\textwidth]{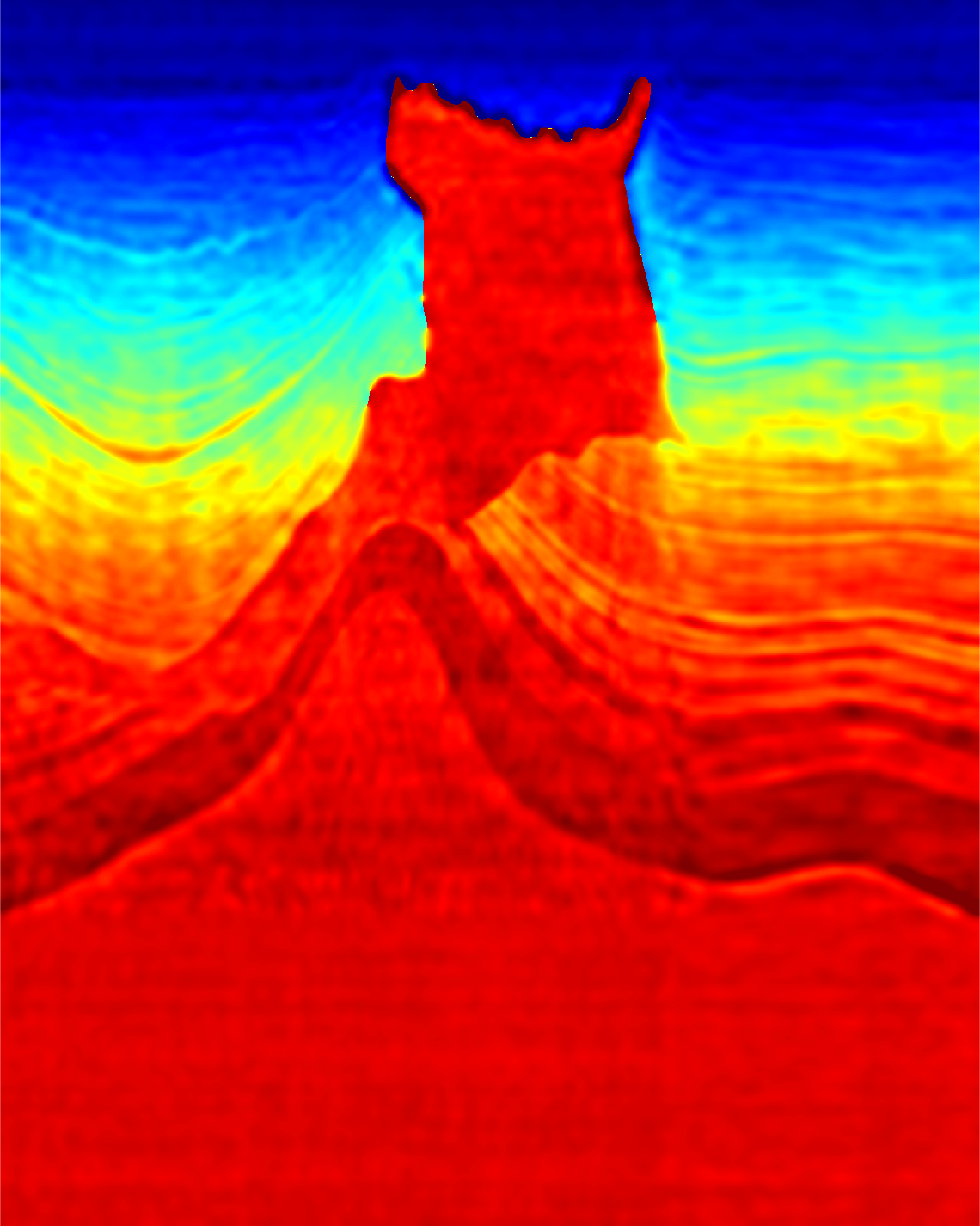}}\\
    \subfloat[]{\includegraphics[height=0.22\textwidth]{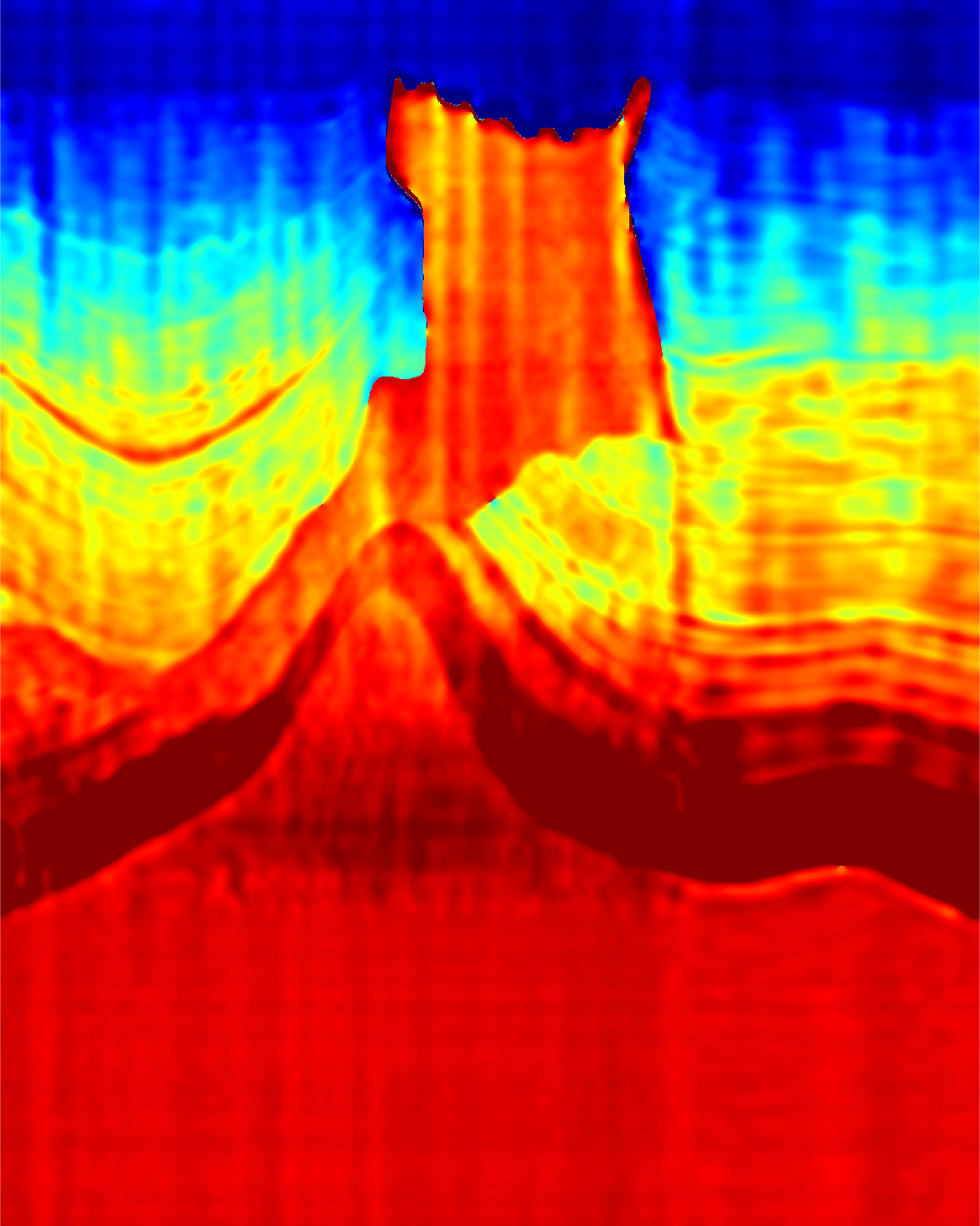}}\ \ 
    \subfloat[]{\includegraphics[height=0.22\textwidth]{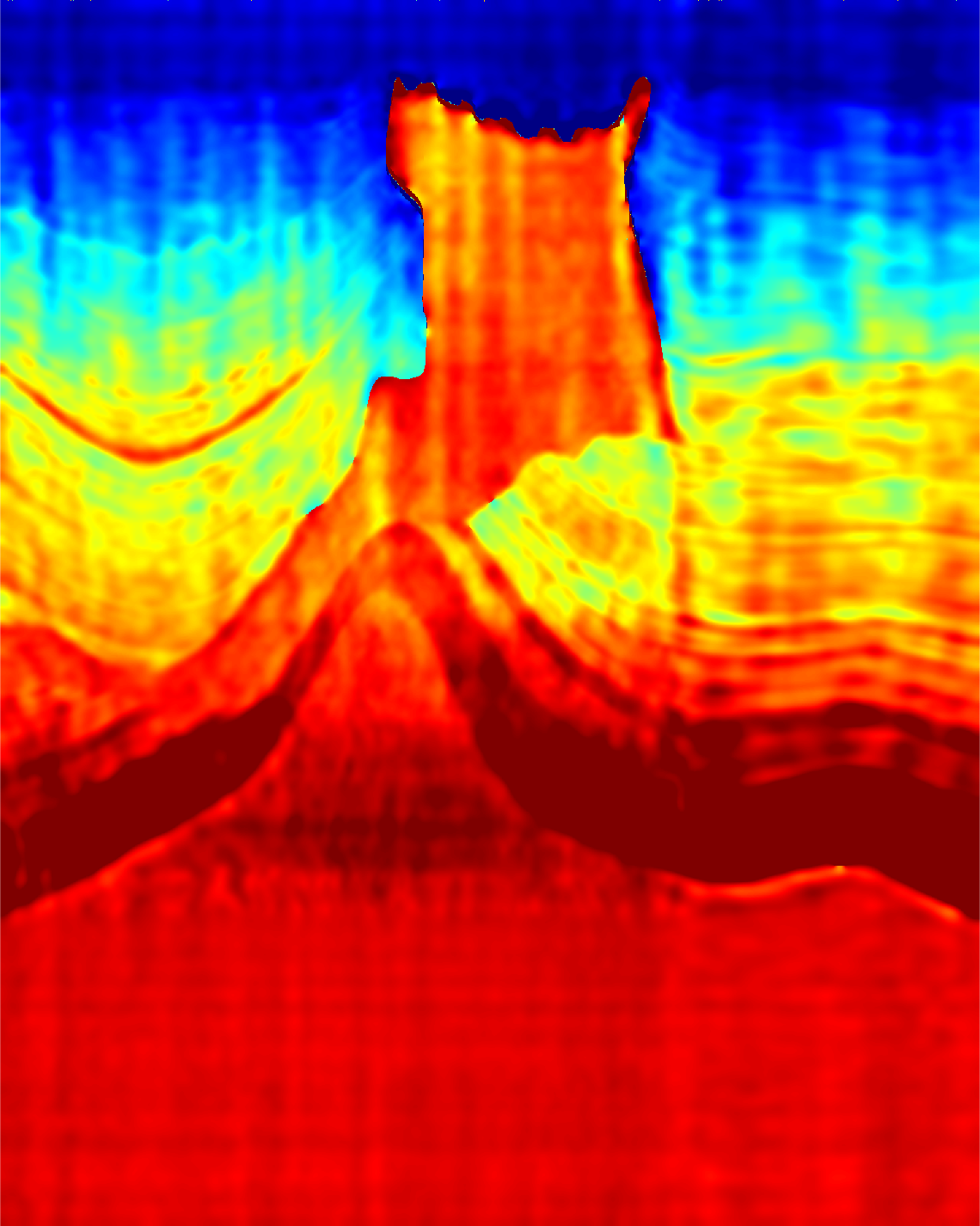}}\ \ 
    \subfloat[]{\includegraphics[height=0.22\textwidth]{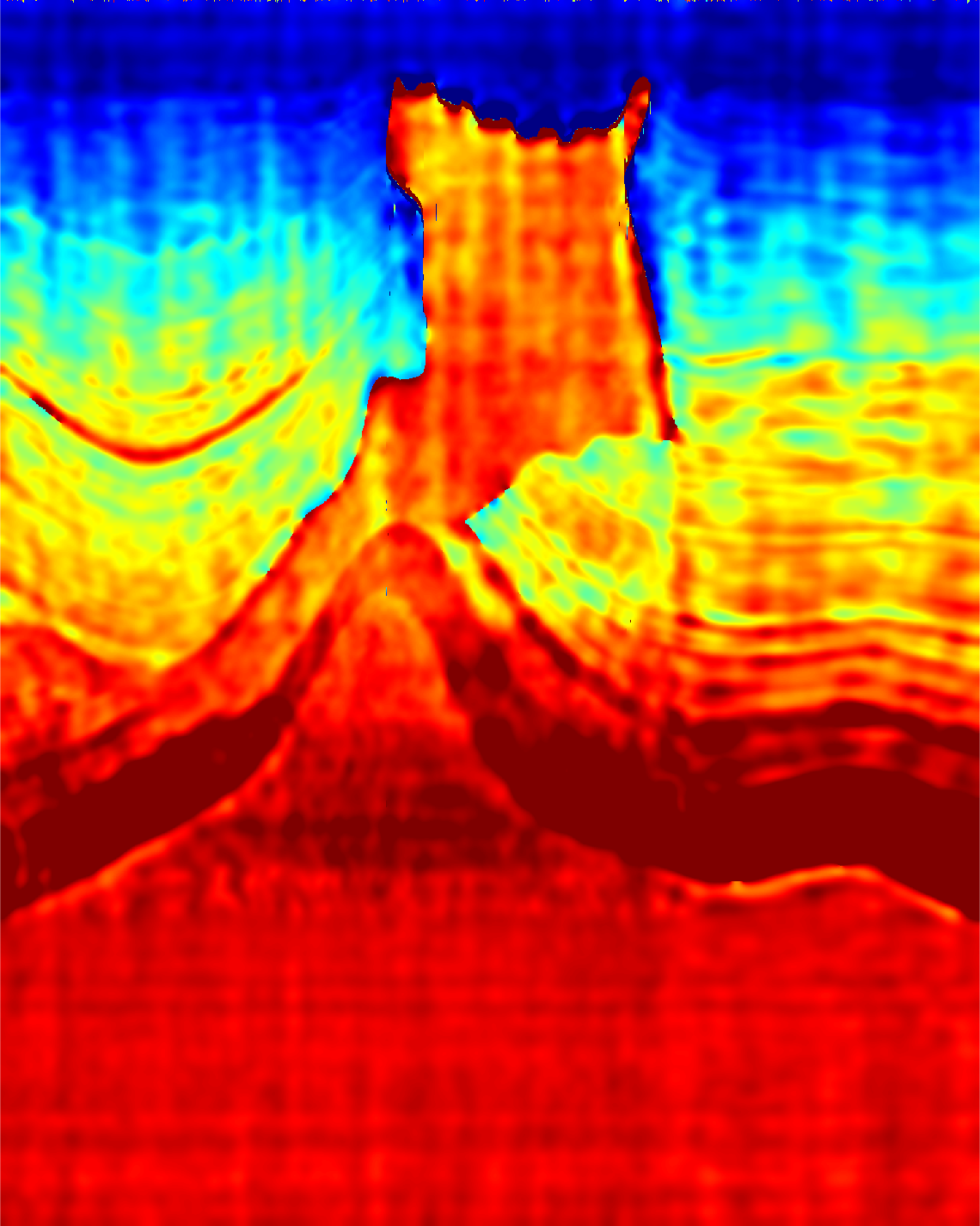}}\\
    \subfloat[]{\includegraphics[height=0.22\textwidth]{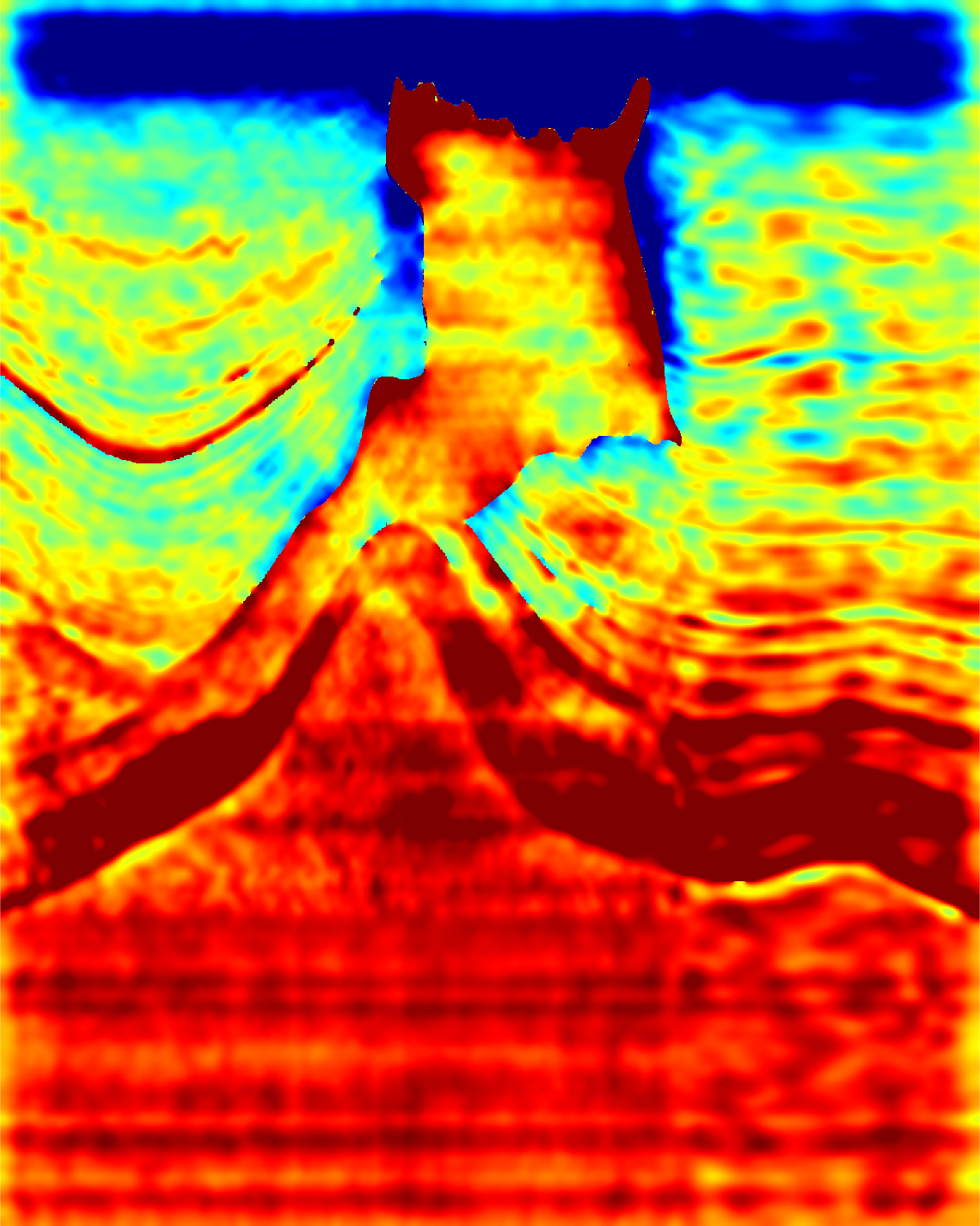}}\ \ 
    \subfloat[]{\includegraphics[height=0.22\textwidth]{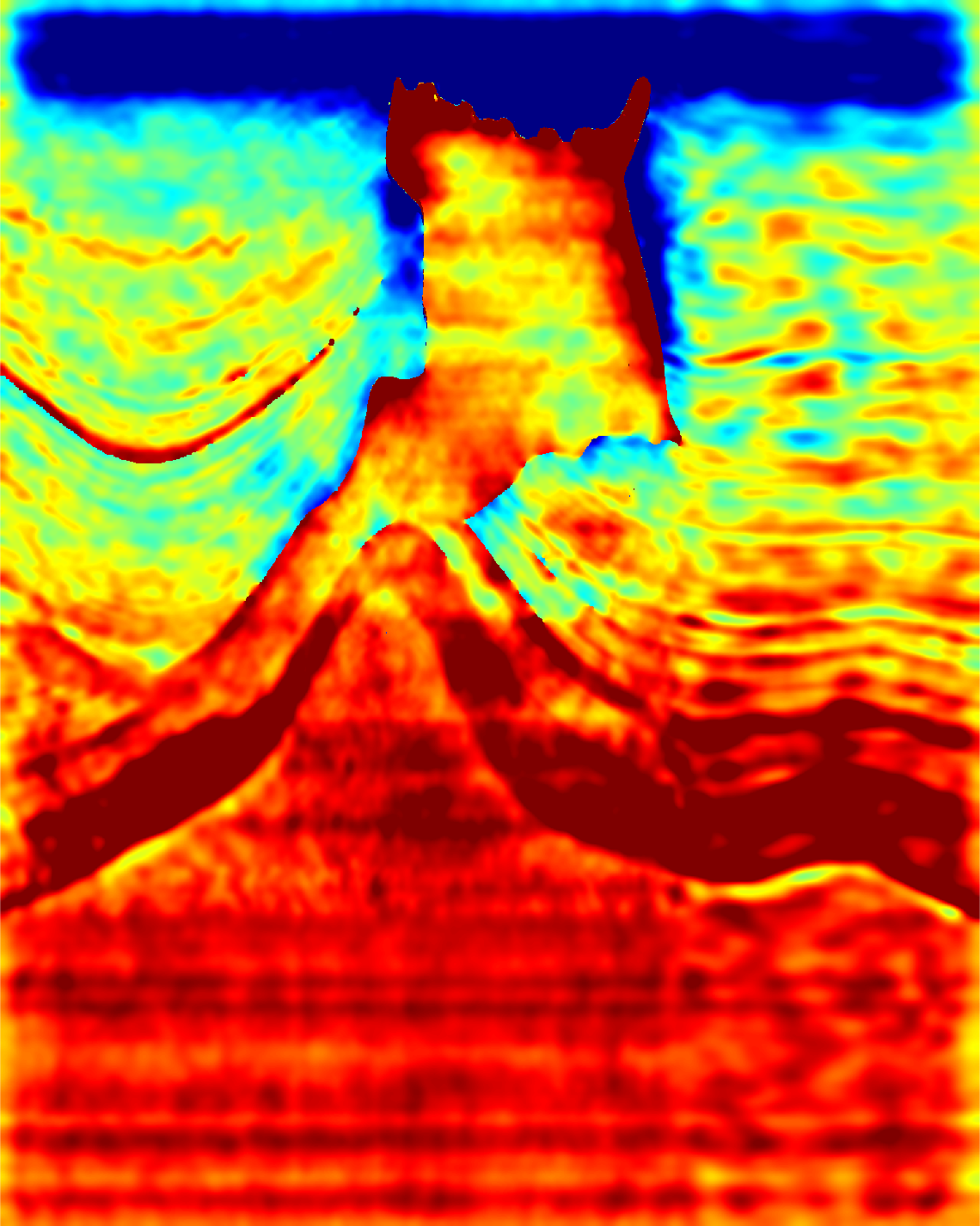}}\ \ 
    \subfloat[]{\includegraphics[height=0.22\textwidth]{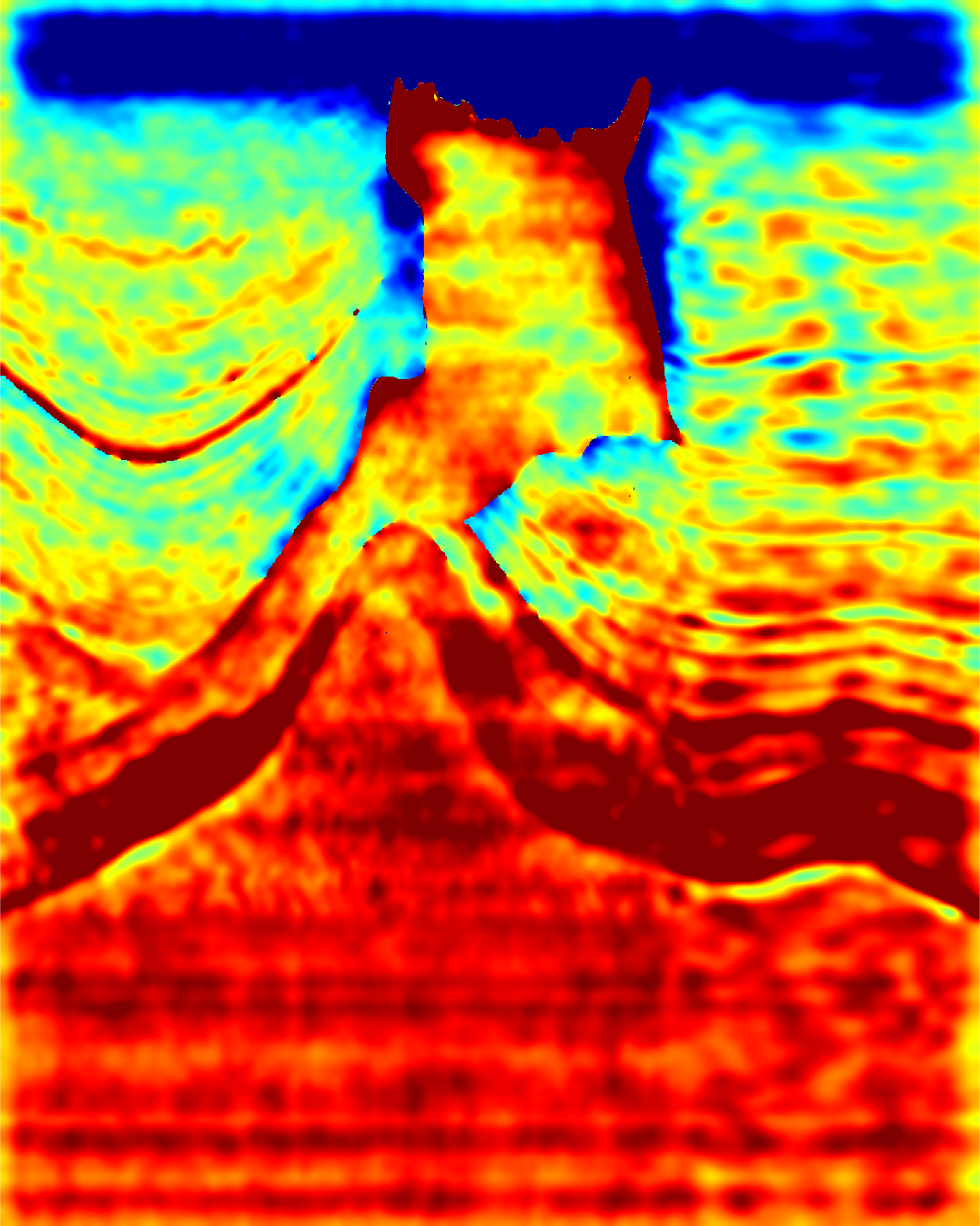}}\\
    \seamCBH{0.33\textwidth}
	\caption{Reconstructed impedance profiles for SEAM from different initialization methods and different noise levels after 10 iterations of \ref{iterated_graphLaNet}: AA (a,b,c), Liu (d,e,f), SSI (g,h,i), SB (j,k,l); without noise (left column: a,d,g,j), medium noise (middle column: b,e,h,k), high noise (right column: c,f,i,l).}
    \label{fig:SEAM_10iters}
\end{figure}

Based on the above considerations, we used the parameter setup $R=2$, $\sigma=0.25$ and 10 iterations for all other experiments. Note that the graph Laplacian matrix has a much higher sparsity for $R=2$ compared to higher settings. The obtained D-MSE and SSIM values for all noise levels and initialization methods are shown in Table \ref{tab:errorsMarmousi} for the Marmousi model and Table \ref{tab:errorsSEAM} for the SEAM model. We show the values obtained by the initial reconstruction as well as after 10 iterations of the proposed method. In our experiments the SSIM value was a good indicator for the network based initialization methods while the D-MSE is much more reliable for the classical initialization methods. We see that we are able to improve the reconstruction quality in almost all cases. Only in the noiseless case or for very low noise the Liu initialization (which obtained the best results of all initialization methods) performed better without \texttt{it-graphLa}$\Psi$. In these cases \texttt{it-graphLa}$\Psi$ tends to overregularize and reduce the reconstruction quality. Figures \ref{fig:Marmousi_10iters} and \ref{fig:SEAM_10iters} show the final reconstructions in the noiseless cases as well as for medium and high noise (compare Figure \ref{fig:Marmousi_init} and \ref{fig:SEAM_init}). Again, the data was normalized and colorcoded in the same manner as all previous images for a better comparison.

\begin{table}
    \centering    
\begin{tabular}{|c|c|c|c|c|c|c|}
    \hline
    \multicolumn{2}{|r|}{PSNR} & noiseless & $39$ & $33$ & $30$ & $27$ \\
    \hhline{=======}
    \multirow{4}{*}{\rotatebox[origin=c]{90}{initialization}}
    & AA &
    $\begin{matrix}0.004397\\0.78311\end{matrix}$ & 
    $\begin{matrix}0.0050761\\0.57689\end{matrix}$ & 
    $\begin{matrix}0.0091736\\0.37902\end{matrix}$ & 
    $\begin{matrix}0.010019\\0.2279\end{matrix}$ & 
    $\begin{matrix}0.016137\\0.17054\end{matrix}$ \\
    \cline{2-7}
    & Liu &
    $\begin{matrix}0.0034844\\0.8898\end{matrix}$ & 
    $\begin{matrix}0.0039259\\0.80405\end{matrix}$ & 
    $\begin{matrix}0.0098691\\0.61324\end{matrix}$ & 
    $\begin{matrix}0.017726\\0.46227\end{matrix}$ & 
    $\begin{matrix}0.014271\\0.4658\end{matrix}$ \\
    \cline{2-7}
    & SSI &
    $\begin{matrix}0.0060146\\0.22754\end{matrix}$ & 
    $\begin{matrix}0.0061296\\0.15793\end{matrix}$ & 
    $\begin{matrix}0.0063865\\0.093704\end{matrix}$ & 
    $\begin{matrix}0.0064695\\0.056462\end{matrix}$ & 
    $\begin{matrix}0.0064667\\0.034648\end{matrix}$ \\
    \cline{2-7}
    & SB &
    $\begin{matrix}0.017705\\-0.011821\end{matrix}$ & 
    $\begin{matrix}0.018717\\-0.012986\end{matrix}$ & 
    $\begin{matrix}0.020986\\-0.014372\end{matrix}$ & 
    $\begin{matrix}0.027666\\-0.014848\end{matrix}$ & 
    $\begin{matrix}0.035281\\-0.01365\end{matrix}$ \\
    \hhline{=======}
    \multirow{4}{*}{\rotatebox[origin=c]{90}{after 10 iterations}}
    & AA &
    $\begin{matrix}0.0062188\\0.83374\end{matrix}$ &
    $\begin{matrix}0.005795\\0.78584\end{matrix}$ & 
    $\begin{matrix}0.0050861\\0.73877\end{matrix}$ & 
    $\begin{matrix}0.0056011\\0.63481\end{matrix}$ & 
    $\begin{matrix}0.0060183\\0.56238\end{matrix}$ \\
    \cline{2-7}
    & Liu &
    $\begin{matrix}0.0053559\\0.87256\end{matrix}$ & 
    $\begin{matrix}0.005641\\0.86061\end{matrix}$ & 
    $\begin{matrix}0.0054828\\0.80637\end{matrix}$ & 
    $\begin{matrix}0.0053539\\0.7578\end{matrix}$ & 
    $\begin{matrix}0.0057657\\0.73285\end{matrix}$ \\
    \cline{2-7}
    & SSI &
    $\begin{matrix}0.0040958\\0.2628\end{matrix}$ & 
    $\begin{matrix}0.0041885\\0.23842\end{matrix}$ & 
    $\begin{matrix}0.0045508\\0.18316\end{matrix}$ & 
    $\begin{matrix}0.0050804\\0.12671\end{matrix}$ & 
    $\begin{matrix}0.0057562\\0.084302\end{matrix}$ \\
    \cline{2-7}
    & SB &
    $\begin{matrix}0.0038658\\-0.013646\end{matrix}$ & 
    $\begin{matrix}0.0039444\\-0.01458\end{matrix}$ & 
    $\begin{matrix}0.0039989\\-0.016986\end{matrix}$ & 
    $\begin{matrix}0.0042669\\-0.01899\end{matrix}$ & 
    $\begin{matrix}0.005145\\-0.019798\end{matrix}$ \\
    \hline
\end{tabular}
    \caption{D-MSE (top value) and SSIM (bottom value) for reconstructions of the Marmousi model with different levels of noise. The top half of the table shows the error values for different initialization methods and the lower half shows the values after 10 iterations of the proposed method \ref{iterated_graphLaNet}.}
    \label{tab:errorsMarmousi}
\end{table}

\begin{table}
    \centering    
\begin{tabular}{|c|c|c|c|c|c|c|}
    \hline
    \multicolumn{2}{|r|}{PSNR} & noiseless & $39$ & $33$ & $30$ & $27$ \\
    \hhline{=======}
    \multirow{4}{*}{\rotatebox[origin=c]{90}{initialization}}
    & AA &
    $\begin{matrix}0.001773\\0.91001\end{matrix}$ & 
    $\begin{matrix}0.0019041\\0.83392\end{matrix}$ & 
    $\begin{matrix}0.0026555\\0.68195\end{matrix}$ & 
    $\begin{matrix}0.0046934\\0.51769\end{matrix}$ & 
    $\begin{matrix}0.0057634\\0.47786\end{matrix}$ \\
    \cline{2-7}
    & Liu &
    $\begin{matrix}0.0016381\\0.92105\end{matrix}$ & 
    $\begin{matrix}0.0018241\\0.89907\end{matrix}$ & 
    $\begin{matrix}0.0051705\\0.78846\end{matrix}$ & 
    $\begin{matrix}0.0090094\\0.67537\end{matrix}$ & 
    $\begin{matrix}0.010531\\0.6205\end{matrix}$ \\
    \cline{2-7}
    & SSI &
    $\begin{matrix}0.0021355\\0.69022\end{matrix}$ & 
    $\begin{matrix}0.0022409\\0.62954\end{matrix}$ & 
    $\begin{matrix}0.002541\\0.5268\end{matrix}$ & 
    $\begin{matrix}0.0029858\\0.42829\end{matrix}$ & 
    $\begin{matrix}0.0035276\\0.3424\end{matrix}$ \\
    \cline{2-7}
    & SB &
    $\begin{matrix}0.014901\\0.54052\end{matrix}$ & 
    $\begin{matrix}0.014905\\0.53939\end{matrix}$ & 
    $\begin{matrix}0.014915\\0.53592\end{matrix}$ & 
    $\begin{matrix}0.018151\\0.49425\end{matrix}$ & 
    $\begin{matrix}0.015148\\0.51241\end{matrix}$ \\
    \hhline{=======}
    \multirow{4}{*}{\rotatebox[origin=c]{90}{after 10 iterations}}
    & AA &
    $\begin{matrix}0.0019756\\0.92203\end{matrix}$ & 
    $\begin{matrix}0.0017359\\0.91855\end{matrix}$ & 
    $\begin{matrix}0.0015104\\0.90535\end{matrix}$ & 
    $\begin{matrix}0.0018774\\0.89058\end{matrix}$ & 
    $\begin{matrix}0.0018206\\0.88341\end{matrix}$ \\
    \cline{2-7}
    & Liu &
    $\begin{matrix}0.0015805\\0.92309\end{matrix}$ & 
    $\begin{matrix}0.0018461\\0.92064\end{matrix}$ & 
    $\begin{matrix}0.0019006\\0.91525\end{matrix}$ & 
    $\begin{matrix}0.0014853\\0.90707\end{matrix}$ & 
    $\begin{matrix}0.0014271\\0.89793\end{matrix}$ \\
    \cline{2-7}
    & SSI &
    $\begin{matrix}0.0017929\\0.78801\end{matrix}$ & 
    $\begin{matrix}0.0018969\\0.78533\end{matrix}$ & 
    $\begin{matrix}0.002096\\0.77983\end{matrix}$ & 
    $\begin{matrix}0.00254\\0.76899\end{matrix}$ & 
    $\begin{matrix}0.0030692\\0.75412\end{matrix}$ \\
    \cline{2-7}
    & SB &
    $\begin{matrix}0.0034704\\0.60788\end{matrix}$ &
    $\begin{matrix}0.0034654\\0.60777\end{matrix}$ & 
    $\begin{matrix}0.0034795\\0.60628\end{matrix}$ & 
    $\begin{matrix}0.0035749\\0.57959\end{matrix}$ & 
    $\begin{matrix}0.0035743\\0.58702\end{matrix}$ \\
    \hline
\end{tabular}
    \caption{D-MSE (top value) and SSIM (bottom value) for reconstructions of the SEAM model with different levels of noise. The top half of the table shows the error values for different initialization methods and the lower half shows the values after 10 iterations of the proposed method \ref{iterated_graphLaNet}.}
    \label{tab:errorsSEAM}
\end{table}

\subsection{Volve field data}

\begin{figure}[h!tb]
	\centering
    \subfloat[]{\includegraphics[height=0.25\textwidth]{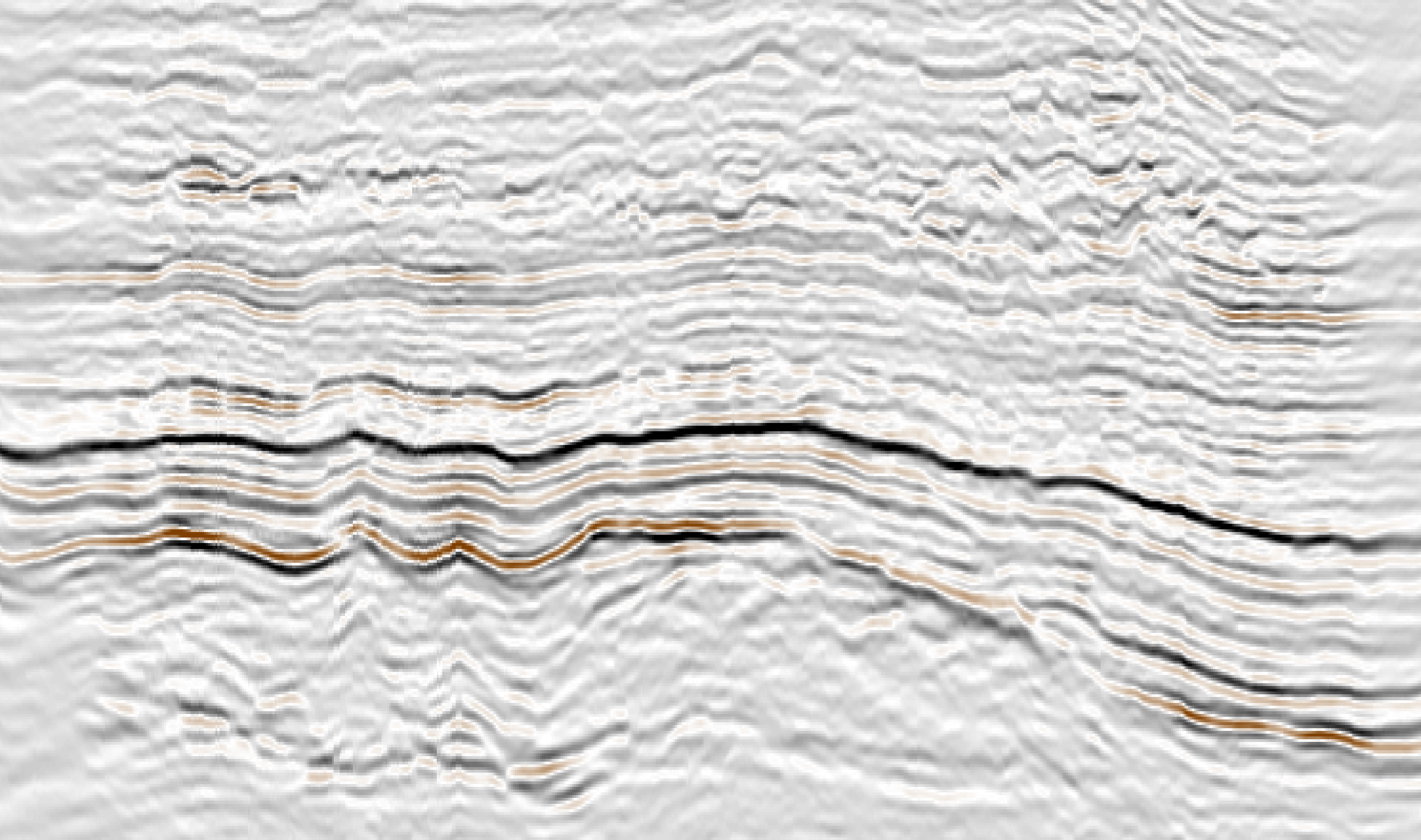}}\ \ 
    \colorbar{seis}{-12.1274}{10.0445}{0.25\textwidth}\ \ 
    \subfloat[]{
    \begin{tikzpicture}
		\begin{axis}[enlargelimits=false,
            x tick label style = {font = \footnotesize},
            y tick label style = {font = \footnotesize},
            xlabel = {seconds},
			axis on top,
            width=0.4\textwidth, height=0.3\textwidth]
			\addplot[blue,line width=2pt] table [x expr=\thisrowno{0}, y expr=\thisrowno{1}, col sep=comma] {Figures/volve/wavelet.csv};
		\end{axis}
	\end{tikzpicture}
    }
	\caption{Volve field seismic data (a) and seismic wavelet extracted from frequency pattern (b).}
    \label{fig:volve_data}
\end{figure}

In our next experiment we apply the algorithm to a post-stack data from the Volve oil field \citep{volvedata}. We use a 2D section of the ST10010ZC11 survey displayed in Figure \ref{fig:volve_data} (a). The data consists of $745$ traces with $850$ time samples each where we picked a $300$ samples wide region of interest that corresponds to a two way travel time of $1.8$ to $3$ seconds. This area contains a well log sample (NO/15-9 19 BT2) which we use for comparison of the different reconstruction methods. However, since one well log is not sufficient to train the discussed neural networks, we use a classical method as the initialization step. Namely, we use SB at it has been superior compared to SSI in the last experiment. This introduces two problems.

First, we have to construct the system matrix $K$ for which we require the seismic wavelet. This can be estimated from the given data. A seismic trace is essentially the sum of scaled and shifted versions of the wavelet and scaling/shifting does not change the frequency distribution. Thus, the frequency distribution of the seismic data is a good approximation of the frequency distribution of the original wavelet. We take the mean absolute value of the Fourier transform over all traces, apply an inverse Fourier transform and use the real part of the obtained signal. This gives the estimated wavelet shown in Figure \ref{fig:volve_data} (b), where we normalized its largest maginute to $10$.

\begin{figure}[h!tb]
	\centering
    \begin{tabular}{lll}
    \subfloat[]{\includegraphics[height=0.25\textwidth]{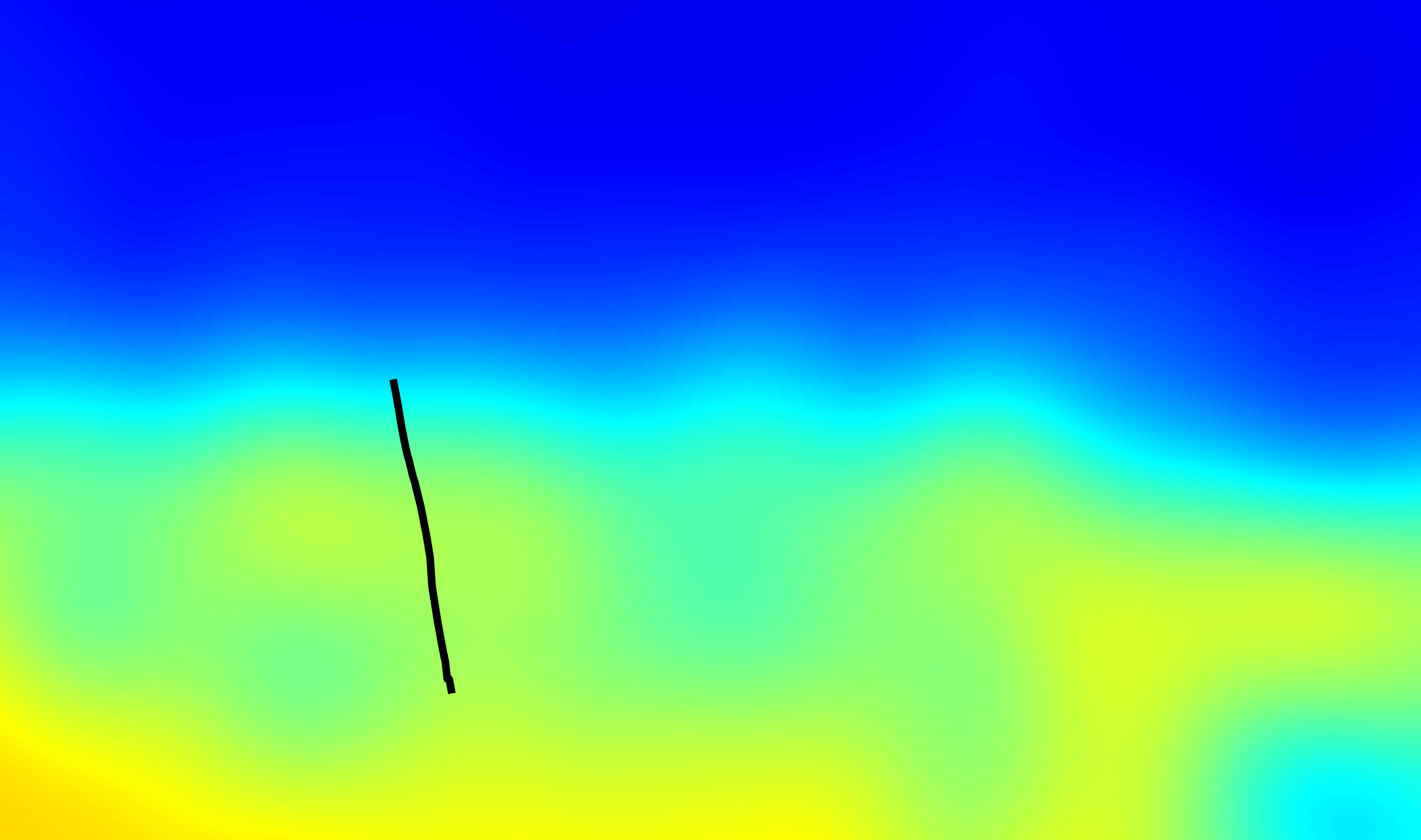}} &
    \subfloat[]{\includegraphics[height=0.25\textwidth]{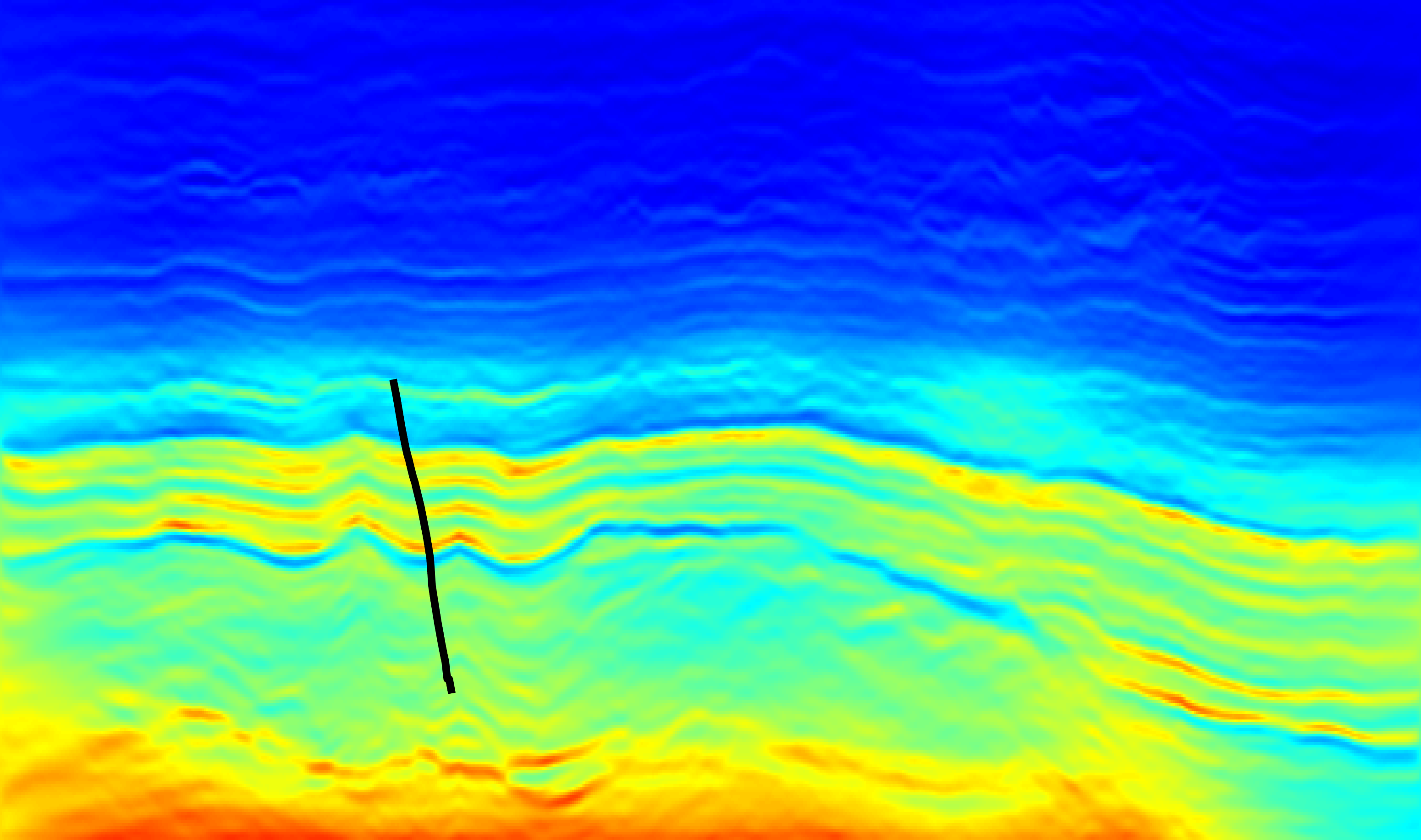}} &
    \colorbar{jet}{2000}{16000}{0.25\textwidth} \\
    \subfloat[]{\includegraphics[height=0.25\textwidth]{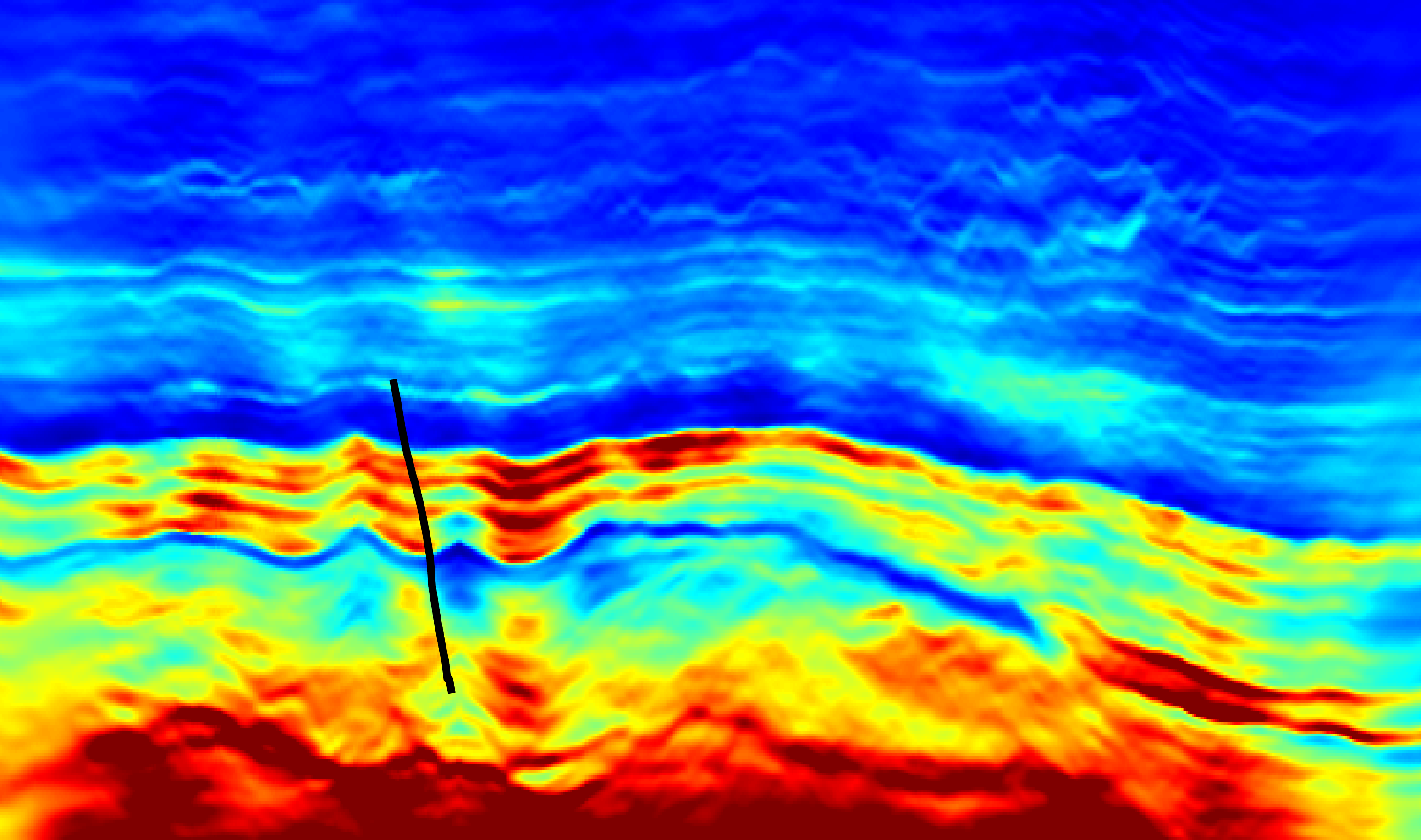}} &
    \subfloat[]{\includegraphics[height=0.25\textwidth]{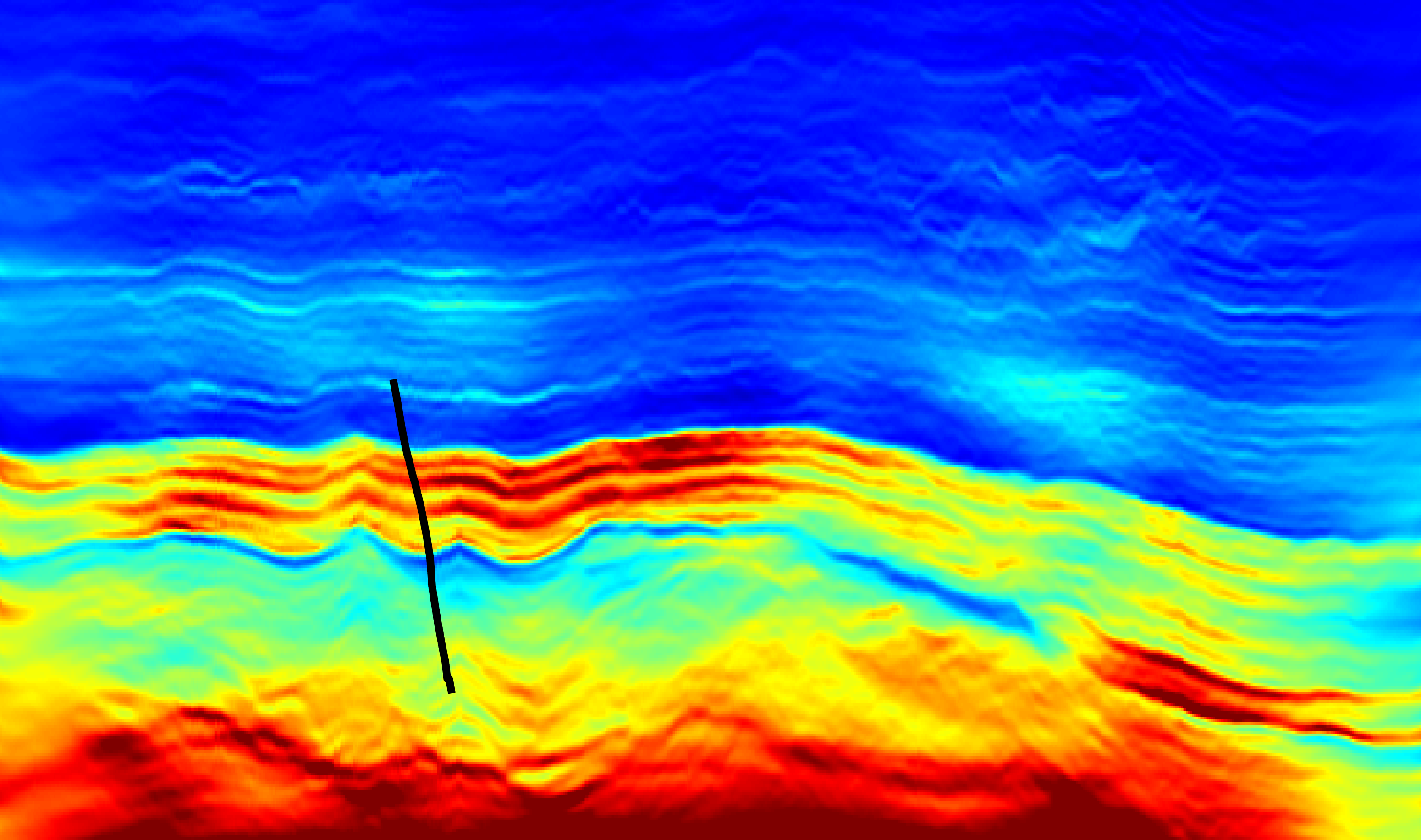}}
    \end{tabular}
	\caption{Reconstructed impedance profiles for Volve field data: background impedance from root mean square velocities (a), initialization using SB (b), \ref{iterated_graphLaNet} approach with $R=3$ (c) and $R=7$ (d). The black line marks the well log sampled data.}
    \label{fig:volve_rec}
\end{figure}

As a second problem, the SB approach assumes sparse differences $\nabla\bx$ and $\nabla\bx^T$. While this is a reasonable assumption for the simplified models seen in the previous experiment, field data is often more complex. Here, the impedance can change continuously within single layers. This leads to the model $\bx=\bx_{\textbf{jump}}+\bx_{\text{cont}}$, i.e., the impedance can be split in two parts. The first part $\bx_{\textbf{jump}}$ is mostly constant only having jump discontinuities at the layer boundaries, which means this term has sparse differences. The second part $\bx_{\text{cont}}$ represents the continuous change throughout the layers, it is also referred to as the background impedance. Fortunately, the background impedance can be approximated very easily. We follow the approach given in \cite{Ravasi2021} where the given root-mean-square velocities are converted into interval velocities and further calibrated with the given well log samples to obtain the background impedance. The result is the convolved with a Gaussian kernel to created a smoothed version. The final result is shown in Figure \ref{fig:volve_rec} (a).

Given the linear operator $K$ and the background impedance $\bx_{\text{cont}}$ we can now solve for $\bx_{\textbf{jump}}$ by minimizing

\begin{align}
    \min\limits_{\bx}\|(\by^\delta-K\bx_{\text{cont}})-K\bx\|_F^2+\alpha\|\nabla\bx\|_{1,1}+\beta\|\nabla\bx^T\|_{1,1}.
\end{align}

The parameter values used in the following calculations are fine tuned by visual comparing the reconstruction to the well log data as well as the taking mean absolute difference as indicator. However, since this data is corrupted by noise and only covers a small part of the area, the used parameter values are not necessarily optimal.

\begin{figure}[h!tb]
	\centering
    \includegraphics[width=\textwidth]{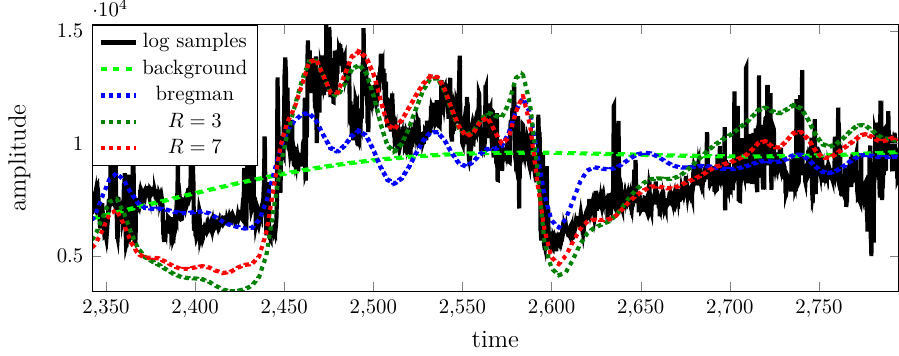}
	\caption{Reconstruction comparison to well log data: sampled well log (black), background impedance (green, dashed), SB reconstruction (blue, dotted), \ref{iterated_graphLaNet} reconstruction with $R=3$ (dark green, dotted) and $R=7$ (red, dotted).}
    \label{fig:volve_well}
\end{figure}

The initial reconstruction is obtained using SB with $\alpha=\beta=200$. Afterwards, we used the proposed method in two different setups. First, we use $R=3$, $\sigma=0.25$ and $10$ iterations which avoids overregularization but is more prone to noise. Second, we use $R=7$, $\sigma=1$ and $5$ iterations which is more stable under the given noise but can overregularize if iterated too long. The obtained reconstructions are shown in Figure~\ref{fig:volve_data}. While the initial reconstruction with SB is able to detect some of the major layer boundaries it still lacks a lot of details. Using the \ref{iterated_graphLaNet} algorithm we are able to improve the result significantly. For $R=3$ the reconstruction contains slightly stronger oscillations due to the noise effect which leads to overshooting or undershooting of the impedance values. This can be seen in more detail when comparing all reconstructions to the well log data (see Figure \ref{fig:volve_well}). The according location of the well log is indicated by the black line in Figure \ref{fig:volve_data}.

\color{black}

\section{Conclusion}\label{sec:conclusions}
In this work we introduced \ref{iterated_graphLaNet}, an iterative version of a Tikhonov-like regularization method for solving an ill-posed impedance inversion problem. \red{The regularizer is based on a data adaptive graph operator constructed from an initial guess of the impedance, which can be obtained using any known inversion technique.} We showcased different numerical experiments with different level of noise and \red{various initialization techniques from simple sparse spike inversion to modern DNN based solvers. Furthermore, we applied our method to real data from the Volve oil field.}

Despite the instabilities of the initialization methods in presence of noise, the iterated graph Laplacian is able to provide stable final reconstructions of much higher quality after a few iterations. The influence of noise is greatly reduced and more details of the impedance profile are recovered. The method maintains stability under many scenarios without the need for significant parameter tuning. 

\red{Limitations: Although the theory of \texttt{graphLa}$\Psi$ is well-developed, rigorous proofs for the regularization properties of the \ref{iterated_graphLaNet} method are still needed.}

\red{Future directions: To improve results even further, the linear forward operator can be replaced with the more accurate non-linear version. However, this requires a more extensive theory and a sophisticated numerical algorithm. Another interesting question is if the proposed method can be used to reduce the number of required data traces and well log samples which are hard to obtain in practice.}


\bibliographystyle{seg}
\bibliography{seismic_impedance}

\end{document}